\input ./style/arxiv-general.cfg
\documentclass[aap,MSNbibl,seceqn,dvips]{arximspdf}
\makeatletter
   \@ifpackageloaded{graphicx}{}{\usepackage{graphicx}}
\makeatother
\usepackage{mathrsfs}

\doi{10.1214/15-AAP1114}
\volume{26}
\issue{2}
\pubyear{2016}
\firstpage{1147}
\lastpage{1207}
\docsubty{FLA}

\makeatletter
\newcommand{\rrVert}{\Vert}
\newcommand{\llVert}{\Vert}
\def\cal{\mathcal}

\newproclaim{Def}{Definition}[section]

\newtheorem{theorem}{Theorem}[section]
\newtheorem{prop}{Proposition}[section]
\newtheorem{cor}{Corollary}[section]
\newtheorem{lemma}{Lemma}[section]
\newproclaim{remark}{Remark}[section]

\newcommand{\ep}{\varepsilon}
\def\mE{\mathbb{E}}

\def\intl{\int}
\def\H{\mathcal H}
\def\wt{\widetilde}
\def\si{{\sigma}}
\def\NN{\mathbb{N}}
\def\EE{\mathbb{E}}

\def\cD{{\cal D}}
\def\cH{{\cal H}}
\def\de{{\delta}}
\def\la{{\lambda}}
\def\si{{\sigma}}
\def\Om{{\Omega}}
\def\al{{\alpha}}
\def\ga{{\gamma}}
\def\de{{\delta}}
\def\si{{\sigma}}
\def\la{{\lambda}}
\def\barU{\overline U}
\def\Om{{\Omega}}
\makeatother

\begin{document}
\begin{frontmatter}

\title{Rate of convergence and asymptotic error distribution of Euler
approximation schemes for fractional diffusions}
\runtitle{Euler schemes of SDE driven by fBm}

\begin{aug}
\author[A]{\fnms{Yaozhong} \snm{Hu}\thanksref{T2}\ead[label=e1]{yhu@ku.edu}},
\author[A]{\fnms{Yanghui} \snm{Liu}\corref{}\ead[label=e3]{y485l911@ku.edu}}
\and
\author[A]{\fnms{David} \snm{Nualart}\thanksref{T3}\ead[label=e2]{nualart@ku.edu}}
\thankstext{T2}{Supported in part by a Grant from the Simons Foundation
\#209206.}
\thankstext{T3}{Supported by  NSF Grant DMS-12-08625 and  ARO Grant FED0070445.}
\runauthor{Y. Hu, Y. Liu and D. Nualart}
\runauthor{Y. Hu, Y. Liu and D. Nualart}
\affiliation{University of Kansas}

\address[A]{Department of Mathematics\\
University of Kansas\\
Lawrence, Kansas 66045\\
USA\\
\printead{e1}\\
\phantom{E-mail:\ }\printead*{e3}\\
\phantom{E-mail:\ }\printead*{e2}}
\end{aug}

%
\received{\smonth{8} \syear{2014}}

%
\begin{abstract}
For a stochastic differential equation(SDE) driven by a fractional
Brownian motion(fBm)
with Hurst parameter $H> \frac{1}2$, it is known
that the existing (naive) Euler scheme has the rate of convergence
$n^{1-2H}$. Since the limit $H \rightarrow\frac{1}2$ of the SDE
corresponds to a Stratonovich SDE driven by standard Brownian motion,
and the naive Euler scheme is the extension of the classical Euler
scheme for It\^o SDEs for $H=\frac{1}2$, the convergence rate of the
naive Euler scheme deteriorates for $H \rightarrow\frac{1}2$.
In this paper we introduce a new (modified Euler) approximation scheme
which is closer to the classical
Euler scheme for Stratonovich SDEs for $H=\frac{1}2$, and it has the
rate of convergence $\gamma_n^{-1}$, where $ \gamma_n=n^{ 2H-{1}/2}$
when $H < \frac{3}4$, $ \gamma_n= n/ \sqrt{ \log
n } $ when $H = \frac{3}4$ and $\gamma_n=n$ if $H> \frac{3}4$.
Furthermore, we study the asymptotic behavior of the fluctuations of
the error. More precisely, if $\{X_t, 0\le t\le T\}$ is the solution of a
SDE
driven by a fBm and
if $\{X_t^n, 0\le t\le T\}$ is its approximation obtained by the new
modified Euler scheme, then we prove that $ \gamma_n
(X^n-X)$ converges stably to the solution of a linear SDE driven by a
matrix-valued Brownian motion, when $H\in( \frac{1}2, \frac{3}4]$.
In the case $H > \frac{3}4$, we show the
$L^p$ convergence of $n(X^n_t-X_t)$, and the limiting process is
identified as the solution of a linear SDE driven
by a matrix-valued Rosenblatt process. The rate of weak convergence is
also deduced for this scheme.
We also apply our approach to the naive Euler scheme.
\end{abstract}

%
\begin{keyword}[class=AMS]
\kwd[Primary ]{60H10}
\kwd[; secondary ]{60H07}
\kwd{26A33}
\kwd{60H35}
\end{keyword}

\begin{keyword}
\kwd{Fractional Brownian motion}
\kwd{stochastic differential equations}
\kwd{Euler scheme}
\kwd{fractional calculus}
\kwd{Malliavin calculus}
\kwd{fourth moment theorem}
\end{keyword}
%
\end{frontmatter}

\section{Introduction}\label{section settings}
Consider the following stochastic differential equation (SDE) on
$\mathbb{R}^d$:
%
%
\begin{equation}
X_t = x + \int_0^t
b(X_s)\,ds +\sum_{j=1}^m \int
^t_0 \sigma^j(X_s)\,dB_s^j,\qquad
t \in[0, T], \label{e.1.1}
\end{equation}
where $x \in\mathbb{R}^d$, $B=(B^1, \ldots, B^m)$ is an $m$-dimensional
fractional Brownian motion (fBm) with Hurst parameter $H\in(\frac{1}2,
1)$ and $b,\sigma^1, \ldots, \sigma^m \dvtx\mathbb{R}^d
\rightarrow
\mathbb{R}^{d }$ are continuous functions.
The above stochastic integrals are pathwise Riemann--Stieltjes integrals.\vadjust{\goodbreak}
If $\si^1, \ldots, \si^m$ are continuously differentiable
and their partial derivatives are bounded and locally H\"older
continuous of
order $\delta> \frac{1}{H} -1$ and $b$ is Lipschitz, then equation
(\ref{e.1.1}) has
a unique solution which is H\"older continuous of order
$\gamma$ for any $0<\gamma<H$. This result was first proved by Lyons
\cite{Ly}, using Young integrals (see \cite{Yo}) and $p$-variation
estimates, and later by Nualart and Rascanu \cite{NuRa}, using
fractional calculus; see \cite{Za}.

We are interested in numerical approximations for the solution to
equation (\ref{e.1.1}).
For simplicity of presentation, we consider uniform partitions of the
interval $[0,T]$,
$t_i = \frac{i T}{n}$, $i = 0, \ldots, n$. For every positive integer $n$,
we define $\eta(t) = t_i $ when $ t_i \leq t < t_i+\frac{T}n$. The
following naive
Euler numerical approximation scheme has been previously studied:
%
%
\begin{equation}
\label{e.1.2} X^n_t = x +\int_0^t
b\bigl(X_{\eta(s)}^n \bigr)\,ds +\sum
_{j=1} ^m \int^t_0
\sigma^j\bigl(X^n_{\eta(s)}\bigr)\,dB^j_s,\qquad
t\in[0,T].
\end{equation}
This scheme can also be written as
\[
X^n_t = X^n_{t_k } + b
\bigl(X^n_{t_k}\bigr) (t-t_k) + \sum
_{j=1}^m \sigma^j \bigl(X^n_{t_k}
\bigr) \bigl(B_t^j - B_{t_k}^j
\bigr), 
\]
for $t_k\le t\le t_{k+1}$, $ k=0, 1, \ldots,n-1 $ and $ X^n_0=x
$.
It was proved by Mishura \cite{Mishura} that for any real number
$\varepsilon>0$
there exists a random variable $C_\varepsilon$ such that almost surely,
\[
\sup_{0\le t\le T} \bigl|X_t^n-X_t\bigr|
\leq C_\varepsilon n^{1-2H +\varepsilon}.
\]
Moreover, the convergence rate $n^{1-2H}$ is sharp for this scheme, in
the sense that $n^{2H-1} [X_t^n -X_t]$ converges almost surely to a
finite and nonzero limit. This has been proved in the one-dimensional
case by Nourdin and Neuenkirch in \cite{Neuenkirch} using the Doss
representation of the solution; see also Theorem~\ref{t.6.1} below.
Notice that while $H$ tends to $\frac{1}2$, the convergence rate $2H-1$
of the numerical scheme
(\ref{e.1.2}) deteriorates, and so it is not a proper extension of the
Euler--Maruyama scheme for the case $H=\frac{1}2$; see, for example,
\cite
{Hu,Kloeden}. This is not surprising because the limit $H \rightarrow
\frac{1}2$ of the SDE (\ref{e.1.1}) corresponds to a Stratonovich SDE
driven by standard Brownian motion, while the Euler scheme (\ref{e.1.2})
is the extension of the classical Euler scheme for the It\^o SDEs.
It is then natural to ask the following question: Can we find a
numerical scheme that generalizes
the Euler--Maruyama scheme to the fBm case?

In this paper we introduce the following new approximation scheme that
we call a \textit{modified Euler scheme}:
%
%
\begin{eqnarray}
\label{e.1.6} 
X^n_t &=& x + \int
^t_0 b(X_{\eta(s)}) \,ds+\sum
_{j=1}^m \int^t_0
\sigma^j\bigl(X^n_{\eta
(s)}\bigr)\,dB^j_s
\nonumber
\\[-8pt]
\\[-8pt]
\nonumber
&&{} + H\sum_{j=1}^{m}\int
^t_0 \bigl(\nabla\sigma^j
\sigma^j\bigr) \bigl(X^n_{\eta(s)}\bigr) \bigl(s-
\eta(s)\bigr)^{2H-1}\,ds,
\end{eqnarray}
or
\begin{eqnarray*}
X^n_t &= &X^n_{t_k} + b
\bigl(X^n_{t_k}\bigr) (t-t_k) + \sum
_{j=1}^m \sigma^j\bigl(X^n_{t_k}
\bigr) \bigl(B_t^j - B_{t_k}^j\bigr)
\\
&&{} + \frac{1}{2}\sum_{j=1}^{m}
\bigl(\nabla\sigma^j \sigma^j\bigr) \bigl(X^n_{t_k}
\bigr) ( t-t_k )^{2H}, 
\end{eqnarray*}
for any $t\in[t_{k},t_{k+1}]$ and $X^n_0=x$. Here
$\nabla\sigma^j$ denotes the $d\times d$ matrix\break $ ( \frac
{\partial\sigma^{j,i}}{
\partial x_k} )_{1\le i,k \le d}$, and
$ (\nabla\sigma^j \sigma^j)^i= \sum_{k=1}^d \frac{\partial\sigma
^{j,i}} {\partial x_k} \sigma^{j,k}$.

Notice that if we formally set $H=\frac{1}2$ and replace $B$ by a
standard Brownian motion $W$,
this is the classical Euler scheme for the Stratonovich SDE,
\begin{eqnarray*}
X_t&=& x + \int_0^t
b(X_s) \,ds+ \sum_{j=1}^m\int
^t_0 \sigma^j( X_s)
\,dW^j_s
\\
&=& x +\int_0^t b(X_s)\,ds + \sum
_{j=1}^m\int_0^t
\sigma^j(X_s) \delta W^j_s +
\frac{1}2 \int_0^t \sum
_{j=1}^m \bigl(\nabla\sigma^j
\sigma^j\bigr) (X_s) \,ds.
\end{eqnarray*}
In the above and throughout this paper, $d$ denotes the Stratonovich
integral, and
$\de$ denotes the It\^o (or Skorohod) integral.

For our new modified Euler scheme (\ref{e.1.6}) we shall prove the
following estimate:
%
%
\begin{equation}
\label{eqn 35} \sup_{0\le t\le T} \bigl( \EE\bigl| X_t-
X^n_t\bigr|^p \bigr)^{{1}/ p} \le C
\ga^{-1}_n,
\end{equation}
for any $p\ge1$, where
%
%
\begin{equation}
\label{e.gamma} \ga_n= %
\cases{ n^{ 2H-{1}/2 },&\quad  $\mbox{if }
\frac{1}2 <H <\frac{3}4$,\vspace*{2pt}
\cr
\displaystyle\frac{n}{ \sqrt{\log n}}, &\quad $
\mbox{if } H =\frac{3}4$,\vspace*{2pt}
\cr
n, &\quad $\mbox{if }
\frac{3}4<H<1 $. } %
\end{equation}
Note that in (\ref{eqn 35}), if we formally set $H = \frac{1}2$, then the
convergence rate is $n^{-{1}/2}$, which is exactly the convergence
rate of the classical Euler--Maruyama scheme in the Brownian motion
case. This suggests that the modified Euler scheme should be viewed as
an authentic modified version of the Euler--Maruyama scheme (\ref
{e.1.2}). The cutoff of the convergence rate for the Euler scheme has
already been observed in a simpler context in \cite{NTU}. The L\'evy
area corresponds to the simple SDE with $b=0$, $\sigma^1(x,y) = (1,0)$,
$\sigma^2(x,y) = (0,x)$. In particular, one has $\nabla\sigma^j
\sigma
^j = 0$, $j=1,2$ here, that is, no diagonal noise.

The proof of this result combines the techniques of Malliavin calculus
with classical fractional calculus. On the other hand, we make use of
uniform estimates for the moments of all orders of the processes $X$,
$X^n$ and their first and second-order Malliavin derivatives, which can
be obtained using techniques of fractional calculus, following the
approach used, for instance, by Hu and Nualart \cite{HuNu}.
The idea of the proof is to properly decompose the error $X_t - X^n_t$
into a weighted quadratic variation term plus a higher order term, that is,
%
%
\begin{equation}
\label{e1.6} X_t-X^n_t = \sum
_{i,j=1}^m \sum_{k=0}^{\lfloor{nt}/T \rfloor}
f^{i,j}(t_k) \int_{t_k}^{t_{k+1 } }
\int_{t_k}^{s} \delta B^i_u
\delta B^j_s + R_t^n,
\end{equation}
where $\lfloor x \rfloor$ denotes the integer part of a real number
$x$. The weighted quadratic variation term provides the desired rate of
convergence in $L^p$.

To further study this new scheme and compare it to the classical
Brownian motion case, it is natural to ask the following questions: Is
the above rate of convergence (\ref{eqn 35})
exact or not? Namely, does the quantity $\ga_n(X_t- X^n_t)$
have a nonzero limit? If yes, how does one identify the limit, and is
there a similarity to the classical Brownian motion case (see \cite
{Jacod,Kurtz})? In the second part of the paper, we give a complete
answer to these questions.
The weighted variation term in (\ref{e1.6}) is still a key ingredient in
our study of the scheme.
As in the Breuer--Major theorem, there is a different behavior in the cases
$H \in(\frac{1}2, \frac{3}4]$ and $H \in(\frac{3}4, 1)$. If $H \in
(\frac{1}2, \frac{3}4]$, we show that $ \ga_n(X_t- X^n_t) $
converges stably to
the solution of a linear stochastic differential equation driven by a
matrix-valued Brownian motion $W$ independent of~$B$. The main tools in
this case are Malliavin calculus and the fourth moment theorem. We will
also make use of a recent limit theorem in law for weighted sums proved
in \cite{CNP}.
In the case $H \in(\frac{3}4, 1)$, we show the convergence of $ \ga
_n(X_t- X^n_t) $ in $L^p $ to the solution of a linear stochastic
differential equation driven by a matrix-valued Rosenblatt process.
Again we use the technique of Malliavin calculus and the convergence in
$L^p$ of weighted sums, which is obtained applying the approach
introduced in
\cite{CNP}. We refer to \cite{NNT} for a discussion on the asymptotic
behavior of some weighted Hermite variations of one-dimensional fBm,
which are related with the results proved here.

We also consider a weak approximation result for our new numerical
scheme. In this case, the rate is $n^{-1}$ for all values of $H$. More
precisely, we are able to show that $n [\mathbb{E} (f(X_t))
-\mathbb{E}( f(X^n_t)) ]$ converges to a finite nonzero limit
which can be explicitly computed. This extends the result of \cite{TT}
to $H > \frac{1}2$. Let us mention that the techniques of Malliavin
calculus also allow us to provide an alternative and simpler proof of
the fact that the rate of convergence of the numerical scheme (\ref
{e.1.2}) is of the order $n^{1-2H}$, and this rate is optimal,
extending to the multidimensional case the results by Neuenkirch and
Nourdin \cite{Neuenkirch}.

If the driven process is a standard Brownian motion,
similar problems have been studied in \cite{Jacod,Kurtz} and the references
therein. See also \cite{CH} for the precise $L^2$-limit and also for a
discussion on the
``best'' partition.
In the case $\frac{1}4 <H<\frac{1}2$ the SDE (\ref{e.1.1}) can be solved
using the theory of rough paths introduced by Lyons; see~\cite{LQ}.
There are also a number of results on the rate of convergence of
Euler-type numerical schemes in this case; see, for instance, the
paper by Deya, Neuenkirch and Tindel \cite{DNT} for a Milstein-type
scheme without L\'evy area in the case $\frac{1}3 <H<\frac{1}2$, the
paper by Friz and Riedel \cite{FR} for the $N$-step Euler scheme
without involving iterated integrals and the monograph by Friz and
Victoir \cite{FV}.

The paper is organized as follows. The next section contains some basic
materials on fractional
calculus and Malliavin calculus that will be used throughout the paper,
and introduces a matrix-valued
Brownian motion and a generalized Rosenblatt process, both of which are
key ingredients in our results on
the asymptotic behavior of the error; see Section~\ref{eqn 43} and
Section~\ref{section thmm2 proof}.
In Section~\ref{subsec estimate
solution of SDE}, we derive the necessary estimates for the uniform norms
and H\"older seminorms of the
processes $X$, $X^n$ and their Malliavin derivatives. In Section~\ref{section strong conv}, we
prove our result on the rate
of convergence in $L^p$ for the numerical scheme (\ref{e.1.6}).
In Section~\ref{s.theorem1}, we prove a central limit theorem for weighted
quadratic sums, and then
in Section~\ref{eqn 43} we apply this result to the study of
the asymptotic behavior of the error $\gamma_n (X_t - X^n_t)$ in case
$H \in(\frac{1}2, \frac{3}4]$.
In Section~\ref{section 7}, we study the $L^p$-convergence of some
weighted random sums.
In Section~\ref{section thmm2 proof}, we apply the results of
Section~\ref{section 7} to establish the $L^p$-limit of $ n (X_t - X^n_t)$ in case $H
\in
(\frac{3}4, 1)$.
The weak approximation result is discussed in Section~\ref{sec9}. In Section~\ref{sec10},
we deal with the numerical scheme (\ref{e.1.2}).
In the \hyperref[app]{Appendix}, we prove some auxiliary results.

\section{Preliminaries and notation}\label{sec2}

Throughout the paper we consider a fixed time interval $[0, T]$. To
simplify the presentation we only deal
with the uniform partition of this interval; that is, for each $n\ge1$
and $i=0,1,\ldots, n$, we set
$t_i= \frac{iT}n$.
We use $C$ and $K$ to
represent constants that are independent of $n$ and
whose values may change from line to line.

\subsection{Elements of fractional calculus}\label{sec2.1}
In this subsection we introduce the definitions of the fractional
integral and derivative operators, and we review some properties of
these operators.

Let $a, b \in[0,T]$ with $a < b$, and let $\beta\in( 0, 1)$. We
denote by $C^{\beta} (a, b)$ the space of $\beta$-H\"older continuous
functions on the interval $[a, b]$. For a function $x\dvtx[0,T]
\rightarrow
\mathbb{R}$, $\|x\|_{a, b, \beta} $ denotes the $\beta$-H\"older
seminorm of $x$ on $[a, b]$,
that is,
\[
\|x\|_{a,b, \beta} = \sup\biggl\{\frac{|x_u-
x_v|}{(v-u)^{\beta}}; a \leq u < v \leq b
\biggr\}.
\]
We will also make use of the following seminorm:
%
%
\begin{equation}
\label{e2.1} 
\|x\|_{a,b, \beta, n} = \sup\biggl\{\frac{|x_u-
x_v|}{(v-u)^{\beta}};
a \leq u < v \leq b, \eta(u) = u \biggr\}. 
\end{equation}
Recall that for each $n\ge1$ and $i=0,1, \ldots, n$, $t_i = \frac{
iT}{n}$ and $\eta(t) = t_i $ when $ t_i \leq t < t_i+\frac{T}n$.

We will denote the uniform norm of $x$ on the interval $[a, b]$ as $\|
x\|_{a, b, \infty}$. When $a=0$ and $b=T$, we will simply write $\|x\|
_\infty$ for $\|x\|_{0,T, \infty}$ and $\|x\|_\beta$ for $\|x\|_{0,T,
\beta}$.

Let $f \in L^1 ([a, b])$ and $\alpha>0$. The left-sided and
right-sided fractional Riemann--Liouville integrals of $f$ of order
$\alpha$ are defined, for almost all $t \in(a, b)$, by
\[
I^{\alpha}_{a+} f(t) = \frac{1}{\Gamma(\alpha)} \int
^t_a (t -s)^{\alpha-1} f(s) \,ds
\]
and
\[
I^{\alpha}_{b-} f(t) = \frac{(-1)^{-\alpha} }{\Gamma(\alpha)} \int
^b_t (s -t)^{\alpha-1} f(s) \,ds,
\]
respectively, where $(-1)^{\alpha} = e^{-i \pi\alpha}$ and $\Gamma
(\alpha) = \int^{\infty}_0 r^{\alpha-1} e^{-r} \,dr$ is the Gamma
function. Let $I^{\alpha}_{a+} (L^p)$ [resp., $I^{\alpha}_{b-} (L^p)$]
be the image of $L^p ([a, b])$ by the operator $I^{\alpha}_{a+} $
(resp., $I^{\alpha}_{b-} $). If $f \in I^{\alpha}_{a+} (L^p)$ [resp., $f
\in I^{\alpha}_{b-} (L^p)$] and $0 < \alpha<1$, then the fractional
Weyl derivatives are defined as
%
%
\begin{equation}
\label{e.2.1} D^{\alpha}_{a+} f(t) = \frac{1}{\Gamma(1-\alpha)}
\biggl(
\frac
{f(t)}{(t-a)^{\alpha}} + \alpha\int^t_a
\frac{f(t) -
f(s)}{(t-s)^{\alpha+1}} \,ds \biggr)
\end{equation}
and
%
%
\begin{equation}
\label{e.2.2} D^{\alpha}_{b-} f(t) = \frac{(-1)^{\alpha} }{\Gamma
(1-\alpha)} \biggl(
\frac{f(t)}{(b-t)^{\alpha}} + \alpha\int^b_t
\frac{f(t) -
f(s)}{(s-t)^{\alpha+1}} \,ds \biggr),
\end{equation}
where $a <t< b$.

Suppose that $f \in C^{\lambda} (a, b)$ and $g \in C^{\mu} (a, b)$ with
$\lambda+ \mu>1$. Then, according to Young \cite{Yo}, the
Riemann--Stieltjes integral $\int^b_a f\,dg$ exists. The following
proposition can be regarded as a fractional integration by parts
formula, and provides an explicit expression for the integral $\int^b_a
f\,dg$ in terms of fractional derivatives. We refer to \cite{Za} for
additional details.

%
\begin{prop}\label{prop2.1}
Suppose that $f \in C^{\lambda}(a, b)$ and $g\in C^{\mu} (a, b)$ with
$\lambda+\break  \mu>1$. Let $\lambda> \alpha$ and $\mu> 1- \alpha$. Then
the Riemann--Stieltjes integral $\int^b_a f \,dg$ exists, and it can be
expressed as
%
%
\begin{equation}
\label{fracibp} \int^b_a f\,dg =
(-1)^{\alpha} \int^b_a
D^{\alpha}_{a+} f(t) D^{1-\alpha
}_{b-}
g_{b-}(t) \,dt,
\end{equation}
where $g_{b-} (t) =\mathbf{1}_{(a, b)}(t) ( g(t) - g(b-) )$.
\end{prop}

The notion of H\"older continuity and the above result on the existence
of Riemann--Stieltjes integrals can be generalized to functions taking
values in some normed spaces.
We fix a probability space $(\Omega, \mathscr{F}, P)$ and denote by
$\|
\cdot\|_p $ the norm in the space $L^p:=L^p(\Omega)$, where $p\ge1$.

%
\begin{Def}
Let $f=\{ f(t), t \in[0, T]\}$ be a stochastic process such that $f(t)
\in L^p$ for all $t \in[0, T]$. We say that $f$ is H\"older continuous
of order $\beta>0$ in $L^p$ if
%
%
\begin{equation}
\label{eqn2.4} \bigl\|f(t)- f(s)\bigr \|_{p} \le C |t-s|^{ \beta},
\end{equation}
for all $s,t\in[0,T]$.
\end{Def}

The following result shows that with proper H\"older continuity
assumptions on $f$ and $g$, the
Riemann--Stieltjes integral $\int_0^T f\,dg$ exists, and equation
(\ref{fracibp}) holds.

%
\begin{prop}\label{prop.1}
Let the positive numbers
$p_0$, $\lambda$, $\mu$, $p$, $q$ satisfy $p_0 \geq1$, $\lambda+\mu>1$,
$\frac{1}{p}+\frac{1}{q}=1$ and $p_0p> \frac{1}{\mu}$, $p_0q> \frac
{1}{\lambda}$.
Assume that $f=\{f(t), t\in[0,T\} $ and $g=\{g(t), t\in[0, T]\}$ are
H\"older continuous stochastic processes of order $\mu$ and $\lambda$
in $L^{p_0 p} $ and $L^{p_0 q} $, respectively, and $f(0) \in L^{p_0 p} $.
Let $\pi: 0=t_0 < t_1 < \cdots< t_N =T$ be a partition on $[0, T]$,
and $\xi_i\dvtx t_{i-1} \leq\xi_i \leq t_i $.
Then the sum $\sum_{i=1}^N f(\xi_i) [g(t_i) - g(t_{i-1})]$ converges in
$L^{p_0}$ to the Riemann--Stieltjes integral $\int_0^T f \,dg$ as
$|\pi|$
tends to zero, where $|\pi| = \max_{1 \leq i \leq N} | t_i - t_{i-1} |
$, and equation (\ref{fracibp}) holds.
\end{prop}

Proposition~\ref{prop.1} can be proved through a slight modification
of Z\"ahle's proof in the real-valued case
\cite{Za} using H\"older's inequality.

\subsection{Elements of Malliavin calculus}\label{sec2.2}
We briefly recall some basic facts about the stochastic calculus of
variations with respect to an fBm. We refer the reader to~\cite
{Nualart2} for further details. Let $B = \{(B^1_t, \ldots, B^m_t ), t
\in[0, T] \}$ be an $m$-dimensional
fBm with Hurst parameter $H\in(\frac{1}2, 1)$, defined on some
complete probability space $(\Omega, \mathscr{F}, {P})$. Namely, $B$ is
a mean zero Gaussian process with covariance
\[
\mE\bigl(B^i_t B^j_s\bigr) =
\tfrac{1}2 \bigl(t^{2H} + s^{2H} - |t-s|^{2H}
\bigr) \delta_{ij}, \qquad i,j=1,\ldots, m,
\]
for all $s, t \in[0, T]$, where $\delta_{ij}$ is the Kronecker symbol.

Let
$\mathcal{H}$ be the Hilbert space defined as the closure of the set
of step functions on $[0, T]$ with respect to the scalar product
\[
\langle\mathbf{1}_{[0, t]}, \mathbf{1}_{[0, s ]}
\rangle_{\mathcal{H} } = \tfrac{1}2 \bigl(t^{2H} +
s^{2H} - |t-s|^{2H} \bigr).
\]
It is easy to see that the covariance of fBm can be written as
\[
\alpha_H \int_0^t\int
_0^s |u-v|^{2H-2}\,du\,dv,
\]
where $\alpha_H = H(2H-1)$. This implies that
\[
\langle\psi, \phi\rangle_{\cH} = \alpha_H \int
_0^T\int_0^T
\psi_u \phi_v |u-v|^{2H-2}\,du\,dv
\]
for any pair of step functions $\phi$ and $\psi$ on $[0, T]$.

The elements of the Hilbert space $\cH$, or more generally, of the
space $\cH^{\otimes l }$
may not be functions, but distributions; see \cite{PT1} and \cite{PT2}.
We can find a linear space of functions contained in $\cH^{\otimes l} $
in the following way:
Let $|\cH|^{\otimes l}$ be the linear space of measurable functions
$\phi$ on $[0, T]^l \subset\mathbb{R}^l$ such that
\[
\|\phi\|^2_{|\cH|^{\otimes l} }:= \alpha_H^l
\int_{[0,T]^{2l} } |\phi_{\mathbf{u}}| |\phi_{\mathbf{v}}|
|u_1-v_1|^{2H-2} \cdots|u_l-v_l|^{2H-2}
\,d{\mathbf{u}}\,d{\mathbf{v}} < \infty,
\]
where $ {\mathbf{u}}= (u_1, \ldots, u_l), {\mathbf{v}} = (v_1,
\ldots,
v_l) \in[0, T]^l$.
Suppose $\phi\in L^{{1}/{H} } ([0, T]^l ) $. The following
estimate holds:
%
%
\begin{equation}
\label{eq2.5} \|\phi\|_{ {|\cH|^{\otimes l} } } \leq b_{H,l} \|\phi
\|_{L^{{1}/H}
([0, T]^l)}
\end{equation}
for some constant $b_{H,l}>0$; the case $l=1$ was proved
in \cite{MMV}, and the extension to general case
is easy; see \cite{hunu09}, equation (2.5).

The
mapping $\mathbf{1}_{[0, t_1]} \times\cdots\times\mathbf{1}_{[0, t_m]}
\mapsto( B_{t_1}^1, \ldots, B_{t_m}^m) $ can be extended to a
linear isometry between $\mathcal{H}^m $ and the Gaussian space spanned
by $B$. We denote this isometry by $h \mapsto B(h)$. In this way,
$\{B(h), h \in\mathcal{H}^m\}$ is an isonormal Gaussian process
indexed by the Hilbert space $\mathcal{H}^m$.

Let $\mathcal{S}$ be the set of smooth and cylindrical random
variables of
the form
\[
F=f ( B_{s_1},\ldots,B_{s_N} ),
\]
where $N\geq1$ and $f\in C_{b}^{\infty} ( \mathbb{R}^{m\times
N} ) $.
For each $j=1,\ldots, m$ and $t\in[0,T]$, the derivative operator
$D^jF$ on $F\in\mathcal{S}$ is defined as the $\mathcal{H} $-valued
random variable
\[
D_t^jF =\sum_{i=1}^{N}
\frac{\partial f}{\partial x_{i}^j} ( B_{s_1} ,\ldots,B_{s_N} )
\mathbf{1}_{[0,s_i]}(t),\qquad t\in[0, T].
\]
We can iterate this procedure to define higher order derivatives
$ D^{j_1, \ldots,j_l} F$ which take values on $ \cH^{\otimes l}$.
For any $p\geq
1$ and any integer $k\ge1$, we define the Sobolev space $\mathbb
{D}^{k,p}$ as the closure of $%
\mathcal{S}$ with respect to the norm
\[
\llVert F\rrVert_{k,p}^{p}=\mathbb{E} \bigl[ \vert F
\vert^{p} \bigr] + \mathbb{E} \Biggl[ \sum
_{l=1}^{k} \Biggl(\sum_{j_1, \ldots, j_l=1}^m
\bigl\llVert D^{j_1, \ldots, j_l} F\bigr\rrVert^2 _{\mathcal
{H}^{\otimes l}}
\Biggr)^{{p}/2} \Biggr].
\]
If $V$ is a Hilbert space, $\mathbb{D}^{k,p}(V)$ denotes the
corresponding Sobolev space of $V$-valued random variables.

For any $j=1,\ldots, m$, we denote by $\delta^j$ the adjoint of the
derivative operator $D^j$. We say $u\in\operatorname{Dom} \delta^j$ if
there is a $ \delta^j(u) \in L^2(\Om)$ such that for any $F\in
\mathbb{D}^{1,2}$ the following duality relationship holds:
%
%
\begin{equation}
\label{ibp} \mathbb{E} \bigl( \bigl\langle u,D^jF \bigr
\rangle_{\H} \bigr)= \mathbb{E} \bigl( \delta^j(u) F\bigr).
\end{equation}
The random variable $\delta^j(u)$ is also called the Skorohod integral
of $u$ with
respect to the fBm $B^j$, and we use the notation
$ \delta^j(u)= \int^T_0 u_t \delta B^{j}_t$.

Let $F \in\mathbb{D}^{1, 2}$ and $u$ be in the domain of $\delta^j$
such that $Fu \in L^2 ( \Omega; \cH)$. Then (see~\cite{Nualart}) $Fu$
belongs to the domain of $\delta^j$, and the following equality holds:
%
%
\begin{equation}
\label{e1} \delta^j (Fu) = F \delta^j (u) - \bigl
\langle D^j F, u \bigr\rangle_{\mathcal{H
} },
\end{equation}
provided the right-hand side of (\ref{e1}) is square integrable.

Suppose that $u=\{u_t, t\in[0,T] \} $ is a stochastic process
whose trajectories are H\"older continuous of order $\gamma>1-H$.
Then, for any $j=1,\ldots, m$, the Riemann--Stieltjes
integral $ \int^T_0 u_t \,d B^{j}_t$ exists. On the other hand, if
$u\in\mathbb{D}^{ 1,2}(\mathcal{H})$ and the derivative $D^j_su_t$
exists and satisfies almost surely
\[
\intl^T_0 \intl^T_0
\bigl|D^j_s u_t\bigr| |t -s|^{2H-2} \,ds\,dt <
\infty,
\]
and $ \mathbb{E} ( \|D^j u\|^2_{L^{{1}/H}([0,T]^2)}
)<\infty$, then (see
Proposition~5.2.3 in \cite{Nualart}) $\intl^T_0 u_t \delta B^{j}_t$
exists, and we
have the following relationship between these two stochastic
integrals:
%
%
\begin{equation}
\label{e.2.3} \intl^T_0 u_t\, d
B^{j}_t = \intl^T_0
u_t \delta B^{j}_t + \alpha_H
\intl_0^T \intl^T_0
D^j_s u_t |t-s|^{2H-2} \,ds\,dt.
\end{equation}

The following result is Meyer's inequality for the Skorohod integral;
see, for example, Proposition~1.5.7 of \cite{Nualart}. Given $p>1$ and
an integer $k \geq1$, there
is a
constant $c_{k, p}$ such that
%
%
\begin{equation}
\label{eqn meyer}\bigl \| \delta^k(u)\bigr \|_{p} \leq
c_{k,p} \|u\| _{\mathbb{D}^{k, p}
(\mathcal
{H}^{\otimes k}
)} \qquad\mbox{for all $u \in\mathbb{D}^{k, p}
\bigl(\mathcal{H} ^{ \otimes
k}\bigr)$}.
\end{equation}
Applying (\ref{eq2.5}) and then the Minkowski inequality to the
right-hand side of (\ref{eqn meyer}) yields
%
%
\begin{eqnarray}
\label{eqn 2.10} \bigl\| \delta^k(u) \bigr\|_{p} &\leq& C \bigl
\llVert\llVert u\rrVert_p \bigr\rrVert_{L^{{1}/H} ([0, T]^{p })}
\nonumber
\\[-8pt]
\\[-8pt]
\nonumber
&&{} +
C\sum
_{l=1}^{k} \sum
_{j_1, \ldots, j_l=1}^m \bigl\llVert\bigl\llVert D^{j_1, \ldots, j_l}
u \bigr\rrVert_p\bigr\rrVert_{L^{{1}/H}
([0, T]^{p+l})}
\end{eqnarray}
{for all $u \in\mathbb{D}^{k, p} (\mathcal{H} ^{ \otimes
k})$}, provided $pH \geq1$.

\subsection{Stable convergence}
Let $Y_n$, $n\in\mathbb{N}$ be a sequence of random variables defined
on a probability space
$(\Omega, \mathscr{F}, {P})$
with values in a Polish space $(E, \mathscr{E})$. We say that $Y_n$
\textit{converges stably} to
the limit $Y$, where $Y$ is defined on an extension of the original
probability space
$(\Omega', \mathscr{F}', {P}')$, if and only if for any bounded
$\mathscr{F}$-measurable random variable $Z$, it holds that
\[
(Y_n, Z) \Rightarrow(Y, Z)
\]
as $n \rightarrow\infty$, where $\Rightarrow$ denotes the convergence
in law.

Note that stable convergence is stronger than weak convergence but
weaker than convergence in probability. We refer to \cite
{JacodShiryaev} and \cite{Aldous} for more details on this concept.

\subsection{A matrix-valued Brownian motion}\label{sec2.3}

The aim of this subsection is to define a matrix-valued Brownian motion
that will play a
fundamental role in our central limit theorem. First, we introduce two
constants $Q$ and $R$ which depend on $H$.

Denote by $\mu$ the measure on $\mathbb{R}^2 $ with density
$|s-t|^{2H-2}$. Define, for each $p \in\mathbb{Z} $,
\[
Q( p) =T^{4H} \intl_{ {0}}^{ { 1}}
\intl_{ { p} }^{ {{ p} +1}} \intl_{ {0}}^{ t }
\intl_{ { p} }^{s} \mu(dv\,du) \mu(ds \,dt)
\]
and
\[
R(p) = T^{4H} \intl_{ {0}}^{ { 1}}
\intl_{ p}^{ { p+1}} \intl_{ {t}}^{ 1 }
\intl_{ p}^{s} \mu(dv \,du) \mu(ds \,dt).
\]
It is not difficult to check that for $\frac{1}2 <H< \frac{3}4$, the
series $ \sum_{p \in\mathbb{Z} } Q(p)$ and
$\sum_{p \in\mathbb{Z} } R(p)$ are convergent, and for $H=\frac{3}4$,
they diverge at the rate $\log n$. Then we set
(we omit the explicit dependence of $Q$ and $R$ on $H$ to simplify the notation)
%
%
\begin{equation}
\label{e2.6} Q = \sum_{p \in\mathbb{Z} } Q(p),\qquad R = \sum
_{p \in\mathbb{Z} } R(p),
\end{equation}
for the case $ H \in(\frac{1}2, \frac{3}4)$,
and
\[
Q= \lim_{n \rightarrow\infty} \frac
{ \sum_{|p| \leq n} Q(p) }{ \log n} = \frac{T^{4H}}{2} ,\qquad R=
\lim_{n \rightarrow\infty} \frac{ \sum_{|p| \leq n} R(p) }{
\log n} = \frac{T^{4H}}{2},
\]
for the case $ H = \frac{3}4$.
%

\begin{lemma}
The constants $Q$ and $R$ satisfy $R\le Q$.
\end{lemma}

\begin{pf}
If $H = \frac{3}4$, we see from (\ref{e2.6}) that these two constants
are both equal to $\frac{T^{ 4H } }{ 2 }$.
Suppose $H \in(\frac{1}2, \frac{3}4)$. Consider the functions on
$\mathbb
{R}^2$ defined by $\varphi_p(v,s)= \mathbf{1}_{\{p\le v\le s \le p+1\}
}$, $\psi_p(v,s)= \mathbf{1}_{\{p\le s\le v \le p+1\}}$, $p \in
\mathbb
{Z}$. Then
\begin{eqnarray*}
&&\frac{1}n \Biggl\llVert\sum_{p=0}^{n-1}
(\varphi_p -\psi_p) \Biggr\rrVert^2_{L^2(\mathbb{R}^2,\mu)}
\\
&&\qquad=\frac{2}n \sum_{p,q=0}^{n-1}
\bigl( \langle\mathbf{1}_{\{
p\le
v\le s\le p+1\}}, \mathbf{1}_{\{ q\le v\le s\le q+1\}} \rangle
_{L^2(\mathbb{R}^2,\mu)}
\\
&&\hspace*{69pt}{}- \langle\mathbf{1}_{\{ p\le v\le s\le p+1\}}, \mathbf{1}_{\{
q\le s\le v\le q+1\}}
\rangle_{L^2(\mathbb{R}^2,\mu)} \bigr)
\\
&&\qquad= \frac{2}n \sum_{p,q=0}^{n-1}
\bigl( Q(p-q) -R(p-q) \bigr).
\end{eqnarray*}
It is easy to see that the above is equal to
\[
\frac{2}n \sum_{j=0}^{n-1} \sum
_{k=-j}^{ j } \bigl( Q(k )- R(k ) \bigr).
\]
It then follows from a Ces\`aro limit argument that the quantity in the
right-hand side of the above converges to $2(Q-R)$ as $n$ tends to
infinity. Therefore, $Q\ge R$.
\end{pf}

Let $ \wt{W}^{0,ij} =\{ \wt{W}^{0,ij}_t, t\in[0,T]\} $, $ i\leq j,
i,j = 1, \ldots, m $ and $
\wt{W}^{1,ij}= \{\wt{W}^{1,ij}_t, t\in[0,T]\}$, $i, j = 1, \ldots,
m $
be independent
standard Brownian motions. When $i > j$, we define
$\wt{W}^{0,ij}_t= \wt{W}^{0,ji }_t $. The matrix-valued Brownian motion
$(W^{ij})_{1\le i,j \le m}$,
$i,j=1, \ldots, m$ is defined as follows:
\[
W^{ii} = \frac{\alpha_H}{\sqrt{T } } \bigl( \sqrt{Q+R} \wt{W}^{1,ii}
\bigr)
\]
and
\[
W^{ij} =\frac{\alpha_H}{\sqrt{T } } \bigl( \sqrt{Q-R} \wt
{W}^{1,ij} +
\sqrt{R} \wt{W}^{0,ij} \bigr) \qquad\mbox{when } i\neq j.
\]
Notice that this definition makes sense because $R\le Q$. The random
matrix $W_t$ is not symmetric when $H<\frac{3}4$; see the plot and table
below. For $i,j,i',j'=1,\ldots, m$, the covariance $\mE(W^{ij}_t
W^{i'j'}_s )$ is equal to
\[
\frac{ \alpha_H^2 (t \wedge s) }{ T} ( R \delta_{ji'} \delta_{i j'}
+ Q
\delta_{j j'} \delta_{ii'} ) ,
\]
where $\delta$ is the Kronecker function.

\begin{figure}

\includegraphics{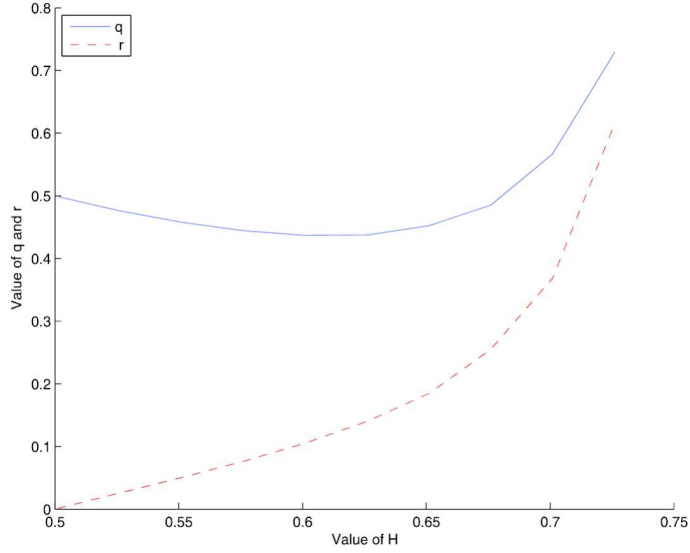}

\caption{Simulation of $q$ and $r$.}\label{fig1}
\end{figure}

In the Figure \ref{fig1} and Table \ref{tab1}, we consider two quantities for $H\in
(\frac{1}2, \frac{3}4)$,
\[
q=\frac{\al_H^2}{T^{4H} }Q\quad \mbox{and}\quad r=\frac{\al_H^2}{T^{4H}
} R.
\]
We see that the values of $q$ and $r$ approach $0.5$ and $0$ as $H$
tends to $\frac{1}2$, respectively, and both of them tend to infinity
when $H$ gets closer to $\frac{3}4$.
%
%

%

\subsection{A matrix-valued generalized Rosenblatt process}\label
{section 2.3}

In this subsection we
introduce a generalized Rosenblatt process which will appear in the
limiting result proved in Section~\ref{section thmm2 proof} when
$H>\frac{3}4$.
Consider an $m$-dimensional fBm
$B_t= (B^1_t,
\ldots, B^m_t)$ with Hurst parameter $H \in(\frac{3}4, 1)$. Define for
$i_1,i_2 \in1, \ldots, m$,
\[
Z^{i_1,i_2}_n(t):= n\sum_{j=1}^{\lfloor{nt}/{T} \rfloor}
\intl_{t_j}^{t_{j+1} }\bigl( B^{i_1}_s -
B^{i_1}_{t_j}\bigr) \delta B^{i_2}_s.
\]
When $i_1=i_2=i$, we can write
\[
Z^{i,i}_n(t)= \frac{T^{2H } }{ 2n^{2H-1 } } \sum
_{j=1}^{\lfloor
{nt}/{T} \rfloor} H_2\bigl({
\xi}^{n,i}_j\bigr),
\]
where $H_2(x)=x^2-1$ is the second degree Hermite polynomial and $ {\xi
}^{n,i}_j = T^{-H} n^H (B^i_{t_{j+1} } - B^i_{t_j} )$.
It is well known (see \cite{NNT}) that for each $i=1,\ldots, m$, the
process $Z^{i,i}_n(t)$ converges in $L^2$ to the \textit{Rosenblatt
process} $R(t)$. We refer the reader to \cite{Rosenblatt} and \cite
{Tudor} for further details on the Rosenblatt process.

%
\begin{table}
\caption{Simulation of $q$ and $r$}\label{tab1}
\begin{tabular*}{\textwidth}{@{\extracolsep{\fill}}lcccccccc@{}}
\hline
$H$ & 0.5010 & 0.5260 & 0.5510 & 0.6010 & 0.6260 & 0.6510 & 0.7010
&0.7260 \\
$q$ & 0.4990 & 0.4763 & 0.4580& 0.4369 & 0.4375 & 0.4522 & 0.5669 &
0.7290 \\
$r$ & 9.9868$\times10^{-4}$ & 0.0256 & 0.0503 & 0.1053 & 0.1400
& 0.1845 & 0.3689 & 0.6149 \\
\hline
\end{tabular*}
%
\end{table}

When $i_1 \neq i_2$, the stochastic integral $\intl_{t_j}^{t_{j+1} }(
B^{i_1}_s -
B^{i_1}_{t_j}) \delta B^{i_2}_s$ cannot be written as the second
Hermite polynomial of a Gaussian random variable.
Nevertheless, the process $Z^{i_1,i_2}_n(t)$ is still\ convergent in
$L^2$. Indeed, for any positive integers $n$ and~$n'$, we have
\begin{eqnarray*}
& &\mE\bigl( Z_n^{i_1 i_2}(t) Z_{n'}^{i_1 i_2}(t)
\bigr)
\\
&&\qquad= nn' \sum_{k=0}^{ \lfloor{nt}/{T} \rfloor} \sum
_{k'=0}^{\lfloor
{ { n'} t}/T\rfloor} \mE\biggl[
\intl_{({k }/{n}) T}^{ ({(k+1)}/{n}) T } \bigl( B^{i_1 }_s-
B^{i_1 }_{ ({k }/{n}) T} \bigr) \delta B_s^{i_2}
\\
&&\hspace*{116pt}{}\times \intl_{({k }/{ n' }) T}^{ ({(k+1)}/{ n' }) T } \bigl( B^{i_1 }_s-
B^{i_1 }_{ ({k }/{n'}) T} \bigr) \delta B_s^{i_2}
\biggr]
\\
&&\qquad= n{ n'} \alpha_H^2 \sum
_{k=0}^{ \lfloor{nt}/{T} \rfloor} \sum_{k'=0}^{\lfloor{ {n' }
t}/T\rfloor}
\intl^{({(k+1)}/{n})T}_{( {k}/n )T} \intl^{({(k'+1)}/{{ {n'}}}) T
}_{({k'} /{ {n}'}) T }
\intl^{t}_{({k}/n) T} \intl^{s}_{({k'} /{ {n}' }) T} \mu(dv
\,du)\\
&&\hspace*{300pt}{}\times \mu(ds \,dt)
\\
&&\qquad\rightarrow\frac{ T^2 \alpha_H^2 }{4} \intl^t_0
\intl^t_0 |u-v|^{4H-4} \,du\,dv
\\
&&\qquad= c_H t^{4H-2},
\end{eqnarray*}
as $ {n}', n \rightarrow+\infty$, where $c_H = \frac{
T^2H^2(2H-1) } {4 (4H-3)} $. This allows us to conclude that $Z_n^{i_1 i_2}(t)$
is a Cauchy sequence in $L^2$.
We denote by $Z^{i_1 i_2}_t $ the $L^2$-limit of
$Z_n^{i_1 i_2}(t)$. Then $Z^{i_1 i_2}_t $ can be considered a
$\mathit{generalized}$ $\mathit{Rosenblatt}$ $\mathit{process}$.

It is easy to show that
\[
\mE\bigl[\bigl|Z^{i_1 i_2}_t - Z^{i_1 i_2}_s
\bigr|^2\bigr] \leq C |t-s|^{ 4H-2 },
\]
and by the hypercontractivity property, we deduce
%
%
\begin{equation}
\label{eqn 33} \mE\bigl[\bigl|Z^{i_1 i_2}_t - Z^{i_1 i_2}_s
\bigr|^p\bigr] \leq C_p |t-s|^{p(2H-1) }
\end{equation}
for any $p\geq2$ and $s,t\in[0,T]$. By the Kolmogorov continuity
criterion this implies that
$Z^{i_1 i_2}$ has a H\"older continuous version of exponent $\la$
for any $ \lambda< 2H-1$.

\section{Estimates for solutions of some SDEs}\label{subsec estimate
solution of SDE}

The purpose of this section is to provide upper bounds for the H\"older
seminorms of solutions
of two types of SDEs. The first type [see (\ref{e 3.1})] covers equation
(\ref{e.1.1}) and its
Malliavin derivatives, as well as all the other SDEs involving only
continuous integrands which
we will encounter in this paper. The second type [see (\ref{e3.2})]
deals with the case where the
integrands are step processes. These SDEs arise from approximation
schemes such as (\ref{e.1.2}) and (\ref{e.1.6}).

For any integers $k, N, M\ge1$, we denote by $C^k_b(\mathbb
{R}^M;\mathbb{R}^N)$
the space of $k$ times continuously differentiable functions
$f\dvtx\mathbb{R}^M\rightarrow\mathbb{R}^N$ which are bounded together
with their first $k$ partial derivatives.
If $N=1$, we simply write $C^k_b(\mathbb{R}^M)$.

In order to simplify the notation we only consider the case when the
fBm is one-dimensional, that is, $m=1$.
All results of this section can be generalized to the case $m >1$.
Throughout the remainder of the paper we let $\beta$ be any number
satisfying $\frac{1}{2} < \beta< H$.
The first two lemmas are path-wise results, and they will still hold
when $B$ is
replaced by general H\"older continuous functions of index $\gamma
>\beta$.
The constants appearing in the lemmas depend on $\beta$, $H$, $T$ and the
uniform and H\"older seminorms of the coefficients.
We fix a time interval $[\tau, T]$, and to simplify we omit the
dependence on
$\tau$ and $T$ of the uniform norm and $\beta$-H\"older seminorm on
the interval $ [\tau, T]$.

%
\begin{lemma}\label{lem3.3}
Fix $\tau\in[0,T)$.
Let $V=\{V_t, t\in[\tau,T]\} $ be an $\mathbb{R}^{M }$-valued
processes satisfying
%
%
\begin{equation}
\label{e 3.1} V_t = S_{t } + \int^t_{ \tau}
\bigl[g_{1 } (V_u) + U^1_u
V_u \bigr] \,d u + \int^t_{ \tau}
\bigl[g_{2 }(V_u) + U^2_u
V_u \bigr] \,dB_u,
\end{equation}
where $g_1 \in C_b ( \mathbb{R}^M; \mathbb{R}^M), g_2 \in C^1_b (
\mathbb{R}^M; \mathbb{R}^M)$ and $U^i=\{U^i_t, t\in[\tau,T]\}$,
$i=1,2$, and $ S =\{S_t, \in[\tau, T] \}$ are $\mathbb{R}^{M \times
M}$-valued and $\mathbb{R}^{M }$-valued processes, respectively. We
assume that $S$ has $\beta$-H\"older continuous trajectories, and the
processes $U^i$, $i=1,2$, are uniformly bounded by a constant $C$.

\begin{longlist}[(ii)]
\item[(i)] If $U^1=U^2=0$, then
we can find constants $K $ and $K'$ such that $(t-s)^\beta\|B\|_{\beta
} \le K $, $\tau\le s<t\le T$ implies
\[
\|V \|_{s, t, \beta} \leq K' \bigl(\|B \|_{\beta} +1\bigr) + 2\|S
\|_{\beta}.
\]

\item[(ii)] Suppose that there exist constants $K_0 $ and $K'_0$ such that
$(t-s)^\beta\|B\|_{\beta} \le K_0 $, $\tau\le s<t\le T$ implies
%
%
\begin{equation}
\label{eqn3} \bigl\|U^2\bigr \|_{s, t, \beta} \leq K'_0
\bigl(\|B \|_{\beta} +1\bigr).
\end{equation}
Then there exists a positive constant $K $
such that
%
%
\begin{equation}
\label{est1} \max\bigl\{ \| V \|_{ \infty}, \|V \|_{ \beta} \bigr\} \leq
Ke^{K \|B\|_{\beta}^{ {1}/{\beta}}} \bigl(|S_{\tau}| + \|S\|_\beta+1\bigr).
\end{equation}
\end{longlist}
\end{lemma}

\begin{pf} The proof
follows the approach used, for instance, by Hu and Nualart~\cite{HuNu}.
Let $ \tau\leq s < t \leq T$. By the definition of $V$,
%
%
\begin{equation}\qquad
\label{e.3.16} %
V_{t } -V_s =
S_{t }-S_s + \int^t_{ s}
\bigl[g_1 (V_u) + U^1_uV_u
\bigr] \,d u + \int^t_{ s} \bigl[g_{2}(V_u)
+ U^2_uV_u \bigr] \,dB_u.
\end{equation}
Applying Lemma~\ref{lem7.1}(ii) to the vector valued function $f\dvtx
( u
,v ) \rightarrow g_{2}(v) + u v $
and the integrator $z=B$, and taking $\beta'=\beta$ yields
%
%
\begin{eqnarray}
\label{e3.4} |V _t - V _s| & \leq& \|S
\|_\beta(t-s)^{\beta} + \bigl(\|g_{1}\|_{\infty} + C
\|V \|_{s, t, \infty}\bigr) (t-s)
\nonumber\\
&&{} + K_1 \bigl( \|g_{2}\|_{\infty} +C \|V
\|_{s,t,\infty}\bigr) \|B \|_{\beta} (t-s)^{\beta}
\nonumber
\\[-8pt]
\\[-8pt]
\nonumber
&&{} + K_2 \bigl(\| \nabla g_{2}\|_{\infty} + C\bigr) \|V
\|_{s,t, \beta} \|B \|_{\beta} (t-s)^{2\beta}
\nonumber
\\
& &{}+ K_2 \|V \|_{s,t, \infty} \bigl\|U^2
\bigr\|_{s,t, \beta} \|B \|_{\beta} (t-s)^{2\beta}.
\nonumber
\end{eqnarray}

\textit{Step} 1.
In the case $U^1 = U^2=0$ (which means that we can take $C=0$ and $ \|
U^2 \|_{s,t, \beta} =0$),
dividing both sides of (\ref{e3.4}) by $(t-s)^{\beta}$ and taking the
H\"
older seminorm
on the left-hand side, we obtain
%
%
\begin{eqnarray}
\label{e3.5} \|V\|_{s,t, \beta} &\leq&\|S\|_\beta+ c_1
(t-s)^{1 - \beta}+ K_1 c_1 \|B \|_{\beta}
\nonumber
\\[-8pt]
\\[-8pt]
\nonumber
&&{} +
K_2c_1 \|V \|_{s,t, \beta} \|B\|_{\beta}
(t-s)^{\beta},
\end{eqnarray}
where (and throughout this section) we denote
%
%
\begin{equation}
\label{eqn3.6} c_1 = \max\bigl\{C, \|g_{1}\|_{\infty},
\|g_{2}\|_{\infty}, \| \nabla g_{2}
\|_{\infty} \bigr\}.
\end{equation}
Take $K = \frac{1}2 (K_2 c_1)^{-1}$. Then for any $\tau\le s<t \le T$
such that $(t-s)^\beta\|B\|_{\beta} \le K $, we have
\[
\|V\|_{s,t, \beta} \leq2 \|S\|_\beta+ 2 c_1
(t-s)^{1 - \beta}+ 2 K_1 c_1 \|B \|_{\beta},
\]
which implies (i).

\textit{Step} 2.
As in step 1, we divide (\ref{e3.4}) by $(t-s)^{\beta}$ and then take
the H\"older seminorm
on the left-hand side to obtain
%
%
\begin{eqnarray}
\label{eqn7} 
\|V \|_{s,t, \beta} & \leq& \|S \|_\beta+
c_1\bigl( 1+ \|V \|_{s, t, \infty}\bigr) (t-s)^{1 - \beta}
\nonumber\\
&&{}+ K_1 c_1\bigl(1+ \|V \|_{s,t,\infty}\bigr) \|B
\|_{\beta}
\nonumber
\\[-8pt]
\\[-8pt]
\nonumber
&&{}+ 2K_2c_1 \|V \|_{s,t, \beta} \|B\|_{\beta}
(t-s)^{\beta}
\\
&&{}+ K_2 \|V \|_{s,t, \infty} \bigl\|U^2
\bigr\|_{s,t, \beta} \|B\|_{\beta} (t-s)^{\beta}.
\nonumber
\end{eqnarray}
If $(t-s)^\beta\|B\|_{\beta} \le\frac{1}4 (K_2 c_1)^{-1}$, then
the coefficient of $\|V \|_{s,t, \beta}$ on the right-hand side of
(\ref{eqn7}) is less or equal than $\frac{1}2$.
Thus we obtain
\begin{eqnarray*}
\|V\|_{s,t, \beta} & \leq& 2\|S\|_\beta+2c_1 \bigl(1+ \|V
\|_{s, t, \infty}\bigr) (t-s)^{1-\beta}
\\
&&{}+ 2 K_1 c_1\bigl(1+ \|V \|_{s,t,\infty} \bigr)\|B
\|_{\beta}
\nonumber
\\
&& {}+2 K_2 \|V \|_{s,t, \infty} \bigl\|U^2
\bigr\|_{s,t, \beta} \|B\|_{\beta} (t-s)^{\beta}.
\end{eqnarray*}
On the other hand, assuming $(t-s)^\beta\|B\|_{\beta} \le K_0$ and
applying (\ref{eqn3}), we obtain
%
%
\begin{equation}
\label{a2} 
\|V\|_{s,t, \beta} \leq2 \|S \|_\beta+
C_1 \bigl(1+ \| B \|_\beta\bigr) \bigl(1+ \|V\|_{s,t,\infty}\bigr),
\end{equation}
for some constant $C_1$. This implies
\[
\|V\|_{s,t, \infty} \leq|V_s| + 2 (t-s)^\beta\|S
\|_\beta+ C_1(t-s) ^\beta\bigl(1+ \| B
\|_\beta\bigr) \bigl(1+ \|V\|_{s,t,\infty}\bigr).
\]
Assuming $(t-s) ^\beta\| B \|_\beta\le\frac{1}{4C_1}$ and $(t-s)
^\beta\le\frac{1}{4C_1}\wedge\frac{1}2$, we obtain
%
%
\begin{equation}
\label{a1} \|V\|_{s,t, \infty} \leq2|V_s|+ 2\|S
\|_\beta+ 1.
\end{equation}
Take $\Delta= [ \| B\|_\beta^{-1 } \min( \frac{1}{ 4K_2
c_1}, K_0, \frac{1}{4C_1} ) ] ^{1/\beta} \wedge
( \frac{1}{4C_1} \wedge\frac{1}2 ) ^{1/\beta}$.
We divide the interval $[\tau, T]$ into $N = \lfloor\frac{ T-\tau}{
\Delta} \rfloor+1$
subintervals and denote by $s_1, s_2, \ldots, s_{ N}$ the left
endpoints of these intervals and $s_{N+1}=T$.
Applying inequality (\ref{a1}) to each interval $[s_i, s_{i+1}]$ for $
i=1,\ldots, N$ yields
%
%
\begin{equation}
\label{e.3.21} \|V \|_{ \infty} \leq2^{ N +1} \bigl(|S_{\tau}|
+ 2\|S\|_\beta+1 \bigr).
\end{equation}
From the definition of $\Delta$ we get
%
%
\begin{equation}
\label{e.3.23} N \leq1+ \frac{T} \Delta\leq1+T \max\bigl(
C_2, C_3 \|B\|^{ {1}/{\beta} }_{\beta} \bigr),
\end{equation}
for some constants $C_2$ and $C_3$.
From inequalities (\ref{e.3.21}) and (\ref{e.3.23}) we obtain the
desired estimate for $ \|V\|_{ \infty}$.

If $t, s \in[\tau, T] $ satisfy $ 0 \leq t-s \leq\Delta$, then from
(\ref{a2}) and
from the upper bound of $ \|V\|_{\infty}$ we can estimate $\frac
{V_t-V_s}{(t-s)^\beta}$
by the right-hand side of (\ref{est1}) for some constant~$K$. On the
other hand, if
$t-s > \Delta$,
then
\[
\frac{ |V_t - V_s |}{(t - s)^{\beta}} \leq2\| V\|_{ \infty} \Delta^{-1}.
\]
We can obtain a similar estimate from the upper bound of $ \|V\|_{
\infty}$ and from the definition of $\Delta$.
This gives then the desired estimate for
$\|V\|_{ \beta}$, and hence we complete the proof of (ii).
\end{pf}

For the second lemma we fix $n$ and consider the partition of $[0,T]$ given
by $t_i = i \frac{T}{n} $, $ i =0, 1, \ldots, n$. Define $\eta(t)=
t_i $ if
$t_i\leq t < t_i+ \frac{T}n $ and $\varepsilon(t)= t_i+\frac{T}n $ if
$t_i< t \le t_i+\frac{T}n$.

%
\begin{lemma} \label{lemma3.2}
Suppose that
$ S $, $g_{i}$, $ U ^i$, $i=1,2$ are the same as in Lemma~\ref
{lem3.3}. Let $g \in C([0, T])$.
Let $V=\{V_t, t\in[\tau,T]\} $ be an $\mathbb{R}^{M}$-valued processes
satisfying the equation
%
%
\begin{eqnarray}
\label{e3.2} V_t &=& S_{t } + \int^{t \vee\ep(\tau) }_{ \ep(\tau
) }
\bigl[g_{1} ( V _{\eta(u ) }) + U ^1_{\eta(u ) }
V _{\eta(u ) } \bigr] g \bigl(u-\eta(u) \bigr) \,d u
\nonumber
\\[-8pt]
\\[-8pt]
\nonumber
&&{} + \int^{t \vee\ep(\tau) }_{ \ep( \tau)} \bigl[g_{2} ( V
_{\eta(u ) }) + U _{\eta(u ) } ^2 V _{\eta(u ) } \bigr]
\,dB_u.
\nonumber
\end{eqnarray}

\begin{longlist}[(ii)]
\item[(i)] If $U^1= U^2= 0 $, then
we can find constants $K $ and $K'$ such that $(t-s)^\beta\|B\|_{\beta
} \le K$, $\tau\le s<t\le T$ implies
\[
\| V \|_{s, t, \beta, n} \leq K' \bigl(\|B \|_{\beta} +1\bigr) + 2\|S
\|_{\beta
}.
\]

\item[(ii)] Suppose that there exist constants $ K _0 $ and $ K _0'$ such that
$(t-s) ^\beta\|B\|_{\beta} \leq K_0$, $\tau\le s<t\le T$ implies
%
%
\begin{equation}
\label{eqn3'}\bigl \| U^2 \bigr\|_{s, t, \beta, n} \leq K'
_0 \bigl(\|B \|_{\beta} +1\bigr).
\end{equation}
Then there exists a constant $K $
such that
\[
\max\bigl\{ \| V \|_{ \infty},\| V \|_{ \beta} \bigr\} \leq Ke^{K \|B\|_{\beta
}^{ {1}/{\beta}}}
\bigl( |S_{\tau}| + \|S\|_\beta+1 \bigr).
\]
\end{longlist}
\end{lemma}

%
\begin{remark} The proof of this result is similar to that of
Lemma~\ref{lem3.3}.
Nevertheless, since the integral is discrete,
we need to replace
the H\"older seminorm $\| \cdot\|_{s,t,\beta}$ by the seminorm $\|
\cdot\|_{s, t, \beta, n}$ introduced in (\ref{e2.1}).
\end{remark}

\begin{pf*}{Proof of Lemma \ref{lemma3.2}}
Let $s, t \in[ \tau, T]$ be such that $ s < t$ and $s=\eta(s)$. This
implies $s\ge\varepsilon(\tau)$.
As in the proof of (\ref{e3.4}), applying Lemma~\ref{lem7.1}(i)
[instead of Lemma~\ref{lem7.1}(ii)] yields
\begin{eqnarray*}
&&|V _t - V _s| \\
&&\qquad\leq \|S \|_\beta(t-s)^{\beta}
+ \bigl(\|g_{1}\|_{\infty} + C \|V \|_{s, t, \infty}\bigr)\|g
\|_{\infty} (t-s)
\\
&&\qquad\quad{}+ K_1 \bigl( \|g_{2}\|_{\infty} +C \|V
\|_{s,t,\infty}\bigr) \|B \|_{\beta} (t-s)^{\beta}
\\
&&\qquad\quad{} + K_3\bigl[ \bigl(\| \nabla g_{2}\|_{\infty} + C\bigr)
\|V \|_{s,t, \beta,n} + \|V \|_{s,t, \infty} \bigl\|U^2
\bigr\|_{s,t, \beta
,n}\bigr]\|B \|_{\beta} (t-s)^{2\beta}.
\end{eqnarray*}
Dividing both sides of the above inequality by $(t-s)^{\beta}$ and taking
the H\"older seminorm on the left-hand side, we obtain
%
%
\begin{eqnarray}
\label{eqn3.14} 
\|V\|_{s,t,\beta, n} &\leq& \|S \|_\beta+ \bigl(
\|g_{1}\|_{\infty} + C \|V \|_{s, t, \infty}\bigr) \| g
\|_{\infty} (t-s)^{1-\beta}\nonumber
\\
&&{}+ K_1 \bigl( \|g_{2}\|_{\infty} +C \|V
\|_{s,t,\infty}\bigr) \|B \|_{\beta}
\nonumber
\\[-8pt]
\\[-8pt]
\nonumber
&&{} + K_3 \bigl(\| \nabla g_{2}\|_{\infty} + C\bigr) \|V
\|_{s,t, \beta,n} \|B \|_{\beta} (t-s)^{ \beta}
\\
&&{} + K_3 \|V \|_{s,t, \infty} \|U^2
\|_{s,t, \beta,n} \|B \|_{\beta} (t-s)^{ \beta}.
\nonumber
\end{eqnarray}

\textit{Step} 1. In the case $U^1=U^2= 0$, (\ref{eqn3.14}) becomes
\begin{eqnarray*}
&&\| V \|_{s, t, \beta, n}\\
&&\qquad \leq\|S \|_\beta+ c_1 \| g
\|_{\infty} (t-s)^{1 - \beta} + K_1 c_1 \|B
\|_{\beta} + K_3c_1 \| V \|_{s,t, \beta, n} \|B
\|_{\beta} (t-s)^{ \beta},
\end{eqnarray*}
where $c_1$ is defined in (\ref{eqn3.6}).
Taking $K = \frac{1}2 ( K_3 c_1)^{- {1} }$,
for any $\tau\le s<t\le T$ such that $(t-s)^\beta\|B\|_{\beta} \le
K$, we have
\[
\| V \|_{s,t, \beta, n } \leq2 \|S\|_\beta+ 2 c_1 \| g
\|_{\infty} (t-s)^{1 - \beta}+ 2 K_1 c_1 \|B
\|_{\beta}.
\]
This completes the proof of (i).

\textit{Step} 2.
In the general case, we follow the proof of Lemma~\ref{lem3.3},
except that we assume $s= \eta(s)$ and use the seminorm $\|\cdot\|
_{s,t, \beta, n}$
instead of $\| \cdot\|_{s, t, \beta}$. We also apply
(\ref{eqn3'}) instead of (\ref{eqn3}).
In this way we obtain inequality (\ref{a2}) with $\|V\|_{s, t, \beta}$
replaced by $\|V\|_{s, t, \beta, n}$, that is,
%
%
\begin{equation}
\label{a2'} 
\|V\|_{s,t, \beta, n } \leq 2 \|S \|_\beta+
C_1 \bigl(1+ \| B \|_\beta\bigr) \bigl(1+ \|V\|_{s,t,\infty}\bigr)
\end{equation}
for some constant $C_1$. Inequality
(\ref{a1}) remains the same,
%
%
\begin{equation}
\label{est2} \|V\|_{s,t, \infty} \leq 2|V_s|+ 2\|S
\|_\beta+ 1,
\end{equation}
provided $s=\eta(s)$, and both $t-s$ and $(t-s)^\beta\|B\|_\beta$ are
bounded by some constant~$C_4$.

Take $\Delta= (C_4^{1/\beta} \|B\|_{\beta}^{-1/\beta}) \wedge C_4$.
We are going to consider two cases
depending on the relation between $\Delta$ and $\frac{2T}n$.

If $ \Delta> \frac{2T}{n} $, we take $N = \lfloor\frac{2(T-\ep
(\tau
) ) }{ \Delta}\rfloor$ and divide the interval $[\ep(\tau), \ep
(\tau
)+ N\frac{\Delta}{2}] $ into ${N} $
subintervals of length $\frac{\Delta}{2}$. Since the length of each of
these subintervals is larger than $\frac{T}n $, we are able to choose
$N$ points $s_1, s_2, \ldots, s_{N}$ from each of these intervals such
that $s_1 =\ep(\tau)$ and $\eta(s_i) = s_i$, $i=1, 2, \ldots, {N}$. On
the other hand, we have $s_{i+1} -s_i \le\Delta$ for all $i=1,\ldots
, N-1$.
Applying inequality (\ref{est2}) to each of the intervals $[s_1,
s_{2}]$, $[s_2, s_{3}], \ldots, [s_{N-1}, s_{N }]$, $[s_{N}, T]$ yields
%
%
\begin{equation}
\label{e 3.21} \|V \|_{\ep(\tau),T, \infty} \leq 2^{ N+1 }
\bigl(|S_{\ep(\tau
)}| + 2\|S\|_\beta+1 \bigr).
\end{equation}
From the definition of $\Delta$ we have
%
%
\begin{equation}
\label{e 3.15} N \leq\frac{ 2 T}{ \Delta} \leq K + K \|B\|_{\beta
}^{ {1}/{\beta}},
\end{equation}
for some constant $K$ depending on $T$ and $C_4$.
From (\ref{e 3.21}) and (\ref{e 3.15}) and taking into account that
%
%
\begin{equation}
\label{eqn3.18} \|V\|_{\tau, \ep(\tau), \infty}=\|S \|_{\tau, \ep
(\tau), \infty
} \leq
|S_{\tau}| +T^\beta\|S\|_{\beta},
\end{equation}
we obtain the desired estimate for $\|V\|_{\infty}$.

If $ \Delta\leq\frac{2T}{n}$, that is, when $ n \leq\frac{2T}{
\Delta} \leq K + K \|B\|_{\beta}^{ {1}/{\beta}} $, then by equation
(\ref{e3.2}) we have
\begin{eqnarray*}
| {V}_t| &\leq& | {V}_{\eta(t)} |+ |
S_t - S_{\eta(t)} | + \bigl( c_1+C | {V}_{\eta(t)}
| \bigr) \|g\|_{\infty} (T/n)
\\
&&{} + \bigl( c_1 +C | {V}_{\eta(t)} | \bigr) \|B\|_{\beta} (T/n
)^{\beta}
\\
& \leq&A_n + B_n \bigl| {V}_{\eta(t)} \bigr|, 
\end{eqnarray*}
for any $t\in[\tau,T]$, where
\[
A_n= \| S\|_\beta(T/n)^\beta+ c_1\|g
\|_{\infty} (T/n) +c_1 \|B\| _{\beta}
(T/n)^{\beta}
\]
and
\[
B_n= 1 + C \|g\|_{\infty} ( T/n) + C \|B\|_{\beta}
(T/n)^{\beta}.
\]
Iterating this estimate, we obtain
%
%
\begin{eqnarray}
\label{e.3.33} \|{V} \|_{\ep(\tau), T, \infty} &\leq& |S_{\ep
(\tau)} |
B_n^n + n A_nB_n^{n-1}
\nonumber
\\[-8pt]
\\[-8pt]
\nonumber
&\leq& K\bigl( |S_{\ep(\tau)} | + \|S\|_\beta+1 \bigr) e^{K \|B\|_\beta^{ 1 /
\beta} },
\end{eqnarray}
for some constant $K$ independent of $n$, where we have used the inequality
\[
B_n^n \le e^{K(1+ \|B\|_\beta) n^{1-\beta}},
\]
and the fact that $ n \leq K + K \|B\|_{\beta}^{ {1}/{\beta}} $ for
some constant $K$.
Taking (\ref{eqn3.18}) into account, we obtain the desired upper bound
for $\|V\|_{\infty}$.

In order to show the upper bound for $ \| {V} \|_{\tau, T, \beta}$,
we notice that if $ 0 \leq t-s \leq\Delta$,
then from (\ref{a2'}) and from the upper bound of $\| {V} \|_{\tau, T,
\infty}$, we have
\[
\| {V}\|_{\ep(s),t, \beta,n } \leq K \bigl( |S_\tau| + \|S\|_\beta+1
\bigr)e^{K
\|B\|_\beta^{ 1/\beta} },
\]
for some constant $K$. Thus
\begin{eqnarray*}
\frac{|V_t-V_s| }{(t-s)^{\beta} }& \leq& \| {V}\|_{\ep(s),t, \beta
,n } + \frac{|V_{\ep(s) }-V_s| }{(\ep(s)-s)^{\beta} }
\\
& \leq& K \bigl( |S_\tau| + \|S\|_\beta+1 \bigr)e^{K \|B\|_\beta^{ 1/\beta} }.
\end{eqnarray*}
If $t-s \geq\Delta$, we can obtain the upper bound of $ \| {V} \|_{
\beta}$ by an argument similar to that in the proof of Lemma~\ref{lem3.3}.
The proof of (ii) is now complete.
\end{pf*}

The following result gives upper bounds for the norm of Malliavin
derivatives of the
solutions of the two types of SDEs, (\ref{e 3.1}) and (\ref{e3.2}).
Given a
process $P=\{P_t, t\in[\tau, T]\}$ such that $P_t \in\mathbb
{D}^{N,2}$, for each $t$ and some $N\ge1$,
we denote by $ \mathscr{D}^*_{ N } P$ the maximum of the supnorms of
the functions
$ P_{r_0}$, $ D_{r_1} P_{r_0 } , \ldots, D^N_{r_1, \ldots, r_N}
P_{r_0} $
over $ r_0, \ldots, r_N\in[\tau, T] $,
and denote by $\mathscr{D}_{ N } P$ the maximum of the random variable
$ \mathscr{D}^*_{ N } P$ and the supnorms of $ \|P\|_{\beta}$, $\|
D_{r_1} P\|_{r_1, T, \beta}, \ldots, \|D^N_{r_1, \ldots, r_N} P\|_{r_1
\vee\cdots\vee r_N, T, \beta}
$
over $ r_0, \ldots, r_N\in[\tau, T] $.
If $N=0$, we simply write $\mathscr{D}^*_{ 0 } P= \|P\|_{\infty} $ and
$\mathscr{D}_{ 0 } P=\max( \|P\|_{\infty}, \|P\|_{\beta}
)$.

%
\begin{lemma} \label{lemma3.4}
\textup{(i)} Let $V $ be the solution of equation (\ref{e 3.1}). Assume that
$g_{1}=g_{2}=0$.
Suppose that $U^1$ are $U^2$ are uniformly bounded by a constant $C$,
and assume that there exist constants $ K _0 $ and $ K _0'$ such that
$(t-s) ^\beta\|B\|_{\beta} \leq K_0$, $\tau\le s<t\le T$ implies
%
%
\begin{equation}
\label{eqn3.21} \bigl\| U^{2 }\bigr\|_{s, t, \beta} \leq K'
_0 \bigl(\|B \|_{\beta} +1\bigr).
\end{equation}
Suppose that $S, U^{1 }, U^{2 } \in\mathbb{D}^{N,2}$, where $N\ge0$ is
an integer, and $D_rS_t =D_rU^{i}_t =0$, $i=1,2,$ if $0 \leq t < r \leq
T$, and suppose that
there exists a constant $K>0$ such that the random variables $ \mathscr
{D}_{ N } S $,
$ \mathscr{D}^*_{ N } U^{1} $, $ \mathscr{D}_{ N } U^{2} $
are less than or equal to
$ Ke^{K \|B\|^{ {1}/{\beta} }_{\beta} }$.
Then there exists a constant $K' >0$ such that
$\mathscr{D}_{ N } V $
is less than
$ K'e^{K' \|B\|^{ {1}/{\beta} }_{\beta} }$.

\textup{(ii)} Let $V$ be the solution of equation (\ref{e3.2}).
Then the conclusion in \textup{(i)} still holds true under the same
assumptions, except
that in (\ref{eqn3.21}) we replace $ \| U^{2}\|_{s, t, \beta} $ by $
\| U^{2 }\|_{s, t, \beta, n} $.
\end{lemma}

\begin{pf}
We first show point (i). The upper bounds of $ \|V \|_{\infty}$ and $
\|
V \|_{\beta} $ follow from Lemma~\ref{lem3.3}(ii).
The Malliavin derivative $D_r V _t$ satisfies the equation (see
Proposition~7 in \cite{NuSau})
\[
D_rV _t = S^{(1)}_t +
\int^t_{ r} U^{1} _u
D_r V _u\, d u + \int^t_{ r}
U^{2} _u D_r V _u
\,dB_u
\]
while $t \in[r \vee\tau, T]$ and $D_r V _t = 0 $ otherwise, where
%
%
\begin{equation}
\label{eqn11} S^{(1)}_t:= D_r S
_{t } + U^{2}_r V _r+ \int
^t_{ r} \bigl[ D_r
U^{1}_u \bigr] V_u\, d u + \int
^t_{ r} \bigl[ D_r
U^{2}_u \bigr] V _u \,dB_u
\end{equation}
for $t\in[ r\vee\tau,T]$.
Lemma~\ref{lem3.3}(ii) applied to the time interval $[r,T]$, where
$r\ge\tau$, implies that
\[
\max\bigl\{ \| D_rV \|_{r, T, \infty}, \|D_rV
\|_{ r,T,\beta} \bigr\} \leq K e^{ K \|B\|_{\beta}^{ {1}/{\beta}}}
\bigl( \bigl| S^{(1)}_r\bigr|
+ \bigl\|S^{(1)} \bigr\|_{r,T, \beta} +1\bigr).
\]
Therefore, to obtain the desired upper bound it suffices to show that
there exists
a constant $K$ independent of $r$ such that both $\|S^{(1)} \|_{r,T,
\infty}$
and $ \|S^{(1)}\|_{r,T, \beta} $ are less than or equal to $ K e^{K \|
B\|_{\beta}^{{1}/{\beta} } }$.
Applying Lemma~\ref{lem7.1}(ii) to the second integral in (\ref
{eqn11}) and
noticing that $\|D_rU^2 \|_{\infty}$, $\|D_r U^2 \|_{r, T,\beta}$, $\|
V \|_{\infty}$,
$\|V \|_{r,T, \beta}$ are bounded by $K e^{K \|B\|_{\beta
}^{{1}/{\beta
} } }$, we see
that the upper bound of $\|S^{(1)}\|_{\infty} $ is bounded by $K e^{K
\|B\|_{\beta}^{{1}/{\beta} } }$.
On the other hand, in order to show the upper bound for $\|S^{(1)}\|
_{r,T, \beta} $, we
calculate $ \frac{S^{(1)}_t - S^{(1)}_s}{(t-s)^{\beta}} $ using (\ref
{eqn11}) to obtain
\begin{eqnarray*}
\frac{S^{(1)}_t - S^{(1)}_s}{(t-s)^{\beta}} &\leq& \| D_r S \|
_{r,T,\beta} +
(t-s)^{-\beta}\int^t_{ s} \bigl[
D_r U^{1}_u \bigr] V_u\, d u
\\
&&{}+ (t-s)^{-\beta}\int^t_{ s} \bigl[
D_r U^{2}_u \bigr] V _u
\,dB_u.
\end{eqnarray*}
Now we can estimate each term of the above
right-hand side as before.
Taking the supremum over $s,t \in[r, T]$ yields the upper bound of $ \|
S^{(1)}\|_{r,T, \beta} $.

We turn to the second derivative. As before, we are able to find the
equation of $D^2_{r_1,r_2} V _t$; see Proposition~7 in \cite{NuSau}.
The estimates of $ D^2_{r_1,r_2}V _t $ can then be obtained in the same
way as above by applying Lemma~\ref{lem3.3}(ii) and the estimates that
we just obtained for $V_t $ and $D_sV_t $, as well as the assumptions
on $S $
and $U^i $.
The estimates of the higher order derivatives of $V $ can be obtained
analogously.

The proof of (ii) follows along the same lines, except that we use
Lem\-ma~\ref{lemma3.2}(ii) and Lemma~\ref{lem7.1}(i) instead of
Lemma~\ref
{lem3.3}(ii) and Lemma~\ref{lem7.1}(ii).
\end{pf}

%
\begin{remark}\label{remark1} Since $\beta> \frac{1}2$,
from Fernique's theorem we know that\break $ Ke^{K \|B\|^{ {1}/{\beta}
}_{\beta} }$ has finite moments of any order. So Lemma~\ref{lemma3.4}
implies that the uniform norms and H\"older seminorms of the solutions
of (\ref{e 3.1}) and (\ref{e3.2}) and their Malliavin derivatives have
finite moments of any order. We will need this fact in many of our arguments.
\end{remark}

The next proposition is an immediate consequence of Lemma~\ref{lemma3.4}.
Recall that the random variables $\mathscr{D}^*_N P $ and $\mathscr
{D}_N P $
are defined in Section~\ref{subsec estimate solution of SDE}.
%

\begin{prop}\label{prop 3.6}
Let $X$ be the solution of equation (\ref{e.1.1}), and let
$X^n$ be the solution of the Euler scheme (\ref{e.1.2}).
Fix $N\ge0$, and suppose that $b \in C^N_b (\mathbb{R}^d, \mathbb
{R}^d ),
\sigma\in C^{N+1}_b(\mathbb{R}^d, \mathbb{R}^d) $ (recall that we
assume $m=1$).
Then there exists a positive constant $K$ such that the random variables
$\mathscr{D}_N X $ and $\mathscr{D}_{N} X^n $
are bounded by
$ Ke^{K \|B\|^{{1}/{\beta} }_{\beta} }$
for all $n \in\mathbb{N}$. If we further assume $\sigma
\in C^{N+2}_b(\mathbb{R}^d, \mathbb{R}^d) $, then the same
upper bound holds for the modified Euler scheme (\ref{e.1.6}).
\end{prop}

\begin{pf}We first consider the process $X$, the solution
to equation (\ref{e.1.1}).
The upper bounds for $\|X\|_{\infty}$ and $\|X\|_{\beta}$ follow from
Lemma~\ref{lem3.3}(ii).
The Malliavin derivative $D_r X_t$ satisfies the following linear
stochastic differential equation:
%
%
\begin{equation}
\label{eqn12} D_rX_t = \sigma(X_r) + \int
^t_{ r} \nabla b(X_u)
D_r X_u\, d u + \int^t_{ r}
\nabla\sigma(X_u) D_r X_u
\,dB_u,
\end{equation}
while $0 <r \leq t \leq T$, and $D_r X_t = 0 $ otherwise.
Then it suffices to show that
%
%
\begin{equation}
\label{k1} \sup_{r\in[0, T ] } \mathscr{D}_{M}(D_rX)
\leq Ke^{K \|B\|
^{{1}/{\beta
} }_{\beta} },
\end{equation}
for $M=N-1$.
We can prove estimate (\ref{k1}) by induction on $N\ge1$.
Set $S_t= \sigma(X_r)$, $U^1_t= \nabla b (X_t)$ and $U^2_t= \nabla
\sigma(X_t)$.
Applying Lemma~\ref{lem3.3}(i) to $X$ we obtain
that $U^{2}$ satisfies (\ref{eqn3.21}).
Therefore, Lemma~\ref{lemma3.4} implies that (\ref{k1}) holds for $M=0$.
Now we assume that
\[
\sup_{r\in[0, T ] } \mathscr{D}_{M}(D_rX)
\leq Ke^{K \|B\|
^{{1}/{\beta
} }_{\beta} }
\]
for some $0\le M \le N-2$.
It is then easy to see that
\[
\mathscr{D}^*_{M +1}\bigl(U^{1 }\bigr) \vee
\mathscr{D}_{M+1}\bigl(U^{2 }\bigr) \vee
\mathscr{D}_{M+1}(S) \leq Ke^{K \|B\|^{{1}/{\beta} }_{\beta} },
\]
taking into account that $b \in C^N_b (\mathbb{R}^d; \mathbb{R}^d)$,
$\sigma\in C^{N+1}_b (\mathbb{R}^d; \mathbb{R}^d) $, which enables us
to apply Lemma~\ref{lemma3.4} to (\ref{eqn12}) to obtain the upper bound
of the quantity $ \sup_{r\in[0, T ] } \mathscr{D}_{M+1}(D_rX) $.

The estimates of the Euler scheme and the modified Euler scheme and
their derivatives can be obtained in the same way.
We omit the proof, and we only point out that one more derivative of
$\sigma$ is needed for the modified
Euler scheme because the function $\nabla\sigma$ is involved in its
equation.
\end{pf}

\section{Rate of convergence for the modified Euler scheme and related
processes} \label{section strong conv}
The main result of this section is the convergence rate of the scheme
defined by
(\ref{e.1.6}) to the solution of the SDE (\ref{e.1.1}).
Recall that $\gamma_n$
is the function of $n$ defined in (\ref{e.gamma}).

%
\begin{theorem}\label{thmm4.1}
Let $X$ and $X^n$ be solutions to equations (\ref{e.1.1}) and (\ref
{e.1.6}), respectively.
We assume $b \in C_b^3( \mathbb{R}^d; \mathbb{R}^{d} )$, $\sigma\in
C_b^4 (\mathbb{R}^d; \mathbb{R}^{d \times m})$. Then
for any $p \geq1$ there exists a constant $C$ independent of $n$ (but
dependent on $p$) such that
\[
\sup_{0 \leq t \leq T}\mathbb{E} \bigl[ \bigl|X^n_t -
X_t\bigr|^p \bigr]^{
{1}/p} \leq C
\gamma_n^{-1}.
\]
\end{theorem}

\begin{pf}
Denote $Y:=X - X^n$. Notice that $Y$ depends on $n$, but for
notational simplicity we shall omit
the explicit dependence on $n$ for $Y$ and some other processes
when there is no ambiguity. The idea of the proof is to decompose $Y$ into
seven terms [see (\ref{eqn 36}) below] and then study their convergence
rate individually.

\textit{Step} 1. By the definitions of the processes
$X$ and $X^n$, we have
\begin{eqnarray*}
Y_t &=& \intl^t_0 \bigl[
b(X_s) -b\bigl(X^n_s\bigr) + b
\bigl(X^n_s\bigr) - b\bigl(X^n_{\eta(s) }
\bigr) \bigr] \,ds
\\
&&{} + \sum_{j=1}^m \intl^t_0
\bigl[\sigma^j(X_s) -\sigma^j
\bigl(X^n_s\bigr) + \sigma^j
\bigl(X^n_s\bigr) - \sigma^j
\bigl(X^n_{\eta(s) }\bigr) \bigr] \,d B^j_s
\\
&& {}-H \sum_{j=1}^m \intl^t_0
\bigl(\nabla\sigma^j \sigma^j\bigr) \bigl(X^n_{\eta
(s)}
\bigr) \bigl(s- \eta(s)\bigr)^{2H-1}\,ds.
\end{eqnarray*}
By denoting
\begin{eqnarray*}
\sigma^j_0(s) &=& \bigl( \nabla\sigma^j
\sigma^j \bigr) \bigl(X^n_{\eta(s) } \bigr),\qquad
b_1(s) = \int^1_0 \nabla b\bigl(
\theta X_s + (1 - \theta) X^n_s\bigr) \,d
\theta,
\\
\sigma_1^j(s)&=& \int^1_0
\nabla\sigma^j\bigl( \theta X_s + (1 - \theta)
X^n_s\bigr) \,d \theta,
\end{eqnarray*}
we can write
\begin{eqnarray*}
Y_t&=& \intl^t_0 b_1 (s)
Y_s\, d s + \sum_{j=1}^m
\intl^t_0 \sigma_1^j(s)
Y_s \,dB^{j}_s + \intl^t_0
\bigl[ b\bigl(X^n_s\bigr) - b\bigl(X^n_{\eta(s) }
\bigr) \bigr] \,ds
\\
& &{}+ \sum_{j =1}^m \int
_0^t \bigl[ \sigma^j
\bigl(X^n_s\bigr) - \sigma^j
\bigl(X^n_{\eta
(s) } \bigr) \bigr] \,d B^j_s
-H \sum_{j=1}^m \intl^t_0
\sigma_0^j(s) \bigl(s- \eta(s)\bigr)^{2H-1}\,ds.
\end{eqnarray*}

Let $\Lambda^n= \{ \Lambda^n_t, t\in[0,T]\}$ be the
$d\times d$ matrix-valued solution of the following linear SDE:
%
%
\begin{equation}
\label{e.lambdan} \Lambda^n_t =I+ \intl^t_0
b_1 (s) \Lambda^n_s\, d s + \sum
_{j=1}^m \intl^t_0
\sigma_1^j(s) \Lambda^n_s
\,dB^{j}_s,
\end{equation}
where $I$ is the $d\times d$ identity matrix. Applying the chain rule
for the Young integral to $\Gamma^n_t \Lambda^n_t $, where $ \Gamma
^n_t$, $t \in[0, T ]$ is the unique solution of the equation
%
%
\begin{equation}
\label{eqn 2} \Gamma^n_t =I - \intl^t_0
\Gamma^n_s b_1 (s) \,d s - \sum
_{j=1}^m \intl^t_0
\Gamma^n_s \sigma_1^j(s)\,dB^{j}_s,
\end{equation}
for $t \in[0, T]$, we see that
$\Gamma^n_t \Lambda^n_t = \Lambda^n_t \Gamma^n_t= I$ for all $t \in[0,
T]$. Therefore, $(\Lambda^n_t)^{-1}$
exists and coincides with $\Gamma^n_t$.

We can express the process $Y_t$ in terms of $\Lambda^n_t $ as follows:
%
%
\begin{eqnarray}
\label{eqn 5.3} Y_t &=& \intl^t_0
\Lambda_t^n \Gamma_s^n \bigl[ b
\bigl(X^n_s\bigr) - b\bigl(X^n_{\eta(s) }
\bigr) \bigr] \,ds
\nonumber
\\
&&{}+ \sum_{j=1}^m \intl^t_0
\Lambda_t^n \Gamma_s^n \bigl[
\sigma^j\bigl(X^n_s\bigr) -
\sigma^j\bigl(X^n_{\eta(s) } \bigr) \bigr]
\,dB^j_s
\\
&&{} -H \sum_{j=1}^m \intl^t_0
\Lambda_t^n \Gamma_s^n
\sigma_0^j(s) \bigl(s- \eta(s)\bigr)^{2H-1}\,ds.
\nonumber
\end{eqnarray}
The first two terms in the right-hand side of equation (\ref{eqn 5.3})
can be further decomposed as follows:
%
%
\begin{eqnarray}\label{ow1}
&& \intl^t_0 \Lambda_t^n
\Gamma_s^n \bigl[ \sigma^j
\bigl(X^n_s\bigr) - \sigma^j
\bigl(X^n_{\eta(s ) } \bigr) \bigr] \,d B^j_s
\nonumber
\\
&&\qquad = \intl^t_0 \Lambda_t^n
\Gamma_s^n b^j_2 (s) \bigl(s-
\eta(s)\bigr) \,dB^j_s \nonumber\\
&&\qquad\quad{}+ \sum_{i =1}^m
\intl^t_0 \Lambda_t^n
\Gamma_s^n \sigma_2^{j,i }(s)
\bigl(B^i_s - B^i_{\eta(s)}\bigr)
\,dB^{j}_s
\\
&&\qquad\quad + \intl^t_0 \Lambda_t^n
\Gamma_s^n \sigma_3^j(s)
\bigl(s- \eta(s)\bigr)^{2H} \,dB^{j}_s\nonumber
\\
&&\qquad:= I_{2,j}(t) + \sum_{i =1}^m
I_{3,j,i } (t)+ I_{4,j}(t), \nonumber
\end{eqnarray}
where
\begin{eqnarray*}
b^j_2 (s) &= & \intl^1_0
\nabla\sigma^j\bigl( \theta X^n_s + (1 -
\theta) X^n_{\eta(s)}\bigr) b \bigl(X^n_{\eta(s) }
\bigr)\,d \theta,
\\
\sigma_2^{j,i} (s) &=& \intl^1_0
\nabla\sigma^j\bigl( \theta X^n_s + (1 -
\theta) X^n_{\eta
(s)}\bigr) \sigma^i
\bigl(X^n_{\eta(s) } \bigr)\,d \theta,
\\
\sigma^j_3(s) &= & \frac{1}{2}
\intl^1_0 \nabla\sigma^j\bigl( \theta
X^n_s + (1 - \theta) X^n_{\eta(s)}
\bigr) \sum_{l=1}^m \sigma^l
_0(s) \,d \theta
\end{eqnarray*}
and
%
%
\begin{eqnarray}\label{ow2}
&& \Lambda^n_t \intl^t_0
\Gamma^n_s \bigl[ b\bigl(X^n_s
\bigr) - b\bigl(X^n_{\eta(s)
} \bigr) \bigr] \,ds
\nonumber
\\
&&\qquad= \Lambda^n_t \intl^t_0
\Gamma^n_s b_3 (s) \Biggl[ b
\bigl(X^n_{\eta(s)}\bigr) \bigl(s- \eta(s)\bigr) + \sum
_{j=1} ^m\sigma^j \bigl(X^n_{\eta(s) }
\bigr) \bigl( B^j_s - B^j_{\eta(s)}
\bigr)
\\
&&\hspace*{164pt}\qquad\quad{} + \frac{1}{2} \sum_{j=1}^m
\sigma^j_0 (s) \bigl(s- \eta(s)\bigr)^{2H}
\Biggr] \,ds\nonumber
\\
&&\qquad:= I_{11}(t)+ \sum_{j=1}^m
I_{12, j}(t)+ I_{13}(t),\nonumber
\end{eqnarray}
where
$ b_3(s) = \intl^1_0 \nabla b( \theta X^n_s + (1-\theta) X^n_{\eta(s)
} ) \,d\theta$.
We also denote
%
%
\begin{equation}
\label{ow3} I_{5,j}(t) = -H \Lambda^n_t
\intl^t_0 \Gamma^n_s
\sigma_0^j(s) \bigl(s- \eta(s)\bigr)^{2H-1}\,ds.
\end{equation}
Substituting equations (\ref{ow1}), (\ref{ow2}) and (\ref{ow3}) into
(\ref{eqn 5.3}) yields
%
%
\begin{equation}\quad
\label{eqn 36} %
Y = I_{11}+ \sum
_{j=1}^m I_{12, j} +I_{13} +
\sum_{j=1}^m I_{2,j} + \sum
_{j,i =1}^m I_{3,j,i } + \sum
_{j=1}^m I_{4,j}+ \sum
_{j=1}^m I_{5,j}. %
\end{equation}

\textit{Step} 2.
Denote by $ ( \Lambda^n)_i$, $i=1,\ldots, d$, the $i$th columns of
$\Lambda^n$.
We claim that $ ( \Lambda^n)_i$ satisfy the conditions in Lemma~\ref
{lemma3.4} with
$M=d$, $\tau=0$, $U^{1}_t = b_1(t) $, $U^{2}_t = \sigma^j_1(t)$ and $N=2$.
We first show that $U^{2}$ satisfies (\ref{eqn3.21}).
Taking into account that $b \in C^3_b(\mathbb{R}^d; \mathbb{R}^d)$,
$\sigma\in C^4_b(\mathbb{R}^d; \mathbb{R}^{d\times m})$,
it suffices to show that both $X$ and $X^n$ satisfy~(\ref{eqn3.21}).
This is clear for $X$ because of Lemma~\ref{lem3.3}(i). It follows
from Lemma~\ref{lemma3.2}(i) that there exist
constants $K $ and $K'$ such that $(t-s)^\beta\|B\|_{\beta} \le K$,
$0\le s<t \le T$ implies
\[
\bigl\| X^n\bigr \|_{s, t, \beta, n} \leq K'\bigl (\|B
\|_{\beta} +1\bigr).
\]
Notice that
\begin{eqnarray*}
\frac{|X^n_t -X^n_s|}{(t-s)^\beta} &\leq& \frac{|X^n_{t} -X^n_{\ep
(s)}|}{(t-\ep(s))^\beta} +\frac{|X^n_{\ep(s)} -X^n_s|}{(\ep
(s)-s)^\beta}
\\
&\leq& \bigl\|X^n\bigr\|_{s, t, \beta, n} + \frac{ |X^n_{\ep(s)} - X^n_s|
}{\ep(s) -
s }
\end{eqnarray*}
for $t, s\dvtx t \geq\ep(s)$,
where we recall that $\ep(s)=t_{k+1}$ when $s \in(t_{k}, t_{k+1}]$. Therefore,
to verify (\ref{eqn3.21}) for $X^n$ it suffices to show that
\[
\bigl\| X^n \bigr\|_{s, t, \beta} \leq K' \bigl(\|B
\|_{\beta} +1\bigr)
\]
for $s, t \in[t_k, t_{k+1}]$ for some $k$. But this follows
immediately from (\ref{e.1.6}).
On the other hand,
the fact that
$ \mathscr{D}^*_{ 2 } U^{1} $ and $ \mathscr{D}_{ 2 } U^{ 2 } $ are
less than
$ Ke^{K \|B\|^{ {1}/{\beta} }_{\beta} }$ for some $K$ follows from
Proposition~\ref{prop 3.6}, and the
assumption that $b \in C^3_b (\mathbb{R}^d; \mathbb{R}^d)$, $\sigma
\in C^4_b(\mathbb{R}^d; \mathbb{R}^{d\times m}) $, where
$\mathscr{D}^*_2$ and $ \mathscr{D}_2$ are defined in Section~\ref
{subsec estimate solution of SDE}.

In the same way we can show that the columns of $\Gamma^n$ satisfy the
assumptions of Lemma~\ref{lemma3.4}.
As a consequence, it follows from Lemma~\ref{lemma3.4} that
%
%
\begin{equation}
\label{e4.3} \mathscr{D}_2 \Lambda^n \vee
\mathscr{D}_2\Gamma^n \leq K e^{K \|B\|
_{\beta}^{ {1}/{\beta} } }.
\end{equation}

\textit{Step} 3.
From (\ref{e4.3}) and from the fact that $b\in C^3_b(\mathbb{R}^d;
\mathbb{R}^d) $ and $\sigma\in C^4_b (\mathbb{R}^d;\break   \mathbb
{R}^{d\times m}) $,
it follows that
%
%
\begin{equation}
\label{eqn19} \mE\bigl(\bigl|I_{11}(t)\bigr|^p\bigr)^{{1}/p}
\leq C n^{-1} \quad\mbox{and}\quad \mE\bigl(\bigl|I_{13}(t)\bigr|^p
\bigr)^{{1}/p} \leq Cn^{-2H}.
\end{equation}
Notice that $n^{-1} $ and $n^{-2H}$ are bounded by $\gamma_n^{-1}$.
Applying estimates (\ref{11.a}) and~(\ref{11.b}), inequality (\ref{e4.3})
and Proposition~\ref{prop 3.6}, we have for any $j$
%
%
\begin{eqnarray}
\label{eq19a} \mE\bigl(\bigl| I_{12, j } (t)\bigr|^p
\bigr)^{{1}/p} &\leq& Cn^{-1},\qquad \mE\bigl(\bigl| I_{ 2, j }(t)
\bigr|^p \bigr)^{{1}/p} \leq Cn^{-1},
\nonumber
\\[-8pt]
\\[-8pt]
\nonumber
 \mE
\bigl(\bigl|I_{4, j
}(t)\bigr|^p\bigr)^{{1}/p} &\leq&
Cn^{-2H}.
\end{eqnarray}
Now to complete the proof of the theorem it suffices to show that for
any $j$, $\mE(|\sum_{i =1}^m I_{3, j, i}(t) + I_{5, j}(t)|^p)^{{1}/p}
\leq C \gamma_n^{-1}$. For any fixed $j$ we make the decomposition
%
%
\begin{equation}
\label{eqn20} \sum_{i =1}^m
I_{3, j, i} + I_{5, j}= E_{1,j}+E_{2,j}+E_{3,j},
\end{equation}
where
\begin{eqnarray*}
E_{1,j}(t) &= & \Lambda_t^n \sum
_{i =1}^m \intl^t_0 \bigl[
\Gamma_s^n \sigma_2^{j,i}(s) -
\Gamma_{ \eta(s) }^n \bigl(\nabla\sigma^j
\sigma^i\bigr) \bigl(X^n_{\eta(s)}\bigr) \bigr]
\bigl( B^i_{s}- B^i_{\eta(s)} \bigr)
\,dB^{j}_s,
\\
E_{2,j} (t) &=& \Lambda^n_t \sum
_{i=1}^m \int_0^t
\Gamma_{ \eta(s)
}^n \bigl(\nabla\sigma^j
\sigma^i\bigr) \bigl(X^n_{\eta(s)}\bigr) \bigl(
B^i_{s}- B^i_{\eta(s)}
\bigr)
\,dB^{j}_s
\\
&&{} - H \Lambda^n_t \intl^t_0
\Gamma^n_{ \eta(s) } \sigma_0^j(s)
\bigl(s- \eta(s)\bigr)^{2H-1}\,ds,
\\
E_{3,j} (t) &=& H \Lambda^n_t
\intl^t_0 \bigl( \Gamma^n_{ \eta(s) } -
\Gamma^n_s\bigr) \sigma_0^j(s)
\bigl(s- \eta(s)\bigr)^{2H-1}\,ds.
\end{eqnarray*}
Applying (\ref{e4.3}) for the quantities $\|\Lambda^n\|_{\infty}$ and
$\|
\Gamma^n\|_{\beta}$, it is easy to see that $\mE(|E_{3,j} (t )
|^p)^{{1}/p} \leq Cn^{1 - 2H-\beta}$
for any $\frac{1}2<\beta<H$.
On the other hand, applying estimate (\ref{e 11.10}) from Lemma~\ref
{lem11.5} to $E_{1,j}$, we obtain
$\mE(|E_{1,j }(t)|^p)^{{1}/p} \leq Cn^{1-3\beta}$ for any $\frac
{1}2<\beta<H$.
Notice that the exponents $n^{1 - 2H-\beta}$ and
$n^{1 - 3\beta}$ are bounded by $\gamma_n^{-1}$ if $\beta$ is
sufficiently close to $H$.

Taking into account the relationship between the Skorohod and
path-wise integral, we can express the term $E_{2,j}$ as follows:
%
%
\begin{equation}
\label{4.12e} E_{2,j} (t) = \Lambda^n_t \sum
_{i=1}^m \sum
_{k=0}^{
{ \lfloor{nt}/{T} \rfloor}
} F^{n, i, j}_{t_k} \int
_{t_k}^{t_{k+1}\wedge t} \int_{t_k}^s
\delta B^i_u \delta B^j_s,
\end{equation}
for $t\in[0, T]$,
where $ F^{n, i, j}_{t} = \Gamma^n_t ( \nabla\sigma^j \sigma^i )
(X^n_t) $, and we define $t_{n+1} = (n+1) \frac{T}n$. From (\ref{e4.3})
and Proposition~\ref{prop 3.6}, we have
%
%
\begin{equation}
\label{eqn16} \max\bigl\{ \bigl|F^{n, i,j}_t\bigr|,
\bigl|D_{r_1}F^{n, i,j}_t\bigr|,\bigl| D_{r_2 }
D_{r_1}F^{n, i,j} \bigr| \bigr\} \leq K e^{K \|B\|_{\beta}^{{1}/{\beta} } }.
\end{equation}
Hence, applying estimate (\ref{11.c}) from Lemma~\ref{lem11.4} to
$E_{2,j} (t)$, we obtain\break $\mE(|E_{2,j} (t ) |^p)^{{1}/p} \leq
C\gamma_n^{-1}$.
The proof is now complete.
\end{pf}

The following result provides a rate of convergence for the Malliavin
derivatives of the modified scheme and some related processes. Recall
that $\beta$ satisfies $ \frac{1}2 < \beta< H$.
%

\begin{lemma}\label{lem3.5}
Let $X$ and $X^n$ be the processes defined by (\ref{e.1.1}) and
(\ref{e.1.6}), respectively. Suppose that $\sigma\in
C^5_b(\mathbb{R}^d;\mathbb{R}^{d\times m})$, $b \in
C^4_b(\mathbb{R}^d;\mathbb{R}^{d })$. Let $p \geq1$.
Then:
\begin{longlist}[(ii)]
\item[(i)]
There exists a constant $C$ such that
the quantities
$ \| D_{s}X_t - D_{s}X^n_t\|_{p}$, $ \|D_{r}D_{s}X_t -
D_{r}D_{s}X^n_t\|_{p}$,
$ \|D_u D_{r}D_{s}X_t - D_u D_{r}D_{s}X^n_t\|_{p}$
are less than $ C n^{1-2\beta}$ for all $u, r,s,t \in[0,T]$ and
$n\in\mathbb{N}$.

\item[(ii)] Let $V$ and $V^n$ be $d$-dimensional processes
satisfying the equations
\begin{eqnarray*}
\label{eqn 5.1} %
V_t &=& V_0 +
\intl^t_0 f_1 ( X_u,
X_u ) V_u\, d u + \sum_{j=1}^m
\intl^t_0 f_2^j (
X_u, X_u ) V_u \,dB_u^j,
\\
V_t^n &= & V_0 + \intl^t_0
f_1\bigl( X_u, X^n_u\bigr)
V_u^n\, d u + \sum_{j=1}^m
\intl^t_0 f_2^j \bigl(
X_u, X^n_u\bigr) V_u^n
\,dB_u^j,
\end{eqnarray*}
where $f_1\in C^3_b(\mathbb{R}^d \times\mathbb{R}^d; \mathbb{R}^{ d
\times d }
)$ and $ f_2^j \in C^4_b(\mathbb{R}^d \times\mathbb{R}^d; \mathbb
{R}^{ d \times d }
) $.
Then there exists a constant $C$ such that
the quantities $ \|V_t - V^n_t\|_{p}$, $ \| D_{s}V_t - D_{s}V^n_t\|_{p}
$, $ \|D_{r}D_{s}V_t - D_{r}D_{s}V^n_t\|_{p}$ are less than $ C
n^{1-2\beta} $
for all $r, s, t \in[0,T]$ and $n \in\mathbb{N}$.
\end{longlist}
\end{lemma}

%
\begin{remark}\label{remark4.3}
The above results still hold when the approximation process $X^n$ is
replaced by the one defined by the recursive scheme (\ref{e.1.2}). The
proof follows exactly along the same lines.
\end{remark}
\begin{pf*}{Proof of Lemma \ref{lem3.5}}
(i) Taking the Malliavin derivative in both sides of~(\ref{eqn 5.3}),
we obtain
\begin{eqnarray*}
D_r\bigl(X_t-X^n_t\bigr) &=&
\intl^t_0 D_r \bigl[ \Lambda_t^n
\Gamma_s^n \bigl( b\bigl(X^n_s
\bigr) - b\bigl(X^n_{\eta(s) } \bigr) \bigr) \bigr] \,ds
\\
&&{}+ \sum_{j=1}^m \intl^t_0
D_r \bigl[ \Lambda_t^n \Gamma_s^n
\bigl( \sigma^j\bigl(X^n_s\bigr) -
\sigma^j\bigl(X^n_{\eta(s) } \bigr) \bigr) \bigr]
\,dB^j_s
\nonumber
\\
&&{} + \sum_{j=1}^m \Lambda_t^n
\Gamma_r^n \bigl( \sigma^j
\bigl(X^n_r\bigr) - \sigma^j
\bigl(X^n_{\eta
(r) } \bigr) \bigr)
\\
&&{} -H \sum_{j=1}^m \intl^t_0
D_r \bigl[ \Lambda_t^n \Gamma_s^n
\sigma_0^j(s) \bigr] \bigl(s- \eta(s)
\bigr)^{2H-1}\,ds.
\end{eqnarray*}
Proposition~\ref{prop 3.6} and equation (\ref{e4.3}) imply that the
first, third and last terms of the
above right-hand side have $L^p$-norms bounded by $Cn^{1-2H}$.
Applying estimate (\ref{e 11.11}) from Lemma~\ref{lem11.5} to the
second term
and noticing that $\|X\|_\beta$ and $\sup_{r\in[0, T]} \| D_rX\|
_\beta
$ have finite
moments of any order, we see that its $L^p$-norm is also bounded by
$Cn^{1-2\beta}$.

Similarly, we can take the second derivative in (\ref{eqn 5.3}) and
then estimate each
term individually as before to obtain that the upper bound of $
\| D_rD_{s }X_t - D_r D_sX^n_t \|_{p} $ is bounded by $Cn^{1-2\beta}$.

(ii) Using the chain rule for Young's integral we
derive the following explicit expression for $V_t-V_t^n$:
%
%
\begin{eqnarray}
\label{eqn5.4} V_t - V^n_t &=& \int
_0^t \Upsilon_t
\Upsilon_s^{-1} \bigl(f_1 ( X_s,
X_s ) - f_1 \bigl( X_s,
X_s^n \bigr) \bigr) V^n_s \,ds
\nonumber
\\[-8pt]
\\[-8pt]
\nonumber
&&{}+ \sum_{j=1}^m \int
_0^t \Upsilon_t
\Upsilon_s^{-1} \bigl(f_2^j (
X_s, X_s ) - f_2^j \bigl(
X_s, X_s^n \bigr) \bigr)
V^n_s\, d B^j_s ,
\end{eqnarray}
where $\Upsilon=\{\Upsilon_t$, $t\in[0, T]\}$ is the $\mathbb{R}^{d
\times d}$-valued process that satisfies
\[
\Upsilon_t = I+ \int_0^t
f_1 ( X_s, X_s ) \Upsilon_t \,ds
+ \sum_{j=1}^m \int_0^t
f_2^j ( X_s, X_s )
\Upsilon_t\, d B^j_s .
\]
Lemma~\ref{lemma3.4} implies that there exists a constant $K$ such that
for all $n \in\mathbb{N}$, $u, r, s, t \in[0, T]$, we have
%
%
\begin{equation}
\label{eqn18} \max\bigl\{ \Upsilon_t, D_s
\Upsilon_t, D_rD_s \Upsilon_t,
D_u D_rD_s \Upsilon_t\bigr\} \leq
K e^{K \|B\|_\beta^{ 1/\beta}}.
\end{equation}
Therefore, applying estimate (\ref{11.a})
to the second integral in (\ref{eqn5.4}) with $\nu=0$ and taking into
account the estimate of Lemma~\ref{lem3.5}(i),
we obtain
\[
\bigl\|V-V^n\bigr\|_{p} \leq C n^{1-2\beta}.
\]
Taking the Malliavin derivative on both sides of (\ref{eqn5.4}), and then
applying estimates~(\ref{11.a}) from Lemmas \ref{lem11.3} and \ref
{lem3.5}(i) as before,
we can obtain the desired estimate for $\| D_{s}V_t - D_{s}V^n_t\|_{p} $.
The estimate for $ \|D_{r}D_{s}V_t - D_{r}D_{s}V^n_t\|_{p}$ can be
obtained in a similar way.
\end{pf*}

We define $\{\Lambda_t, t\in[0, T]\}$ as the solution of the
limiting equation of (\ref{e.lambdan}), that is,
%
%
\begin{equation}
\label{eqn 46} \Lambda_t = I + \intl^t_0
\nabla b (X_s) \Lambda_s \,ds + \sum
_{j=1}^m \intl^t_0 \nabla
\sigma^j(X_s) \Lambda_s \,dB_s^{j}.
\end{equation}
The inverse of the matrix $\Lambda_t$, denoted by $\Gamma_t$,
exists and satisfies
\[
\Gamma_t =I - \intl^t_0
\Gamma_t \nabla b(X_s) \,d s - \sum
_{j=1}^m \intl^t_0
\Gamma_t \nabla\sigma^j(X_s)
\,dB^{j}_s.
\]
It follows from Lemma~\ref{lem3.5} that if we assume that $\sigma\in
C_b^5(\mathbb{R}^d; \mathbb{R}^{d\times m})$
and $ b \in C_b^4(\mathbb{R}^d; \mathbb{R}^d)$, then the estimate in
Lemma~\ref{lem3.5}(ii) holds
with the pair $(V, V^n)$ being replaced by $(\Gamma_i, \Gamma^n_i)$ or
$(\Lambda_i, \Lambda^n_i)$, $i=1,\ldots, d$,
where the subindex $i$ denotes the $i$th column of each matrix.

\section{Central limit theorem for weighted sums}\label{s.theorem1}
Our goal in this section is to prove a central limit result for
weighted sums
(see Proposition~\ref{prop 3} below) that will play a fundamental role
in the proof of Theorem~\ref{theorem 1}
in the next section. This result has an independent interest and we
devote this entire section to it.

We recall that $B=\{B_t, t\in[0,T]\}$ is an $m$-dimensional fBm, and
we assume that the Hurst parameter satisfies $H\in(\frac{1}2,
\frac{3}4 ]$. For any $n\ge1$ we set $t_j =\frac{jT}n$,
$j=0,\ldots
, n$. Recall that $\eta(s) =t_k$ if $t_k \le s<t_{k+1}$.
Consider the $d\times d$ matrix-valued process
\[
\Xi_t ^{n,i,j}= \gamma_n \sum
_{k=0}^{\{t\}} \int_{t_k}^{ t_{k+1}}
\bigl(B^i_s- B^i_{\eta(s)}\bigr)
\delta B^j_s,\qquad 1\le i,j \le m,
\]
where we denote $\{t\} = { \lfloor\frac{nt}{T} \rfloor} $ for $t\in
[0, T)$ and $\{T\}=t_{n-1}$.
%

\begin{prop} \label{prop4.2}
The following stable convergence holds as $n$ tends to infinity
\[
\bigl(\Xi^n, B \bigr)\rightarrow(W,B),
\]
where $W=\{W_t, t\in[0,T]\}$ is the matrix-valued Brownian motion,
introduced in Section~\ref{sec2.3}, and $W$ and $B$ are independent.
\end{prop}

\begin{pf}
From inequality (\ref{11.c}) in Lemma~\ref{lem11.4} it follows that
%
%
\begin{equation}
\label{eqn4.1} \mE\bigl( \bigl\vert\Xi^n_{ t_{k}} -
\Xi^n_{ t_{j}} \bigr\vert^4 \bigr) \leq C
\biggl(\frac{k-j}{n} \biggr)^{2 },
\end{equation}
for any $j\le k$. This implies the tightness of $(\Xi^n, B )$.

Then it remains to show the convergence of the finite dimensional distributions
of $(\Xi^n, B ) $ to that of $(W,B)$. To do this, we fix a finite set
of points $r_1, \ldots, r_{L+1} \in[0, T]$ such that $0 = r_1 < r_2
< \cdots< r_{L+1} \leq T$ and define the random vectors
$B_L= ( B_{r_2} - B_{r_1}, \ldots, B_{r_{L+1 } } - B_{r_L} )$, $\Xi
^n_L=( \Xi^n_{r_2} - \Xi^n_{r_1}, \ldots, \Xi^n_{r_{L+1 } } - \Xi
^n_{r_L} )$ and $W_L= (W_{r_2} - W_{r_1}, \ldots, W_{r_{L+1 } } -
W_{r_L}) $.
We claim that as $n$ tends to infinity, the following convergence in
law holds:
%
%
\begin{equation}
\label{ow4} \bigl(\Xi^n_L, B_L \bigr)
\Rightarrow(W_L,B_L).
\end{equation}
For notational simplicity, we add one term to each component of $\Xi
^n_L$, and we define
%
%
\begin{equation}
\label{eqn 15} { \Theta}^n_{l } (i, j): =
\Xi^{n,i,j}_{r_{l+1}} - \Xi^{n,i,j}_{r_l} +
\zeta^{i,j}_{\{r_l\}, n} = \ga_n \sum
_{k=\{r_l\} }^{\{r_{l+1}\}} {\zeta}^{ i,j }_{k,n},
\end{equation}
for $ 1\le l \le L$, $ 1\le i,j \le d$, where
\[
{\zeta}^{ i,j }_{k,n} = \intl_{t_k }^{t_{k+1} }
\bigl( B^{i }_s - B^{i}_{t_k}\bigr)
\delta B^{j }_s.
\]
Then Slutsky's lemma implies that the convergence in law in (\ref{ow4})
is equivalent to
\[
\bigl( \Theta^n_l (i,j), 1 \leq i, j \leq d, 1 \leq l
\leq L, B_L \bigr) \Rightarrow(W_L,
B_L).
\]

According to Peccati and Tudor \cite{Peccati} (see also Theorem~6.2.3
in \cite{Peccati2}), to show
the convergence in law of $(\Theta^n_L, B_L)$, it suffices to show the
convergence
of each component of $(\Theta^n_L, B_L)$ to the correspondent
component of $(W_L,B_L)$
and the convergence of the covariance matrix.

The convergence of the covariance matrix of $\Theta^n_L$ follows from
Propositions
\ref{prop 1} and \ref{prop 2} below. The convergence in law of each
component to
a Gaussian distribution follows from Proposition~\ref{prop3} below and
the fourth moment
theorem; see \cite{Nualart3} and also Theorem~5.2.7 in \cite
{Peccati2}. This completes the proof.
\end{pf}

In order to show the convergence of the covariance matrix and the
fourth moment of
$\Theta_n$ we first introduce the following notation:
%
%
\begin{eqnarray}
\label{eqn 4.2} %
\cD_{k }&= & \bigl\{ (s,t,v,u)\dvtx
t_{k }\le v\le s \le t_{k+1}, u, t \in[ 0, T] \bigr\},
\nonumber
\\[-8pt]
\\[-8pt]
\nonumber
\cD_{k_1, k_2 }&= & \bigl\{ (s,t,v,u)\dvtx t_{k_2}\le v\le s \le
t_{k_2+1}, t_{k_1}\le u\le t \le t_{k_1+1} \bigr\}.
\end{eqnarray}

The next two propositions provide the convergence of the covariance $
\mE[
\Theta_{l'}^n (i', j')
\Theta_l^n (i, j)
]$ in the cases $l=l'$ and $l\neq l'$, respectively. We denote
$ \beta_{{k}/n }(s) =
\mathbf{1}_{[t_k, t_{k+1} ] } (s)$.

%
\begin{prop}\label{prop 1}
Let $\Theta^n_l(i,j)$ be defined in (\ref{eqn 15}). Then
%
%
\begin{equation}
\label{eqn 11} 
\mE\bigl[ \Theta_l^n
\bigl(i', j'\bigr) \Theta_l^n
(i, j) \bigr] \rightarrow\alpha_H^2 \frac{ r_{l+1} - r_l }{ T}
( R \delta_{ji'} \delta_{i j'} + Q \delta_{j j'}
\delta_{ii'} ) , 
\end{equation}
as $ n \rightarrow+\infty$.
Here $\delta_{ii'}$ is the Kronecker function, $\alpha_H = H(2H-1)$ and
$Q$ and $R$ are the constants defined in (\ref{e2.6}).
\end{prop}

\begin{pf} The proof will involve several steps.

\textit{Step} 1.
Applying twice the integration by parts formula (\ref{ibp}), we have
%
%
\begin{eqnarray}
\label{eqn4.6} && \mE\bigl[ \Theta^n_l
\bigl(i', j'\bigr) \Theta^n_l
(i, j) \bigr]
\nonumber
\\[-8pt]
\\[-8pt]
\nonumber
&&\qquad = \alpha_H^2 \ga_n \sum
_{k=\{r_l\} }^{\{r_{l+1}\}} \int_{\cD_{k}}
D^i_u D^j_t
\Theta^n_{l } \bigl(i', j'
\bigr) \mu(dv \,du)\mu( ds \,dt),
\end{eqnarray}
where we recall that $\{t\} = { \lfloor\frac{nt}{T} \rfloor} $ for
$t\in[0, T)$ and $\{T\}=t_{n-1}$,
and $\cD_k$ is defined in~(\ref{eqn 4.2}). Since
%
%
\begin{eqnarray}
\label{eqn 48} 
&& D^{i}_u
D^{j }_t \Theta^n_{l }
\bigl(i', j'\bigr)
\nonumber
\\[-8pt]
\\[-8pt]
\nonumber
&&\qquad = \gamma_n \sum_{k=\{r_l\} }^{\{r_{l+1}\}}
\bigl( \mathbf{1}_{[t_k, t]}(u) \beta_{{k} /n}(t)
\delta_{j j'} \delta_{ii'} + \mathbf{1}_{[t_k, u ]}(t)
\beta_{{k}/ n}(u) \delta_{ji'} \delta_{i j'} \bigr), 
\end{eqnarray}
the left-hand side of (\ref{eqn 11}) equals
\begin{eqnarray*}
&& \alpha_H^2 \gamma_n^2
\sum_{k, k'=\{r_l\} }^{\{r_{l+1}\}} \int_{\cD_{k}}
\bigl\{ \mathbf{1}_{[t_{k' }, t]}(u) \beta_{{k' }/ n}(t)
\delta_{j j'} \delta_{ii'}\\
&&\hspace*{80pt}{} + \mathbf{1}_{[t_{k' }, u ]}(t) \beta_{{k' } /n}(u)
\delta_{ji'} \delta_{i j'} \bigr\} \mu(dv \,du )\mu(ds \,dt)
\\
&&\qquad:= \alpha_H^2 \gamma_n^2 (
G_1 \delta_{j j'} \delta_{ii'} +G_2
\delta_{ji'} \delta_{i j'} ) . 
\end{eqnarray*}
In the next two steps, we compute the limits of $\gamma_n^2 G_1$ and
$\gamma_n^2 G_2$ as $n$
tends to infinity in the case $H \in(\frac{1}2, \frac{3}4)$ and in the
case $H = \frac{3}4$
separately.

\textit{Step} 2. In this step, we consider the case $H
\in(\frac{1}2, \frac{3}4) $.
Recall that
\begin{eqnarray*}
Q( p)&=& T^{4H} \intl_{ {0}}^{ { 1}}
\intl_{ { p} }^{ {{ p} +1}} \intl_{ {0}}^{ t }
\intl_{ { p} }^{s} \mu(dv \,du) \mu(ds \,dt)
\\
&=& n^{4H} \int_{\cD_{k', k'+p}} \mu(dv \,du )\mu(ds \,dt) ,
\end{eqnarray*}
which is independent of $n$, where the set $\mathcal{D}_{k_1, k_2}$ is
defined in (\ref{eqn 4.2}).
We can express $ \gamma_n^2 G_1$ in terms of $Q(p)$ as follows:
\begin{eqnarray*}
\gamma_n^2 G_1 &=& n^{4H-1 } \sum
_{k, k' =\{r_l\} }^{\{r_{l+1}\}} \int_{\cD_{k', k}}
\mu(dv \,du)\mu( ds \,dt)
\\
& =& \frac{1}n \sum_{p= \{r_l\} -
\{r_{l+1}\}
}^{ \{ r_{l+1}\} - \{r_{l} \} }
\sum_{k'=( \{r_l\} -p ) \vee\{r_l\} }^{
( \{r_{l+1}\} -p ) \wedge\{ r_{l+1} \}
} Q( p)
\\
&=& \sum_{p=-\infty} ^\infty\Psi^n_l(p)
Q(p),
\end{eqnarray*}
where
\[
\Psi^n_l(p) =\frac{ ( \{r_{l+1}\} -p ) \wedge\{ r_{l+1} \}- ( \{
r_l\} -p ) \vee\{r_l\} }n \mathbf{1}_{[\{r_l\} - \{r_{l+1}\}, \{
r_{l+1}\} - \{r_{l} \} ]}
(p).
\]
The term $ \Psi^n_l(p)$ is uniformly bounded and converges\vspace*{1pt} to $ \frac{
{r_{l+1}} -{r_l} }{T} $ as $n$ tends to infinity for any fixed $p$.
Therefore, taking into account that $ \sum_{p=-\infty} ^\infty
Q(p)=Q<\infty$, the dominated convergence theorem implies
\[
\lim_{n\rightarrow\infty} \gamma_n^2
G_1 = \frac{ {r_{l+1}} -{r_l} }{T} Q.
\]
Similarly, we can show that
\[
\lim_{n\rightarrow\infty} \gamma_n^2
G_2 = \frac{ {r_{l+1}} -{r_l} }{T} R.
\]

\textit{Step 3.} In the case $H=\frac{3}4$,
we can write
\begin{eqnarray*}
\gamma_n^2 G_1 &=& \frac{ n^{2 } }{\log n } \sum
_{k, k'=\{r_l\} }^{\{r_{l+1}\}} \int_{\cD_{k',k}}
\mu(dv \,du )\mu(ds \,dt)
\\
& =& \frac{1}{n\log n} \sum_{p= \{r_{l } \} -
\{r_{l+1 } \}
}^{\{r_{l+1 } \} - \{r_{l } \} }
\sum_{k'= \{r_{l } \} }^{
\{r_{l+1 } \} } Q( p)
\\
&&{} - \frac{1}{n\log n} \Biggl\{ \sum_{p= \{r_{l } \} -
\{r_{l+1 } \}
}^{0 }
\sum_{k'= \{r_{l } \} }^{
\{r_{l } \} -p-1
} + \sum
_{p= 1
}^{\{r_{l+1 } \} - \{r_{l } \} } \sum_{k'= \{r_{l+1 } \} -p +1 }^{
\{r_{l+1 } \}
}
\Biggr\} Q( p)
\\
&:=& G_{11}+G_{12}.
\end{eqnarray*}
Taking into account that $Q(p)$ behaves like $1/|p|$ as $|p|$ tends to
infinity, it is then easy to see that
$G_{12}$ converges to zero. On the other hand, recall that $
Q = \lim_{n \rightarrow+\infty} \frac{ \sum_{|p| \leq n} Q(p)
}{ \log n}
$. This implies that $G_{11}$ converges to $
\frac{Q} T ({r_{l+1}} - {r_l} ) $.
This gives the limit of $ \gamma_n^2 G_1$. The limit of $ \gamma_n^2
G_2$ can be obtained similarly.
\end{pf}

%
\begin{prop}\label{prop 2}
Let $l, l' \in\{1, \ldots, L \}$
be such that $l \neq l'$. Let $\Theta^n$ be defined as in (\ref{eqn 15}).
Then
%
%
\begin{equation}
\label{eqn 4.17} \lim_{n\rightarrow\infty} \mE\bigl[ \Theta_{l'}^n
\bigl(i', j'\bigr) \Theta_l^n
(i, j) \bigr] =0.
\end{equation}
\end{prop}

\begin{pf} Without any loss of generality, we assume
$l' < l$. As in (\ref{eqn4.6}) we have
\[
\mE\bigl[ \Theta_{l'}^n \bigl(i',
j'\bigr) \Theta_l^n (i, j) \bigr] =
\alpha_H^2 \ga_n \sum
_{k=\{r_l\} }^{\{r_{l+1}\}} \int_{\cD_{k}}
D^i_u D^j_t
\Theta^n_{l ' } \bigl(i', j'
\bigr) \mu(dv \,du )\mu(ds \,dt).
\]
Taking into account (\ref{eqn 48}), we can write
\begin{eqnarray*}
&&\mE\bigl[ \Theta_{l'}^n \bigl(i',
j'\bigr) \Theta_l^n (i, j) \bigr] \\
&&\qquad=
\alpha_H^2 \gamma_n^2 \sum
_{k =\{r_l\} }^{\{r_{l+1}\}} \sum
_{ k'=\{r_{l' } \} }^{\{r_{l' +1}\}} \int_{\cD_{k}} \bigl\{
\mathbf{1}_{[t_{k' }, t]}(u) \beta_{{k' }/ n}(t) \delta_{j j'}
\delta_{ii'} \\
&&\hspace*{135pt}{}+ \mathbf{1}_{[t_{k' }, u ]}(t) \beta_{{k' }/ n}(u)
\delta_{ji'} \delta_{i j'} \bigr\} \mu(dv \,du )\mu(ds \,dt)
\\
&&\qquad:= \alpha_H^2 \gamma_n^2 (
\wt{G}_1 \delta_{j j'} \delta_{ii'} +
\wt{G}_2 \delta_{ji'} \delta_{i j'} ).
\end{eqnarray*}
In the case $H \in(\frac{1}2, \frac{3}4)$ we have
\begin{eqnarray*}
\gamma_n^2 \wt{G}_1 & = &n^{4H-1 }
\sum_{k=\{r_l\} }^{\{r_{l+1}\}} \sum
_{ k'=\{r_{l' } \}
}^{\{r_{l' +1}\}} \int_{ \cD_{k}}
\mathbf{1}_{[t_{k' }, t]}(u) \beta_{{k' } /n}(t) \mu(dv \,du ) \mu(
ds \,dt)
\\
&=& \frac{1}{n} \sum_{p= \{r_{l } \} - \{r_{l'+1 } \} }^{ \{r_{l+1 }
\} - \{r_{l' } \} }
\sum_{k'=( \{r_{l } \} -p) \vee\{r_{l' } \} }^{
\{r_{l'+1 } \}
\wedge( \{r_{l+1 } \} -p ) } Q(p)
\\
&=& \sum_{p=-\infty} ^\infty\Phi^n_l(p)
Q(p),
\end{eqnarray*}
where $\Phi^n_l(p)$ is equal to
\[
\frac{ \max\{ ( \{r_{l'+1}\} -p ) \wedge\{ r_{l+1} \}- ( \{r_l\} -p
) \vee\{r_l'\},0\} }n \mathbf{1}_{[\{r_l\} - \{r_{l'+1}\}, \{
r_{l+1}\} - \{r_{l'} \} ]} (p).
\]
The term $ \Phi^n_l(p)$ is uniformly bounded and converges to $0$ as
$n$ tends to infinity for any fixed $p$ because $l<l'$. Therefore,
taking into account that $ \sum_{p=-\infty} ^\infty Q(p)=Q<\infty$, the
dominated convergence theorem implies that
$ \gamma_n^2 \wt{G}_1$ converges to zero as $n$ tends to infinity.
Similarly, we can show that $ \gamma_n^2 \wt{G}_2$ converges to zero
as $n$ tends to infinity.

In the case $H= \frac{3}4$, since
\begin{eqnarray*}
\gamma_n^2 \wt{G}_1& = & \frac{n^2}{\ln n}
\sum_{k=\{r_l\} }^{\{r_{l+1}\}} \sum
_{k'=\{r_{l'}\} }^{\{r_{l'+1}\}} \int_{\cD_{k', k}} \mu(dv
\,du)\mu( ds \,dt)
\\
&= & \frac{1 }{n \ln n} \sum_{p= \{r_{l } \} - \{r_{l'+1 } \} }^{ \{
r_{l+1 } \} - \{r_{l' } \} }
\sum_{k'=( \{r_{l } \} -p) \vee\{r_{l' } \} }^{
\{r_{l'+1 } \}
\wedge( \{r_{l+1 } \} -p ) } Q(p) ,
\end{eqnarray*}
we have
\begin{eqnarray*}
\gamma_n^2 \wt{G}_1 &\leq&
\frac{1 }{n \ln n} \sum_{p= \{r_{l } \} - \{r_{l'+1 } \} }^{ \{r_{l+1
} \} - \{r_{l' } \} } \sum
_{k'= \{r_{l' } \} }^{
\{r_{l+1 } \} -p } Q(p)
\\
&\leq& \frac{1 }{n \ln n} \sum_{p= -n }^{ 0 }
(p+1) Q(p).
\end{eqnarray*}

Noticing that $Q(p)= O(\frac{1}{|p|})$, we conclude that
$ \gamma_n^2 \wt{G}_1 \leq\frac{C}{\ln n}$.
This shows that
$ \gamma_n^2 \wt{G}_1 $ converges to zero\vspace*{2pt} as $n$ tends to infinity. In
the same way we can show that
$ \gamma_n^2 \wt{G}_2 $ converges to zero.
\end{pf}

The following estimate is needed in the calculation of the fourth
moment of $\Theta_l^n(i, j)$ in Proposition~\ref{prop3}.
%

\begin{lemma}\label{lemma2}
Let $H \in(\frac{1}2,\frac{3}4]$. We have the following estimate:
\begin{eqnarray*}
&&\sum_{ k_1, k_2, k_3, k_4 =0 }^{n-1} \langle
\beta_{{k_1}/{n} } , \beta_{{k_2}/n} \rangle_{\mathcal{H } }
\langle
\beta_{
{k_2 }/n}, \beta_{{k_3 }/n} \rangle_{\mathcal{H } } \langle
\beta_{{k_3 }/n}, \beta_{{k_4 }/n} \rangle_{\mathcal
{H } } \langle
\beta_{{k_1}/n}, \beta_{{k_4 }/n} \rangle_{\mathcal
{H } }\\
&&\qquad \leq{C}
{n^{-2}}\ga_n^{-2} .
\end{eqnarray*}
\end{lemma}
\begin{pf} Since the indices $ k_1,k_2,k_3,k_4$ are
symmetric, it suffices to consider the case $ k_1 \leq k_2 \leq k_3
\leq k_4 $.
By definition of the inner product we have
\begin{eqnarray*}
&& \sum_{ 0\le k_1 \leq k_2 \leq k_3 \leq k_4 \le n-1 } \langle\beta
_{{k_1}/{n} },
\beta_{{k_2}/n} \rangle_{\mathcal{H
} } \langle\beta_{{k_2 }/n},
\beta_{{k_3 }/n} \rangle_{\mathcal
{H } } \langle\beta_{{k_3}/n},
\beta_{{k_4 }/n} \rangle_{\mathcal
{H } } \langle\beta_{{k_1}/n},
\beta_{{k_4 }/n} \rangle_{\mathcal
{H } }
\\
&&\qquad = \frac{T^{8H}}{2^4 n^{8H}} \sum_{k_1 =0}^{n-1}
\sum_{ k_2 =k_1}^{n-1} \sum
_{ k_3 =k_2 }^{n-1} \sum_{ k_4 =k_3}^{n-1}
\bigl(|k_2 -k_1 +1 |^{2H} + |k_2
-k_1 -1|^{2H} \\
&&\hspace*{240pt}{}-2 |k_2 -k_1
|^{2H} \bigr)
\\
&&\hspace*{118pt}\qquad\quad{} \times\bigl(|k_3-k_2+1 |^{2H} +
|k_3-k_2 -1|^{2H}\\
&&\hspace*{252pt}{} -2 |k_3-k_2
|^{2H} \bigr)
\\
&&\hspace*{118pt}\qquad\quad{} \times\bigl(|k_4-k_3+1 |^{2H} +
|k_4-k_3 -1|^{2H} \\
&&\hspace*{252pt}{}-2 |k_4-k_3
|^{2H} \bigr)
\\
&&\hspace*{118pt}\qquad\quad{} \times\bigl(| k_4 - k_1 +1 |^{2H} + |
k_4-k_1 -1|^{2H}\\
&&\hspace*{252pt}{} -2 |k_4
-k_1 |^{2H} \bigr) .
\end{eqnarray*}
Denote $p_i = k_{i+1}-k_i $, $i=1,2,3$. Then the above sum is bounded by
\[
C n^{1-8H} \sum_{p_1, p_2, p_3=1}^{n-1}
p_1^{2H-2} p_2^{2H-2}
p_3^{2H-2} (p_1 +p_2+p_3)^{2H-2},
\]
which is again bounded by
\[
C n^{1-8H} \sum_{p_1, p_2, p_3=1}^{n-1}
p_1^{2H-2} p_2^{2H-2}
p_3^{4H-4} . %
\]
In the case $H \in(\frac{1}2, \frac{3}4)$, the series $\sum_{p_3=1}^{n-1}
p_3^{4H-4}$ is convergent.
When $H=\frac{3}4$, it is bounded by $C\log n$. So the above sum is
bounded by $Cn^{ -4H-1 }$ if $\frac{1}2 <H <\frac{3}4$
and bounded by $ Cn^{-4}{ {\log n}} $ if $H=\frac{3}4$.
The proof is complete.
\end{pf}

The following proposition contains a result on the convergence of the
fourth moment of $\Theta^l_n(i, j)$.
%

\begin{prop}\label{prop3}
The fourth moment of $ \Theta_l^n (i, j)$ and $3\mE(|\Theta_l^n (i, j)
|^2)^2 $ converge to the same limit as $n \rightarrow\infty$.
\end{prop}

\begin{pf}
Applying the integration by parts formula (\ref{ibp}) yields
\begin{eqnarray*}
&&\mE\bigl[ \Theta_l^n (i, j ) ^4 \bigr] \\
&&\qquad=
\alpha_H^2 \ga_n \sum
_{k=\{r_l\} }^{\{r_{l+1}\}} \int_{\cD_{k}} \mE
\bigl[ D^i_u D^j_t \bigl[
\Theta^n_{l } (i, j )^3 \bigr] \bigr] \mu(dv
\,du)\mu( ds \,dt)
\\
&&\qquad= \alpha_H^2 \ga_n \sum
_{k=\{r_l\} }^{\{r_{l+1}\}} \int_{\cD_{k}} \mE
\bigl[ \bigl\{ 3 \Theta^n_{l } (i, j )^2
D^i_u D^j_t \bigl[
\Theta^n_{l } ( i, j ) \bigr]
 \\
 &&\qquad\hspace*{93pt}{}+ 6 \Theta^n_{l } (i, j ) D^j_t
\bigl[ \Theta^n_{l } ( i, j ) \bigr]\\
&&\hspace*{165pt}{}\times D^i_u
\bigl[ \Theta^n_{l } ( i, j ) \bigr] \bigr\} \bigr] \mu(dv
\,du )\mu(ds \,dt)
\\
&& \qquad:=\overline{G}_1 +\overline{G}_2.
\end{eqnarray*}
Since $ D^i_u D^j_t
[ \Theta^n_{l } ( i, j ) ] $ is deterministic,
it is easy to see that $\overline{G}_1 = 3 \mE(|\Theta_l^n(i,\break  j )
|^2 )^2$.
We have shown the convergence of $\mE(|\Theta^l_n(i, j ) |^2 )$ in
Proposition~\ref{prop 1}.
It remains to show that $\overline{G}_2 \rightarrow0$ as $ n
\rightarrow\infty$.

Applying again the integration by parts formula (\ref{ibp}) yields
\begin{eqnarray*}
\overline{G}_2& =& 6 \alpha_H^4
\ga_n^2 \sum_{k, k'=\{r_l\} }^{\{r_{l+1}\}}
\int_{\cD_{k} \times\cD_{ k' } } D^i_{u'}
D^j_{t'} \bigl\{ D^j_t \bigl[
\Theta^n_{l } ( i, j ) \bigr] D^i_u
\bigl[ \Theta^n_{l } ( i, j ) \bigr] \bigr\}
\\
&&\hspace*{102pt}{}\times \mu\bigl(dv' \,du' \bigr)\mu
\bigl(ds' \,dt'\bigr) \mu(dv \,du )\mu(ds \,dt).
\end{eqnarray*}
Using equation (\ref{eqn 48}) we can derive the inequalities
\begin{eqnarray*}
\overline{G}_2 &\le& 6 \alpha_H^4
\ga_n^4 \sum_{k, k', h, h'=\{ r_l \} }^{ \{ r_{l+1} \} }
\int_{\cD_{k} \times\cD_{ k' } } \bigl\{ \beta_{ { h}/{ n} } (t)
\beta_{{h }/{ n} } \bigl(t' \bigr) \beta_{
{ h'}/{ n} } (u)
\beta_{{h '}/{ n} } \bigl(u'\bigr)
\\
&&\hspace*{122pt}{} + \beta_{ {h}/ { n} } (t) \beta_{ {h}/{ n} } \bigl(u'
\bigr) \beta_{
{ h' }/{n} } (u) \beta_{ {h'}/{ n} } \bigl(t'
\bigr) \bigr\} \\
&&\hspace*{122pt}{}\times \mu\bigl( dv'd u' \bigr)\mu\bigl( d
s' \,d t'\bigr) \mu( dv\, d u)\mu( d s\, d t)
\\
&\le& 12 \alpha_H^4 \ga_n^4
\sum_{k, k', h, h'=\{ r_l \} }^{ \{ r_{l+1} \} } \langle
\beta_{ { h }/{ n} }, \beta_{ { {k} }/{ n} } \rangle_{\cH} \langle
\beta_{ { h '}/{ n} }, \beta_{ { {k} }/{ n} } \rangle_{\cH} \langle
\beta_{ { h }/{ n} }, \beta_{ { {k}' }/{ n} } \rangle_{\cH} \\
&&\hspace*{88pt}{}\times\langle
\beta_{ { h' }/{ n} }, \beta_{ { {k}' }/{ n} } \rangle_{\cH}.
\end{eqnarray*}
The convergence of $\overline{G}_2$ to zero now follows from
Lemma~\ref{lemma2}.
\end{pf}

We can now establish a central limit theorem for weighted sums based on
the previous proposition. Recall that $ {\zeta}^{ i,j }_{k,n} = \intl
_{t_k }^{t_{k+1} }( B^{i }_s - B^{i}_{t_k}) \delta B^{j }_s$, $k=0,
\ldots, n-1$ and $ {\zeta}^{ i,j }_{n,n} = 0$.
%

\begin{prop}\label{prop 3}
Let $f=\{f_t, t\in[0,T]\}$ be a stochastic process with values on the
space of $d\times d$ matrices and with H\"older continuous trajectories
of index greater than $\frac{1}2$. Set, for $i,j =1,\ldots,m$,
\[
{ \Psi}_n^{ i, j } (t ) = \sum_{k=0 }^{ \{t\}}
f_{t_k}^{i, j } {\zeta}^{ i,j }_{k,n}.
\]
Then, the following stable convergence in the space $D([0,T])$ holds as
$n$ tends to infinity:
\[
\bigl\{ \ga_n \Psi_n (t), t\in[0, T] \bigr\} \rightarrow
\biggl\{ \biggl( \int_0^t
f^{i, j }_s \,dW^{ij}_s
\biggr)_{1 \leq i, j \leq m}, t\in[0, T] \biggr\},
\]
where $W$ is a matrix-valued Brownian motion independent of $B$ with
the covariance introduced in Section~\ref{sec2.3}.
\end{prop}

\begin{pf}
This proposition is an immediate consequence of the central limit
result for weighted random sums proved in \cite{CNP}.
In fact, the process $\Psi^{i,j}_n (t)$ satisfies the required
conditions due to Proposition~\ref{prop4.2} and the estimate (\ref{eqn4.1}).
\end{pf}

\section{CLT for the modified Euler scheme in the case \texorpdfstring{$H \in(\frac{1}2,\frac{3}4]$}{Hin(1/2,3/4]}}\label{eqn 43}
The following central limit type result shows that in the case $H
\in(\frac{1}2, \frac{3}4 ]$, the process $\gamma_n (X - X^n )$ converges
stably to the solution of a linear stochastic differential equation
driven by a
matrix-valued Brownian motion independent of $B$ as $n$ tends to infinity.

%
\begin{theorem}\label{theorem 1}
Let $H \in(\frac{1}2, \frac{3}4 ]$, and let $X$, $X^n$ be the solutions
of the SDE (\ref{e.1.1}) and recursive
scheme (\ref{e.1.6}), respectively. Let $W=\{W_t, t\in[0,T]\} $ be the
matrix-valued Brownian motion introduced in Section~\ref{sec2.3}.
Assume $\sigma\in
C^5_b(\mathbb{R}^d; \mathbb{R}^{d \times m})$ and $b \in
C^4_b(\mathbb{R}^d; \mathbb{R}^{d })$. Then the following stable
convergence in the space $C([0, T])$ holds as $ n $ tends to infinity:
%
%
\begin{equation}
\bigl\{\ga_n \bigl(X_t-X_t^n
\bigr), t\in[0, T] \bigr\} \rightarrow\bigl\{ U_t, t\in[0, T]
\bigr
\}, \label{stable}
\end{equation}
where $\{U_t, t \in[0, T]\}$ is the solution of the linear
$d$-dimensional SDE
%
%
\begin{eqnarray}
\label{eqn 45} 
U_t &=& \intl^t_0
\nabla b (X_s) U_s \,ds + \sum
_{j=1}^m \intl^t_0 \nabla
\sigma^j(X_s) U_s \,dB_s^{j}
\nonumber
\\[-8pt]
\\[-8pt]
\nonumber
&&{}+ \sum_{i,j=1}^m \intl^{t}_0
\bigl(\nabla\sigma^j \sigma^i\bigr) (X_s)
\,dW^{ij}_s.
\end{eqnarray}
\end{theorem}

%
\begin{remark}
It follows from \cite{Kurtz} that when $B$ is replaced by a
standard Brownian motion,
the process $\sqrt{n} (X-X^n)$ converges in law to the unique
solution of the $d$-dimensional SDE
%
%
\begin{eqnarray}
\label{eqn 44} dU _t &=& \nabla b (X_t) U_s
\,dt + \sum_{j=1}^m \nabla
\sigma^{ j} (X_t ) U_t \,dB^j_t
\nonumber
\\[-8pt]
\\[-8pt]
\nonumber
&&{}+ \sqrt{\frac{T}{
2} } \sum_{j,i=1}^m
\bigl( \nabla\sigma^{ j} \sigma^i \bigr) (X_t
) \,d W^{ij} _t
\end{eqnarray}
with $U_0=0$.
Here $W^{ij}$, $i,j=1, \ldots, m$ are independent one-dimensional
Brownian motions, independent of $B$. To compare our Theorem~\ref
{theorem 1} with this result, we let the Hurst parameter $H$
converge to $\frac{1}2$. Then the constant $R$ will converge to $0$,
and $ \frac{\alpha_H}{\sqrt{T} } {\sqrt{Q-R } } $ converges to
$\sqrt{\frac{T}2}$. This formally recovers equation (\ref{eqn 44}).
\end{remark}

%
\begin{remark}
The process $U$ defined in (\ref{eqn 45}) is given by
%
%
\begin{equation}
U_t = \sum_{i,j=1}^m
\intl^{t}_0\Lambda_{t} \Gamma_s
\bigl( \nabla\sigma^j\sigma^i\bigr) (X_s)
\,dW^{ij}_s,\qquad t\in[0, T], \label{ut}
\end{equation}
where we recall that $\Lambda$ is defined in (\ref{eqn 46}) and
$\Gamma$
is its inverse.
\end{remark}

\begin{pf*}{Proof of Theorem~\ref{theorem 1}}
Recall that $Y_t=X_t-X^n_t$.
We would like to show that the process $\{\gamma_n Y_t, B_t, t \in[0, T]
\}$ converges weakly in $C([0,T]; \mathbb{R}^{d+m})$ to $\{U_t, B_t,
t\in[0, T] \}$.
To do this, it suffices to prove the following:
\begin{longlist}[(ii)]
\item[(i)] convergence of the finite dimensional distributions of $ \{
\gamma_n Y_t, B_t, t \in[0, T]
\} $;
\item[(ii)] tightness of the process $\{\gamma_n Y_t, B_t, t \in
[0, T]
\}$.
\end{longlist}

We first show (i).
Recall the decomposition of $Y_t$ given in (\ref{eqn 36}) and (\ref{eqn20}),
and recall the estimates obtained for each term in the decomposition of $Y_t$.
Since the
other terms converge to zero in $L^p$ for $p\ge1$, from the Slutsky theorem
it suffices to consider the convergence of the finite dimensional
distributions of $\{\gamma_n \sum_{j=1}^m E_{2,j}(t), B_t, t \in[0,
T]\}$, where $E_{2,j}$ is defined in Theorem~\ref{thmm4.1} step 3.
Set
%
%
\begin{equation}
\label{6.5e} F_s^{i,j}: = \Lambda^n_t
\Gamma^n_s \bigl( \nabla\sigma^j
\sigma^i\bigr) \bigl(X^n_s\bigr) -
\Lambda_t \Gamma_{s} \bigl( \nabla\sigma^j
\sigma^i\bigr) (X_{s} ).
\end{equation}
It follows from Lemma~\ref{lem3.5} and Remark~\ref{remark4.3} that
\[
\sup_{r,s,t \in[0, T] } \bigl( \bigl\llVert F^{i,j}_{t}
\bigr\rrVert_p \vee\bigl\llVert D_{s}F^{i,j}_{t}
\bigr\rrVert_p \vee\bigl\llVert D_{r}
D_{s}F^{i,j}_{t} \bigr\rrVert_p
\bigr) \leq Cn^{1-2\beta}.
\]
Denote
%
%
\begin{equation}
\label{e6.5} \wt{E}_{2,j} (t) = \Lambda_t \sum
_{i=1}^m \sum_{k=0}^{
{ \lfloor{nt}/{T} \rfloor}
}
\Gamma_{t_k} \bigl( \nabla\sigma^j \sigma^i
\bigr) (X_{t_k} ) \int_{t_k}^{t_{k+1} } \int
_{t_k}^s \delta B^i_u
\delta B^j_s,
\end{equation}
for $t \in[0,T)$, and $\wt{E}_{2,j}(T) = \wt{E}_{2,j}(T-)$.
Then applying Lemma~\ref{lem11.4} (\ref{11.d}) with
$F^{i,j}$ defined by (\ref{6.5e}), we obtain that
\[
\gamma_n \bigl\llVert E_{2,j}(t) - \wt{E}_{2,j}(t)
\bigr\rrVert_p \leq C \gamma_n n^{-H}
n^{1-2\beta},
\]
which converges to zero as $n\rightarrow\infty$ since $\beta$ can be
taken as close as possible to~$H$. By Slutsky's theorem again, it
suffices to consider the convergence of the finite dimensional
distributions of
%
%
\begin{equation}
\label{e6.6} \Biggl\{\gamma_n \sum_{j=1}^m
\wt{E}_{2,j}(t), B_t, t \in[0, T] \Biggr\}.
\end{equation}

Applying Proposition~\ref{prop 3} to the family of processes $f^{i,j
}_t = \Gamma_{t } ( \nabla\sigma^j \sigma^i )(X_{t } ) $, we obtain
the convergence of the finite dimensional distributions of
\[
\Biggl\{\gamma_n \sum_{j=1}^m
\Gamma_t \wt{E}_{2,j}(t), B_t, t \in[0, T]
\Biggr\}
\]
to those of
$\{\Gamma_t U_t, B_t, t\in[0,T]\}$. This implies the convergence of
the finite dimensional distributions of
\[
\Biggl\{\gamma_n \sum_{j=1}^m
\wt{E}_{2,j}(t), B_t, t \in[0, T]\Biggr\}
\]
to those of
$\{U_t, B_t, t\in[0,T]\}$.

To show (ii), we prove the following tightness condition:
%
%
\begin{equation}
\label{eqn tight} \sup_{n \geq1} \mathbb{E } \bigl(\bigl\vert
\gamma_n \bigl(X_t-X_t^n\bigr) -
\gamma_n \bigl(X_s-X_s^n\bigr)
\bigr\vert^{4} \bigr) \leq C (t-s)^{2}.
\end{equation}
Taking into account (\ref{eqn 36}) and (\ref{eqn20}), we only need to
show the above
inequality for $ \gamma_n I_{11} $, $ \gamma_n I_{12,j} $, $ \gamma_n
I_{13} $,
$ \gamma_n I_{2,j} $, $ \gamma_nI_{4,j}$, $\gamma_n E_{1,j}$,
$\gamma
_n E_{2,j}$ and $\gamma_n E_{3,j}$.
The tightness
for the terms $\gamma_nI_{11}$, $\gamma_nI_{13}$ and $\gamma_n
E_{3,j}$ is clear.
Now we consider the tightness of the term $I_{2,j}$. We write
\begin{eqnarray*}
I_{2, j} (t)-I_{2, j} (s)&=& \bigl( \Lambda^n_t
- \Lambda^n_s \bigr) \intl^t_0
\Gamma^n_s b^j_2 (s) \bigl(s-
\eta(s)\bigr) \,dB^j_s
\\
&&{} + \intl^t_s \Lambda^n_s
\Gamma^n_u b^j_2 (u) \bigl(u-
\eta(u)\bigr) \,dB^j_u.
\end{eqnarray*}
Then it follows from Lemma~\ref{lem11.3} (\ref{11.a}) that
\begin{eqnarray*}
\mathbb{E} \bigl( \bigl\vert\gamma_n \bigl(
\Lambda^n_t - \Lambda^n_s \bigr)
\intl^t_0 \Gamma^n_s
b^j_2 (s) \bigl(s- \eta(s)\bigr) \,d B^j_s
\bigr\vert^4 \bigr) &\leq& C (t-s)^{4 \beta} \bigl(
\mathbb{E} \bigl\| \Lambda^n \bigr\|_{\beta}^8 \bigr)
^{{1}/2}
\\
&\leq& C(t-s)^{4\beta}. 
\end{eqnarray*}
Lemma~\ref{lem11.3} (\ref{11.a}) also implies that the fourth moment of
the second term is bounded by
$
C(t-s)^{4H}$.
The tightness for $\gamma_n I_{12,j}$, $\gamma_n I_{4,j}$, $\gamma_n
E_{1,j}$, $\gamma_n E_{2,j} $ can be
obtained in a similar way by applying the estimates (\ref{11.b}) and
(\ref{11.a}) from Lemma~\ref{lem11.3},
(\ref{e 11.10}) from Lemma~\ref{lem11.5}, and
(\ref{11.c}) from Lemma~\ref{lem11.4}, respectively.
\end{pf*}

\section{A limit theorem in $L^p$ for weighted sums}\label{section 7}
Following the methodology used in \cite{CNP}, we can show the
following limit result for random weighted sums.
The proof uses the techniques of fractional calculus and the classical
decompositions in large and small blocks.

Consider a double sequence of random variables $\zeta= \{ \zeta_{k, n}
, n \in\mathbb{N}, k=0, 1, \ldots,  n\}$, and for each $t \in[0, T]$,
we denote
%
%
\begin{equation}
\label{e7.1} g_n(t): = \sum_{k=0}^{ \lfloor{nt}/T \rfloor}
\zeta_{k, n}.
\end{equation}
%

\begin{prop}\label{prop7.1}
Fix $ \lambda> 1- \beta$, where $ 0 < \beta<1$.
Let $p\geq1$ and $p', q' > 1$ such that $\frac{1}{p'} +\frac{1}{q'}
=1$ and $pp'>\frac{1}{ \beta}$, $pq' > \frac{1}{\lambda}$. Let $g_n $
be the sequence of processes defined in (\ref{e7.1}). Suppose that the
following conditions hold true:
\begin{longlist}[(ii)]
\item[(i)] for each $t \in[0, T]$, $g_n(t ) \rightarrow z(t)$ in $L^{pq'}$;

\item[(ii)] for any $j,k =0, 1, \ldots, n $ we have
\[
\mE\bigl(\bigl| g_n( kT/n ) - g_n ( jT/n)\bigr |^{pq'}
\bigr) \leq C \bigl( |k-j|/n\bigr)^{\lambda
pq' }.
\]
Let $f=\{ f(t), t \in[0, T]\}$ be a process such that $\mE(\|f\|
_{\beta}^{pp'} ) \leq C$ and $\mE(|f( 0 )|^{pp'}) \leq C$.
Then for each $t\in[0, T]$,
%
%
\begin{equation}
\label{eqn7.3} F^n(t):= \sum_{k=0}^{ \lfloor{nt}/T \rfloor}f(t_k)
\zeta_{k, n} \rightarrow\int_0^t
f(s) \,dz(s) \qquad\mbox{in } L^p \mbox{ as } n \rightarrow\infty.
\end{equation}
\end{longlist}
\end{prop}

%
\begin{remark}
The integral $\int_0^t f(s) \,dz(s)$ is interpreted as a Young integral
in the sense of Proposition~\ref{prop.1},
which is well defined because $f$ and $z$, as functions on $[0,T]$ with
values in $L^{pp'}$ and $L^{pq'}$,
are H\"older continuous [conditions (i) and (ii) together imply the H\"
older continuity of $z$] of
order $\beta$ and $\lambda$, respectively. Recall that the H\"older
continuity of a function with values
in $L^p$ is defined in (\ref{eqn2.4}).
\end{remark}
%
%
\begin{remark} Convergence (\ref{eqn7.3}) still holds true if the
condition \break $\mE(\|f\|_{\beta}^{pp'} )
\leq C$ is weakened by assuming that $f$ is H\"older continuous of
order $\beta$ in $L^{pp'}$.
The proof will be similar to that of Proposition~\ref{prop7.1}.
\end{remark}
\begin{pf*}{Proof of Proposition \ref{prop7.1}}
Given two natural numbers $m <n$ we
consider the associated partitions of the interval $[0,T]$ given by
$t_k= \frac{kT}n$, $k=0,1,\ldots, n$ and $u_l= \frac{lT}m$,
$l=0,1,\ldots, m$.
Then we have the decomposition
%
%
\begin{equation}
F^n(t) = \sum_{l=0}^{ \lfloor{m t}/T \rfloor}
f(u_l) \sum_{k \in I_m(l)} \zeta_{k, n}
+ \sum_{l=0}^{\lfloor{m t }/T \rfloor} \sum
_{k\in I_m(l)} \bigl[f(t_k) - f(u_l) \bigr]
\zeta_{k, n},
\end{equation}
where $I_m(l):= \{ k\dvtx0\le k \le\lfloor\frac{nt} T \rfloor, t_k
\in
[ u_{l}, u_{l+1} )\}$.

Because of condition (i) and the assumption that $\mE(|f(t)|^{pp'})
\leq C$ for all $t\in[0, T]$,
the first term on the right-hand side of the above expression converges
in $L^{p}$, as $n$ tends to infinity, to
\[
\sum_{l=0}^{ \lfloor{mt}/ T \rfloor} f(u_l)
\bigl[ z(u_{l+1})- z(u_l) \bigr].
\]
Applying Proposition~\ref{prop.1} to $f$ and $z$ we obtain that the
above Riemann--Stieltjes sum converges to
the Young integral $\int_0^t f(s) \,dz(s)$ in $L^p$ as $m$ tends to infinity.
To show convergence (\ref{eqn7.3}) it suffices to show that
%
%
\begin{equation}
\label{equ3} \lim_{m\rightarrow\infty} \sup_{n \in\mathbb{N} }
\mE
\Biggl( \Biggl\vert\sum_{l=0}^{ \lfloor{mt}/ T \rfloor} \sum
_{k\in I_m(l)} \bigl[f(t_k) -
f(u_l) \bigr] \zeta_{k, n} \Biggr\vert^p
\Biggr)=0.
\end{equation}

Notice that $k$ belongs to $I_m(l)$ if and only if $u_l \le t_k <
\varepsilon(u_{l+1} ) $ and $t_k \le\eta(t)$. Recall that
$\varepsilon(u)=
t_{k+1} $ if
$t_k < u \le t_{k+1} $ and $\eta(u)= t_{k} $ if
$t_k \le u < t_{k+1} $. As a consequence, we can write
\begin{eqnarray*}
&&\sum_{l=0}^{ \lfloor{mt} /T \rfloor} \sum
_{k\in I_m(l)} \bigl[f(t_k) - f(u_l) \bigr]
\zeta_{k, n}\\
&&\qquad= \sum_{l=0}^{ \lfloor{mt} /T \rfloor} \int
_{ ( {a_l}, {b_l} ) } \bigl[f(s) - f({a_l} ) \bigr] \,d
g_n (s),
\end{eqnarray*}
where ${a_l}=
u_l$ and ${b_l} = \varepsilon(u_{l+1} ) \wedge(\eta(t) +\frac{T}n
) $.
By the fractional integration by parts formula,
%
%
\begin{eqnarray}
\label{e7.3} 
&& \int_{ ( {a_l}, {b_l} ) } \bigl[f(s) -
f({a_l} ) \bigr] \,d g_n (s)
\nonumber
\\[-8pt]
\\[-8pt]
\nonumber
&&\qquad = (-1)^{\alpha} \int_{a_l}^{{b_l} }
D^{\alpha}_{{a_l}+} \bigl[f(s) - f({a_l})\bigr]
D^{1-\alpha}_{ {b_l}-} \bigl[g_n(s) - g_n({b_l}
-)\bigr] \,ds,
\end{eqnarray}
where we take $\alpha\in(1-\lambda, \beta)$.
By (\ref{e.2.1}), it is easy to show that
%
%
\begin{eqnarray}
\label{eq6.6} 
\bigl| D^{\alpha}_{{a_l}+} \bigl[f(s) -
f({a_l})\bigr] \bigr| &\leq& \frac{1}{\Gamma(1-\alpha)} \frac{\beta
}{\beta-\alpha} \|f\|
_{\beta} (s-{a_l})^{\beta-\alpha}
\nonumber
\\[-8pt]
\\[-8pt]
\nonumber
& \leq& C\|f\|_{\beta} m^{\alpha-\beta}.
\end{eqnarray}
On the other hand, by (\ref{e.2.2}) we have
%
%
\begin{eqnarray}
\label{eq6.7} 
&& \bigl\vert D^{1-\alpha}_{{b_l}-}
\bigl[g_n(s) - g_n({b_l}-)\bigr] \bigr\vert
\nonumber
\\[-8pt]
\\[-8pt]
\nonumber
&&\qquad= \frac{1 }{\Gamma(\alpha) } \biggl\vert\frac{g_n(s) -
g_n({b_l}-) }{
({b_l}-s)^{1-\alpha} } +(1-\alpha) \int
_s^{b_l} \frac{g_n(s) - g_n(u)
}{ (u-s)^{2-\alpha}} \,du\biggr\vert.
\end{eqnarray}
We can calculate the integral in the above equation explicitly.
%
%
\begin{eqnarray}
\label{eq6.8} 
&& \int_s^{b_l}
\frac{g_n(s) - g_n(u) }{ (u-s)^{2-\alpha}} \,du\nonumber
\\
&&\qquad= \int_{\ep(s)}^{b_l} \frac{g_n(s) - g_n(u) }{ (u-s)^{2-\alpha}}
\,du
\nonumber
\\[-8pt]
\\[-8pt]
\nonumber
&&\qquad= \sum_{k \dvtx t_k \in[ \ep(s), {b_l} )} \bigl[g_n(s) -
g_n (t_k ) \bigr] \int_{t_k}^{t_{k+1 } }
(u-s)^{\alpha-2} \,du
\\
&&\qquad= \sum_{k \dvtx t_k \in[ \ep(s), {b_l} )} \bigl[g_n(s) -
g_n (t_k ) \bigr] \frac
{1}{1-\alpha}
\bigl[(t_k - s)^{\alpha-1} - (t_{k+1} -
s)^{\alpha-1} \bigr].
\nonumber
\end{eqnarray}
Substituting (\ref{eq6.6}), (\ref{eq6.7}) and (\ref{eq6.8}) into
(\ref{e7.3}), we obtain
\begin{eqnarray*}
&& \biggl\vert\int_{ ( {a_l}, {b_l} ) } \bigl[f(s) -
f({a_l} ) \bigr] \,d g_n (s)\biggr\vert
\\
&&\qquad\leq C\|f\|_{\beta} m^{\alpha-\beta} \int_{ {a_l} }^{b_l}
\bigl\vert D^{1-\alpha}_{{b_l}-} \bigl[g_n(s) -
g_n({b_l}-)\bigr] \bigr\vert\,ds
\\
&&\qquad\leq C\|f\|_{\beta} m^{\alpha-\beta} \sum
_{ k \dvtx t_{ k } \in[ \eta({a_l}), {b_l} ) } \int_{t_{k } }^
{t_{k+1 } } \bigl
\vert D^{1-\alpha}_{{b_l}-} \bigl[g_n(s) -
g_n({b_l}-)\bigr] \bigr\vert\,ds
\\
&&\qquad \leq C\|f\|_{\beta} m^{\alpha-\beta} \sum
_{ k \dvtx t_{k } \in[ \eta({a_l}), {b_l} ) }\bigl | g_n( t_{k } ) -
g_n({b_l}-)\bigr | \int_{t_{k } }^ {t_{k+1 } }
({b_l}-s)^{\alpha-1} \,ds
\\
&&\qquad\quad{} + C\|f\|_{\beta} m ^{\alpha-\beta} \sum_{k, j \dvtx\eta({a_l})
\leq t_{k } < t_j < {b_l} }
\bigl| g_n( t_{k } ) - g_n(t_j ) \bigr|
\\
&&\qquad\quad{}\times\int_{t_{k } }^ {t_{k+1 } } \bigl[(t_j -
s)^{\alpha-1} - (t_{j+1} - s)^{\alpha-1} \bigr] \,ds. 
\end{eqnarray*}
We denote the first term in the right-hand side of the above expression
by $ A_{1,l}$ and the second one by $ A_{2, l}$.

Applying the Minkowski inequality, we see that
the quantity
%
%
\begin{equation}
\label{e 7.9} \mE\Biggl( \Biggl\vert\sum_{l=0}^{ \lfloor{mt}/ T
\rfloor}
A_{1,l} \Biggr\vert^p \Biggr)^{{1}/p}
\end{equation}
is less than
\[
C m^{\alpha-\beta} \Biggl| \sum_{l=0}^{ \lfloor{mt}/ T \rfloor}
\sum_{k \dvtx t_{k } \in[ \eta({a_l}), {b_l} ) } \mE\bigl( \|f\|
_{\beta}^p
\bigl| g_n( t_{k } ) - g_n({b_l}-)
\bigr|^p\bigr)^{{1}/p} \int_{t_{k } }^ {t_{k+1 } }
({b_l}-s)^{\alpha-1} \,ds \Biggr|,
\]
so by applying the H\"older inequality,
condition (ii) and the assumption $\mE(\|f\|_{\beta}^{pp'} ) \leq C$
to the above, we can show that quantity (\ref{e 7.9}) is less than
%
%
\begin{equation}
\label{e 7.10} C m^{\alpha-\beta} \Biggl\vert\sum
_{l=0}^{ \lfloor{mt}/ T \rfloor} \sum_{k \dvtx t_{k } \in[ \eta
({a_l}), {b_l} ) }
({b_l}- t_{k } )^{\lambda} \int_{t_{k } }^ {t_{k+1 } }
({b_l}-s)^{\alpha-1} \,ds \Biggr\vert.
\end{equation}
Since
\begin{eqnarray*}
&& \sum_{ k \dvtx t_{k } \in[ \eta({a_l}), {b_l} ) } ({b_l}- t_{k}
)^{\lambda} \int_{t_{k} }^ {t_{k+1 } }
({b_l}-s)^{\alpha-1} \,ds
\\
&&\qquad= \frac{1}{\alpha} \biggl(\frac{T}n \biggr)^{ \lambda+ \alpha}
+ \sum
_{ k \dvtx t_{k} \in[ \eta({a_l}), {b_l}-{T}/{n}
) } ({b_l}- t_{k }
)^{\lambda} \int_{t_{k } }^ {t_{k+1 } }
({b_l}-s)^{\alpha-1} \,ds
\\
&&\qquad \leq\frac{1}{\alpha} \biggl(\frac{T}n \biggr)^{ \lambda+
\alpha} +
\frac{T}{n} \sum_{k \dvtx t_{k } \in[ \eta({a_l}), {b_l}-{T}/{n}
) } ({b_l}-
t_{k } )^{\lambda} ({b_l}-t_{k+1 })^{\alpha-1}
\\
&&\qquad \leq\frac{1}{\alpha} \biggl(\frac{T}n \biggr)^{ \lambda+
\alpha} +C
\frac{T}{n} \frac{n}{m} m^{ -\alpha+1-\lambda}
\\
&&\qquad \leq C m^{-\alpha- \lambda},
\end{eqnarray*}
where in the second inequality we used the assumption that $\alpha>
1-\lambda$ and the fact that the number of partition points $\{t_k,
k=0, 1,\ldots, n\}$ in $ [ \eta({a_l}), {b_l}-\frac{T}{n}
) $
is bounded by $\frac{n}{m}$, the estimate (\ref{e 7.10}) of (\ref{e 7.9})
implies that
%
%
\begin{equation}
\label{e7.5} \mE\Biggl( \Biggl\vert\sum_{l=0}^{ \lfloor{mt}/T
\rfloor}
A_{1,l} \Biggr\vert^p \Biggr)^{{1}/p} \leq C
m^{\alpha-\beta} \sum_{l=0}^{\lfloor{mt}/T \rfloor}
m^{-\alpha-\lambda} \leq C m^{ 1-\beta-\lambda} \rightarrow0
\end{equation}
as $m$ tends to $ \infty$.

Using an argument similar to the estimate of quantity (\ref{e 7.9}), it
can be shown that
the quantity
\[
\mE\Biggl( \Biggl\vert\sum_{l=0}^{ \lfloor{mt}/T \rfloor}
A_{2,l } \Biggr\vert^p \Biggr)^{{1}/p}
\]
is less than
\begin{eqnarray*}
\label{eq7.12} &&C\Biggl\vert m^{\alpha-\beta} \sum
_{l=0}^{ \lfloor{mt}/T \rfloor} \sum_{ k, j \dvtx\eta({a_l}) \leq
t_{k } < t_j < {b_l} }
| t_{k } - t_j |^{\lambda}
\nonumber
\\[-8pt]
\\[-8pt]
\nonumber
&&\hspace*{142pt}{}\times\int_{t_{k } }^ {t_{k+1 } }
\bigl[(t_j - s)^{\alpha-1} - (t_{j+1} -
s)^{\alpha-1} \bigr] \,ds \Biggr\vert.
\end{eqnarray*}
The summand in the above can be estimated as follows:
\begin{eqnarray*}
&& \sum_{ k, j \dvtx\eta({a_l}) \leq t_{k } <
t_j < {b_l} } | t_{k } - t_j
|^{\lambda} \int_{t_{k } }^ {t_{k+1 } }
\bigl[(t_j - s)^{\alpha-1} - (t_{j+1} -
s)^{\alpha-1} \bigr] \,ds
\\
&&\qquad\leq C\frac{n}{m} \biggl( \frac{T}n \biggr)^{\alpha+ \lambda}\\
&&\qquad\quad{} +
\sum_{ k, j \dvtx\eta(a_l) \leq t_{k +1} < t_j < {b_l} } | t_{k } - t_j
|^{\lambda} \int_{t_{k } }^ {t_{k+1 } }
\bigl[(t_j - s)^{\alpha-1} - (t_{j+1} -
s)^{\alpha-1} \bigr] \,ds
\\
&&\qquad\leq C\frac{n}{m} \biggl( \frac{T}n \biggr)^{\alpha+ \lambda} + C
\biggl(\frac{T}n \biggr)^2 \sum_{k, j \dvtx\eta({a_l}) \leq t_{k
+1} < t_j < {b_l} }
| t_{k+1 } - t_j |^{\lambda} (t_{j } -
t_{k+1 } )^{\alpha-2}
\\
&&\qquad\leq C\frac{n}{m} \biggl( \frac{T}n \biggr)^{\alpha+ \lambda} + C
n^{-2 } n^{2-\lambda- \alpha} \sum_{ k, j \dvtx\eta(a_l ) \leq
t_{k +1} < t_j < {b_l} } (j-k
-1)^{\alpha-2+\lambda}
\\
&&\qquad\leq C\frac{n}{m} \biggl( \frac{T}n \biggr)^{\alpha+\lambda} + C
n^{-2 } n^{2-\lambda- \alpha} \frac{n}{m} \sum
_{p=2 }^{n/m} (p -1)^{\alpha-2+\lambda}
\\
&&\qquad\leq C m^{-\alpha-\lambda}.
\end{eqnarray*}
Therefore, we have
\[
\mE\Biggl( \Biggl\vert\sum_{l=0}^{ \lfloor{mt}/T \rfloor}
A_{2,l} \Biggr\vert^p \Biggr)^{{1}/p} = C
m^{ \alpha- \beta} \Biggl\vert\sum_{l=0}^{ \lfloor{mt}/T \rfloor}
m^{-\al-\lambda} \Biggr\vert\rightarrow0 \qquad\mbox{as } m
\rightarrow\infty.
\]
The above convergence and equality (\ref{e7.5}) together imply
convergence (\ref{equ3}). The proof is now complete.
\end{pf*}

This result has the following two consequences.
%

\begin{cor} \label{prop 5} Let
$B=\{B_t, t\in[0,T]\}$ be an $m$-dimensional fBm with Hurst parameter $H>3/4$.
Define
\[
{\zeta}_{k, n}^{ij} = n\intl_{t_k}^{t_{k+1}}
\bigl(B^i_s - B^i_{\eta(s)}\bigr)
\delta B_s^j,
\]
for $i, j = 1, \ldots, m$ and $k=0,\ldots, n-1$, where we recall that
$t_k =\frac{kT}n$.
Set also ${\zeta}_{n, n}=0$. Let $\lambda= \frac{1}2$, and $\beta$,
$p$, $p'$, $q'$, $f$ satisfy the assumptions in
Proposition~\ref{prop7.1}. Then
\[
\sum_{k=0}^{ \lfloor{nt}/ T \rfloor} f(t_k) {
\zeta}_{k, n} ^{ij} \rightarrow\intl^t_0
f(s) \,dZ^{ij}_s \qquad\mbox{in } L^p,
\]
where $Z^{ij} $ is the generalized Rosenblatt process defined in
Section~\ref{section 2.3}.
\end{cor}

\begin{pf} To prove the corollary, it suffices
to show that the conditions in Proposition~\ref{prop7.1} are all
satisfied here.
We have shown in Section~\ref{section 2.3} the $L^2$ convergence of
$g_n(t) = \sum_{k=0}^{ \lfloor{nt} /T \rfloor} {\zeta}_{k,
n}^{ij} $ to $Z^{ij}_t$.
This convergence also holds in $L^p$ due to the equivalence of all the
$L^p$-norms in a finite Wiener chaos.
Applying (\ref{11.c}) in Lemma~\ref{lem11.4} with $F \equiv1$ and
taking into account that $\gamma_n =n$
when $H > \frac{3}4$, we obtain condition (ii) in Proposition~\ref
{prop7.1} with $\lambda=\frac{1}2$.
\end{pf}

The following result will also be useful later.
%

\begin{cor}\label{prop 6} Let $B=\{B_t, t\in[0,T]\}$ be
one-dimensional fBm with Hurst parameter $H \in(\frac{1}2, 1)$.
Define
%
%
\begin{equation}
\label{7.13e} {\zeta}_{k, n} = \intl_{t_k}^{t_{k+1}}
\bigl(s- \eta(s)\bigr) \,dB_s ,
\end{equation}
for $k=0,\ldots, n-1$. Set also ${\zeta}_{n, n}=0$.
Let $\lambda= H$, and $\beta$, $p$, $p'$, $q'$, $f$ satisfy the
assumptions in
Proposition~\ref{prop7.1}.
Then for each $t\in[0, T]$,
\[
n \sum_{k=0}^{ \lfloor{nt}/ T \rfloor} f(t_k) {
\zeta}_{k, n} \rightarrow\frac{T}{2} \intl^t_0
f(s) \,dB_s,
\]
in $L^p$, as $n $ tends to infinity.
This convergence still holds true when we replace the above $ {\zeta
}_{k, n} $ by
\[
\tilde{\zeta}_{k, n}= \intl_{t_k}^{t_{k+1}}
(B_s-B_{\eta(s)} ) \,d s .
\]
\end{cor}

\begin{pf}As before, to prove the corollary it
suffices to
show that the conditions in Proposition~\ref{prop7.1} are all
satisfied here.
Let us first consider the convergence for $ {\zeta}_{k, n} $.
Set
\[
g _n (t): =n \sum_{k=0}^{ \lfloor{nt}/ T \rfloor}
{\zeta}_{k,
n},
\]
where $ {\zeta}_{k, n} $ is defined in (\ref{7.13e}).
Condition (ii) follows from estimate (\ref{11.a}) in Lemma~\ref{lem11.3}
by taking $F \equiv1$ and $\nu=1 $.
The covariance of the process $g_n $ is given by
\begin{eqnarray*}
\mE\bigl(g _n(t)g_{ n'}(t)\bigr) & = &\alpha_H
n n' \intl^t_0 \intl^t_0
\bigl(u-\eta_n(u)\bigr) \bigl(v- \eta_{ n'} (v)\bigr)
\mu(du\,dv)
\\
&\rightarrow&\frac{T^2}{ 4} \alpha_H \intl^t_0
\intl^t_0 |u-v|^{2H-2} \,du\,dv
\\
&= &\frac{T^{2 }}{4} t^{2H}
\end{eqnarray*}
as $n, n'\rightarrow\infty$,
which implies that $g_n(t)$ is a Cauchy sequence in $L^2$.
Here $\eta_n(t) = \frac{T}n i $ when $ \frac{T}n i \leq t < \frac
{T}n (i+1
)$ and $ \eta_{n'}(t) = \frac{T}{n'} i $ when $ \frac{T}{n'} i \leq
t <
\frac{T}{n'} (i+1)$.
In fact, we can also calculate the kernel of the limit of $z_n (t)$.
Suppose that $\phi_n \in\cH$ satisfies $ g _n (t) = \delta( \phi
_n(t) )$. Then for any $\psi\in\mathcal{H}$,
\begin{eqnarray*}
\langle n \phi_n, \psi\rangle_{\mathcal{H}} &=&n
\alpha_H \intl^T_0 \intl^{\eta(t) }_0
\bigl(u-\eta(u) \bigr) \psi(v)|u-v|^{2H-2}\,du \,dv
\\
& \rightarrow&\frac{T}{2} \langle\psi,\mathbf{1}_{[0, t]}\rangle
_{\cH
},
\end{eqnarray*}
as $n \rightarrow+\infty$. This implies that the kernel of the limit
of $g_n (t)$ is $\frac{T}{2} \mathbf{1}_{[0, t]}$; in other words, the
random variable $g_n(t)$ converges in $L^2$ to $\frac{T}{2} B_t$.

The convergence result for $\tilde{\zeta}_{k, n}$ can be shown by
noticing that
\[
\tilde{\zeta}_{k, n}=\int_{t_k}^{t_{k+1 } }
(B_s - B_{\eta(s) } ) \,ds = \frac{T}{n} (B_{t_{k+1 } }
-B_{t_k}) - \int_{t_k}^{t_{k+1 } } \bigl(s-\eta
(s) \bigr) \,dB_s.
\]
This completes the proof of the corollary.
\end{pf}

\section{Asymptotic error of the modified Euler scheme in case \texorpdfstring{$H \in(\frac{3}{4},1)$}{Hin(3/4,1)}}\label{section thmm2 proof}
The limit theorems for weighted sums proved in the previous section
allow us to derive the $L^p$-limit of
the quantity $n(X_t-X^n_t)$ in the case $H \in(
\frac{3}4, 1)$.

%
\begin{theorem}\label{theorem 2}
Let $H \in( \frac{3}4, 1)$. Suppose that $X$ and $X^n$ are defined by
(\ref{e.1.1}) and (\ref{e.1.6}), respectively. Let $Z^{ij}$, $i, j = 1,
\ldots, m$ be the matrix-valued generalized
Rosenblatt process defined in Section~\ref{section 2.3}. Assume
$\sigma
\in
C^5_b(\mathbb{R}^d; \mathbb{R}^{d \times m})$ and $b \in
C^4_b(\mathbb{R}^d; \mathbb{R}^{d })$. Then
\[
n \bigl(X_t - X^{ n}_t\bigr) \rightarrow
\barU_t
\]
in $L^p (\Omega)$ as $ n $ tends to infinity, where $\{ \bar U_t,
t \in[0, T]\}$ is the solution of the following linear stochastic
differential equation:
%
%
\begin{eqnarray}
\label{eqn 47} %
\barU_t& =& \intl^t_0
\nabla b (X_s) \barU_s \,ds + \sum
_{j=1}^m \intl^t_0 \nabla
\sigma^j(X_s) \barU_s \,dB_s^{j}
\nonumber\\
&&{}+ \sum_{i, j=1}^m \intl^t_0
\bigl( \nabla\sigma^j \sigma^i \bigr) (X_s)
\,dZ^{ij}_s
\nonumber
\\[-8pt]
\\[-8pt]
\nonumber
&&{} + \frac{T}{2} \intl^t_0 ( \nabla b b) (X
_{s}) \,ds + \frac{T}{2} \intl^t_0 (
\nabla b \sigma) (X _{s } ) \,dB _{s}\\
&&{} + \frac{T }{2}
\sum_{j=1}^m \intl^t_0
\bigl( \nabla\sigma^j b \bigr) (X_s)
\,dB^j_s. \nonumber %
\end{eqnarray}
\end{theorem}
\begin{pf}
Recall the decomposition $Y_t = X_t-X^n_t$ given in (\ref{eqn 36}) and
(\ref{eqn20}).
We have shown that $n I_{13} (t)$, $nI_{4,j} (t)$, $nE_{1,j}(t)$ and
$nE_{3,j} (t)$
converge in $L^p$ to zero for each $t \in[0, T]$. It remains to show
the $L^p$ convergence of $ n I_{11} (t) $, $n
I_{12,j}(t)$, $n I_{2,j} (t)$ and $n E_{2,j}(t) $ and identify their limits.

\textit{Step} 1. Recall $\wt{E}_{2,j}(t)$ is defined in
(\ref{e6.5}).
It has been shown in the proof of Theorem~\ref{theorem 1} that
$n(E_{2,j}(t) -\wt{E}_{2,j}(t) )$ converges to zero in $L^p$. On the
other hand,
applying Corollary~\ref{prop 5} to $n\wt{E}_{2,j} (t)$ yields
\[
n \wt{E}_{2,j} (t) \rightarrow\sum_{ i=1}^m
\intl^t_0 \Lambda_t \Gamma_s
\bigl( \nabla\sigma^j \sigma^i\bigr) (X_s)
\,dZ^{ij}_s \qquad\mbox{in } L^p.
\]
Therefore, $nE_{2,j}(t)$ converges in $L^p$, and the limit is the same
as $ n \wt{E}_{2,j} (t)$.

\textit{Step} 2.
Denote
\[
\tilde I_{2,j} (t) = \Lambda_t \sum
_{k=0}^{ { \lfloor{nt}/{T} \rfloor} } \Gamma_{t_k} \bigl(\nabla
\sigma^j b\bigr) (X _{t_k}) \int_{t_k}^{t_{k+1} \wedge t}
\bigl(s- \eta(s)\bigr) \,dB^j_s,
\]
for $t \in[0, T]$ [as before, we define $t_{n+1} = \frac{T}n (n+1)$].
Applying Corollary~\ref{prop 6} to $n \tilde{I}_{2,j}(t)$ yields
\[
n \tilde{I}_{2,j}(t) \rightarrow\frac{T}{2}
\intl^t_0 \Lambda_t \Gamma_s
\bigl( \nabla\sigma^j b\bigr) (X_s) \,dB^j_s
\qquad\mbox{in } L^p.
\]
We want to show that $nI_{2,j}(t)$ and $n \tilde{I}_{2,j}(t)$ have the
same limit in $L^p$.
Write
%
%
\begin{eqnarray}
\label{e 8.2} 
&& n\bigl(I_{2,j}(t) - \tilde{I}_{2,j}
(t)\bigr)\nonumber
\\
&&\qquad= n \intl^t_0 \bigl( \Lambda_t^n
\Gamma_s^n b^j_2 (s) -
\Lambda_t \Gamma_s \tilde{b}^j_2
(s)\bigr) \bigl(s- \eta(s)\bigr) \,dB^j_s
\\
&&\qquad\quad{} + n \intl^t_0 \Lambda_t \bigl(
\Gamma_s \tilde{b}^j_2 (s) -
\Gamma_{\eta(s)} \bigl(\nabla\sigma^j b\bigr) (X
_{\eta(s)}) \bigr) \bigl(s- \eta(s)\bigr) \,dB^j_s,
\nonumber
\end{eqnarray}
where $ \tilde{b}^j_2 (s) = \int_0^1 \nabla\sigma^j (\theta X_s +
(1-\theta) X_{\eta(s ) } )b(X_{\eta(s)}) \,d\theta$.
It suffices to show that the two terms on the right-hand side of (\ref
{e 8.2}) both converge to zero in $L^p$.
The convergence of the second term follows from estimate (\ref{e
11.11}) of Lemma~\ref{lem11.5}.
Lemma~\ref{lem3.5} implies that the $L^p$-norms of $ [ \Lambda_t^n
\Gamma_s^n b^j_2 (s)
- \Lambda_t \Gamma_s \tilde{b}^j_2 (s) ]$ and its Malliavin
derivative converge to
zero as $n \rightarrow\infty$. So applying Lemma~\ref{lem11.3} (\ref{11.a})
with $\nu=1$ and $F_s = \Lambda_t^n \Gamma_s^n b^j_2 (s)
- \Lambda_t \Gamma_s \tilde{b}^j_2 (s)$, we obtain the convergence of
the first term.

\textit{Step} 3. Following the lines in step 2 we can
show that $n I_{12,j} (t)$
converges in $L^p$ to
\[
\frac{T}{2} \intl^t_0 \Lambda_t
\Gamma_{s } \bigl( \nabla b \sigma^j\bigr) (X
_{s } ) \,dB^j_{s}.
\]
Instead of (\ref{11.a}) and { (\ref{e 11.11})} in step 2, we need to use
estimates (\ref{11.b}) and (\ref{e 11.10}) here.

Similarly, it can be shown that $nI_{11}$ converges in $L^p$ to
\[
\frac{T}{2} \intl^t_0 \Lambda_t {
\Gamma_s } ( \nabla b b) (X _{s}) \,ds.
\]

\textit{Step 4.} We have shown that $n(X_t- X^n_t)$
converges in $L^p$ to $\barU_t$, where we define, for each $t \in[0, T]$,
\begin{eqnarray*}
\barU_t &=& \sum_{i, j=1}^m
\intl^t_0 \Lambda_t \Gamma_s
\bigl( \nabla\sigma^j \sigma^i \bigr) (X_s)
\,dZ^{ij}_s + \frac{T}{2} \intl^t_0
\Lambda_t \Gamma_s ( \nabla b b) (X _{s}) \,ds
\\
&&{} + \frac{T}{2} \intl^t_0 \Lambda_t
\Gamma_s ( \nabla b \sigma) (X _{s } ) \,dB _{s}
+ \frac{T }{2} \sum_{j=1}^m
\intl^t_0 \Lambda_t \Gamma_s
\bigl(\nabla\sigma^j b\bigr) (X_s) \,dB^j_s.
\end{eqnarray*}
The theorem follows from the fact that the process $\barU$ satisfies
equation (\ref{eqn 47}).
\end{pf}

\section{Weak approximation of the modified Euler scheme}\label{sec9}
The next result provides the weak rate of convergence for the modified
Euler scheme (\ref{e.1.6}).
%

\begin{theorem}\label{thmm5.1}
Let $X$ and $X^n$ be the solution to equations (\ref{e.1.1}) and (\ref
{e.1.6}),
respectively. Suppose that $b \in C_b^3(\mathbb{R}^d; \mathbb{R}^{d })
$, $\sigma
\in C_b^4(\mathbb{R}^d; \mathbb{R}^{d\times m}) $.
Then for any function $f \in C_b^3( \mathbb{R}^d)$ there exists a
constant $C$ independent of $n$ such that
%
%
\begin{equation}
\label{e.5.3} \sup_{0 \leq t \leq T} \bigl|\mathbb{E} \bigl[f(X_t)
\bigr] - \mathbb{E} \bigl[f\bigl(X^n_t\bigr) \bigr] \bigr|
\leq Cn^{-1}.
\end{equation}
If we further assume that $b \in C^4$, $\sigma\in C^5$ and $f \in
C^4$, then for each $t \in[0, T]$, the sequence
\[
n \bigl\{\mathbb{E} \bigl[f(X_t) \bigr] - \mathbb{E} \bigl[f
\bigl(X^n_t\bigr) \bigr] \bigr\} ,\qquad n\in\NN,
\]
converges as $n$ tends to infinity, and the limit is equal to the sum
of the following two quantities:
%
%
\begin{eqnarray}
\label{e.5.1} %
& &\frac{\alpha_H^2 T}{2} \sum
_{j,i=1}^m \int^t_0
\int^t_0 \int^t_0
\mathbb{E} \bigl\{ D_u^iD_r^j
\bigl[ \nabla f(X_t)\Lambda_t \Gamma_s \bigl(
\nabla\sigma^j \sigma^i \bigr) (X_s) \bigr]
\bigr\}
\nonumber
\\[-8pt]
\\[-8pt]
\nonumber
&&\hspace*{68pt}\qquad{}\times |u-s|^{2H-2}|s-r|^{2H-2} \,du\,ds\,dr %
\end{eqnarray}
and
%
%
\begin{eqnarray}
\label{e.5.1ii} 
&& \frac{T}{2} \mathbb{E} \Biggl\{ \nabla
f(X_t) \Lambda_t \Biggl[ \intl^t_0
\Gamma_s ( \nabla b b) (X _{s}) \,ds +
\intl^t_0 \Gamma_s ( \nabla b \sigma) (X
_{s } ) \,dB _{s}
\nonumber
\\[-8pt]
\\[-8pt]
\nonumber
&&\hspace*{145pt}{} + \sum_{j=1}^m \intl^t_0
\Gamma_s \bigl( \nabla\sigma^j b\bigr)
(X_s) \,dB^j_s \Biggr] \Biggr\}.
\nonumber
\end{eqnarray}
\end{theorem}

\begin{pf}
We use again decompositions (\ref{eqn 36}) and (\ref{eqn20}) of $Y_t =
X_t - X^n_t$, $t \in[0, T]$,
and we continue to use the notation there. Given a function $f \in
C_b^3( \mathbb{R}^d)$, we can write
\[
n \bigl\{ \mathbb{E} \bigl[f(X_t) \bigr] - \mathbb{E} \bigl[f
\bigl(X^n_t\bigr) \bigr] \bigr\}
=
n\int^1_0 \mathbb{E} \bigl[\nabla f\bigl(
Z^{\theta}_t \bigr) Y_t \bigr] \,d \theta,
\]
where we denote $Z^{\theta}_t = \theta X_t + (1 - \theta)X^n_t$, $0
\leq t \leq T $.

\textit{Step} 1. In this step, we show that $\sup_{0 \leq t \leq T}
\vert \mathbb{E} [ \nabla f( Z^{\theta}_t ) Y_t
] \vert
\leq Cn^{-1}$, which implies (\ref{e.5.3}). From estimates (\ref{eqn19})
and (\ref{eq19a})
it follows that this inequality is true when $Y $ is replaced by
$ I_{11} $, $ I_{13} $, $
I_{12, j } $, $ I_{ 2, j } $ or $ I_{4, j } $.
Therefore, it suffices to show that $ \vert \mathbb{E} [
\nabla
f( Z^{\theta}_t ) E_{i,j}(t)
] \vert
\leq Cn^{-1}$ for $i=1,2,3$ and $j=1,\ldots, m$, where $E_{ij}(t)$ are
defined in Theorem~\ref{thmm4.1} step 3.
Consider first the term $i=2$.
The use of expression
(\ref{4.12e}) and an application of the integration by parts formula yield
%
%
\begin{eqnarray}
\label{eqn23} %
& &\mE\bigl[\nabla f\bigl( Z^{\theta}_t
\bigr) E_{2,j}(t) \bigr]\nonumber
\\
&&\qquad = \mE\Biggl[ \nabla f\bigl( Z^{\theta}_t \bigr)
\Lambda^n_t \sum_{i=1}^m
\sum_{k=0}^{ \lfloor{nt}/ T \rfloor} F^{n, i, j}_{t_k}
\int_{t_k}^{t_{k+1}\wedge t} \int_{t_k}^s
\delta B^i_u \delta B^j_s
\Biggr]
\nonumber
\\[-8pt]
\\[-8pt]
\nonumber
&&\qquad= \alpha_H^2 \sum_{i=1}^m
\sum_{k=0}^{ { \lfloor{nt}/{T} \rfloor
} } \mE\biggl[ \int
_0^t \int_{t_k}^{t_{k+1}\wedge t}
\int_0^t \int_{t_k}^s
D^i_v D^j_r \bigl[ \nabla f
\bigl( Z^{\theta}_t \bigr) \Lambda^n_t
F^{n, i,
j}_{t_k} \bigr] \\
&&\hspace*{203pt}{}\times \mu( \,du \,dv)\mu( ds \,dr ) \biggr],\nonumber
\end{eqnarray}
where we recall that
$ F^{n, i, j}_{t} = \Gamma^n_t ( \nabla\sigma^j \sigma^i ) (X^n_t)$.
[As before, in the above equation we set $t_{n+1} = \frac{T}n (n+1)$.]
Therefore,
\begin{eqnarray*}
\bigl\vert\mE\bigl[\nabla f\bigl( Z^{\theta}_t \bigr)
E_{2,j}(t) \bigr] \bigr\vert& \leq &C \sum
_{k=0}^{n-1 } \int_{t_k}^{t_{k+1} }
\int_0^t \int_{t_k}^{t_{k+1}}
\int_0^t \mu(du \,dv)\mu( dr \,ds )
\\
& \le& Cn^{-1}.
\end{eqnarray*}
%
For the term containing $E_{1,j}$ we can write
\[
\mE\bigl[\nabla f\bigl( Z^{\theta}_t \bigr)
E_{1,j}(t) \bigr] = \sum_{i=1}^{m}
\mE\biggl[ \int_0^t H^{n,i,j}_{s}
\bigl(B^i_s- B^i_{\eta
(s)}\bigr)
\,dB^j_s \biggr],
\]
where $H^{n,i,j}_s= \nabla f( Z^{\theta}_t ) \Lambda^n_t [
\Gamma
^n_s \sigma^{j,i}_2(s)- \Gamma^n_{\eta(s)}
(\nabla\sigma^j\sigma^i)(X^n_{\eta(s)} ) ]$. An application of
the relation between the Skorohod and
path-wise integrals (\ref{e.2.3}) yields
\begin{eqnarray*}
&& \mE\biggl[ \int_0^t H^{n,i,j}_{s}
\bigl(B^i_s- B^i_{\eta(s)}\bigr)
\,dB^j_s \biggr]
\\
&&\qquad= \alpha_H \int_0^T \int
_0^t \mE\bigl[ D^j_u
\bigl(H^{n,i,j}_{s} \bigl(B^i_s-
B^i_{\eta(s)}\bigr) \bigr) \bigr] |s-u|^{2H-2} \,ds\,du
\\
&&\qquad= \alpha_H \int_0^T \int
_0^t \mE\bigl[ D^j_u
H^{n,i,j}_{s} \bigl(B^i_s-
B^i_{\eta(s)}\bigr) \bigr] |s-u|^{2H-2} \,ds\,du
\\
&&\qquad\quad{} + \alpha_H \int_0^T \int
_0^t \mE\bigl[ H^{n,i,j}_{s}
\bigr] \mathbf{1}_{[\eta(s),s]}(u) \delta_{ij} |s-u|^{2H-2}
\,ds\,du
\\
&&\qquad:= A_1+A_2.
\end{eqnarray*}
By the integration by parts we see that $A_1$ is equal to
\[
\alpha_H ^2 \int_0^T
\int_0^t \int_0^T
\int_0^s\mE\bigl[ D^i_rD^j_u
H^{n,i,j}_{s} \bigr] \mathbf{1}_{[\eta(s),s]}(v)
|v-r|^{2H-2} |s-u|^{2H-2} \,dv \,dr \,ds\,du.
\]
Using $\sup_{r,u,s} \vert \mE[ D^i_rD^j_u H^{n,i,j}_{s}
] \vert \le C n^{-\beta}$ for any $\frac{1}2<\beta<H$ we obtain
%
%
\begin{equation}
\label{A1} |A_1| \le Cn^{-1-\beta}.
\end{equation}
On the other hand, it is easy to show by the definitions of $\Gamma^n$,
$X^n$ and $X$ that the quantity $ [ \Gamma^n_s \sigma^{j,i}_2(s)-
\Gamma^n_{\eta(s)}
(\nabla\sigma^j\sigma^i)(X^n_{\eta(s)} ) ]$ can be expressed as
the sum of integrals over the interval $[\eta(s), s]$. So by applying
(\ref{e.2.3}) and integration by parts we can show that
$\vert\mE[ H^{n,i,j}_{s} ] \vert \le Cn^{-1}$, which implies
%
%
\begin{equation}
\label{A2} |A_2| \le Cn^{-2H}.
\end{equation}
From (\ref{A1}) and (\ref{A2}) we conclude that $\vert \mE
[\nabla
f( Z^{\theta}_t ) E_{1,j}(t) ] \vert \le Cn^{-1} $.
Finally, for the term containing $E_{3,j}$ we have
\[
\mE\bigl[\nabla f\bigl( Z^{\theta}_t \bigr)
E_{3,j}(t) \bigr] = \int_0^t \mE
\bigl[ J^{n,i,j}_s \bigr] (s-\eta(s) ^{2H-1} \,ds,
\]
where $J^{n,i,j}_s= H\nabla f( Z^{\theta}_t ) \Lambda^n_t (
\Gamma
^n_{\eta(s)} -\Gamma^n_s ) \sigma^j_0(s) $. By expressing the
term $ (\Gamma^n_{\eta(s)} -\Gamma^n_s ) $ as the sum of integrals over
the interval $[\eta(s), s]$ and then applying (\ref{e.2.3}) and
integration by parts, we can show that $ \sup_{s \in[0,T]} \mE
[J^{n,i,j}_s] \leq C n^{-1} $. This implies
%
%
\begin{equation}
\label{3,j} \bigl\vert\mE\bigl[\nabla f\bigl( Z^{\theta}_t
\bigr) E_{3,j}(t) \bigr] \bigr\vert\le C n^{-2H},
\end{equation}
which completes the proof of (\ref{e.5.3}).

\textit{Step} 2. Now we show the second part of the theorem.
From estimates (\ref{eqn19}), (\ref{eq19a}), (\ref{A1}), (\ref{A2}) and
(\ref{3,j}) we see that
the expression $n\int^1_0 \mathbb{E} [\nabla f( Z^{\theta}_t )
Y_t ] \,d \theta$ converges
to zero as $n$ tends to infinity when $Y_t$ is replaced by
$I_{13}(t)$, $I_{4,j}(t)$,
$E_{1,j}(t)$ or $E_{3,j}(t)$. Therefore, it suffices to consider
$n\int^1_0 \mathbb{E} [\nabla f( Z^{\theta}_t ) Y_t ] \,d
\theta
$ when
$Y_t$ is replaced by the remaining terms in the decomposition of $Y_t$.

Consider first the term $E_{2,j}(t)$, and denote
\[
G_{s,r,v}^{i,j} = D^i_v
D^j_r \bigl[ \nabla f( X_t )
\Lambda_t \Gamma_{s} \bigl( \nabla\sigma^j
\sigma^i\bigr) (X_{s} ) \bigr] .
\]
It is clear that
%
%
\begin{eqnarray}\label{9.4e}
& & n \alpha_H^2 \sum
_{i,j=1}^m \sum_{k=0}^{ \lfloor{nt}/ T
\rfloor}
\int_0^t \int_{t_k}^{t_{k+1}\wedge t}
\int_0^t \int_{t_k}^s
G_{t_k,r,v}^{i,j} \mu( du \,dv )\mu( ds \,dr )
\nonumber
\\[-8pt]
\\[-8pt]
\nonumber
 &&\qquad \rightarrow\frac{\alpha_H^2T}{2} \sum
_{i,j=1}^m \int_0^t
\int_0^t \int_0^t
G_{s,r,v}^{i,j} |s-v|^{2H-2} |r-s|^{2H-2} \,ds \,dv
\,dr,
\end{eqnarray}
almost surely. Therefore, by the dominated convergence theorem,
the expectation of the left-hand side of the above expression
converges to the expectation of the right-hand side, which is term
(\ref{e.5.1}). From Lemma~\ref{lem3.5}, we have
%
%
\begin{eqnarray}
\label{eqn25} &&\bigl\llVert D^i_v D^j_r
\bigl[ \nabla f\bigl( Z^{\theta}_t \bigr)
\Lambda^n_t F^{n, i,
j}_{t_k} \bigr] -
D^i_v D^j_r \bigl[ \nabla f(
X_t ) \Lambda_t \Gamma_{t_k} \bigl( \nabla
\sigma^j \sigma^i \bigr) (X_{t_k}) \bigr] \bigr
\rrVert_{p}
\nonumber
\\[-8pt]
\\[-8pt]
\nonumber
&&\qquad \leq C n^{1-2 \beta},
\end{eqnarray}
which, together with equation (\ref{eqn23}), implies that $
n\sum_{j=1}^m \mE[\nabla f( Z^{\theta}_t ) E_{2,j} (t) ]
$ converges to the same limit as the expectation of the left-hand side
of (\ref{9.4e}).

The results in steps 2 and 3 of the proof of Theorem~\ref{theorem 2}
imply that the
terms $ n \mathbb{E} [\nabla f( Z^{\theta}_t ) I_{11} (t) ] $,
$n \mathbb{E} [\nabla f( Z^{\theta}_t ) I_{12,j} (t) ] $ and
$ n \mathbb{E} [\nabla f( Z^{\theta}_t ) \sum_{j=1}^m I_{2,j} (t) ] $
converge to the second, third and fourth term in (\ref{e.5.1ii}),
respectively. For example,
let us consider $ n \mathbb{E} [\nabla f( Z^{\theta}_t ) \sum_{j=1}^m
I_{2,j} (t) ] $.
We have shown in Theorem~\ref{theorem 2} that
\[
n {I}_{2,j}(t) \rightarrow\frac{T}{2} \Lambda_t
\intl^t_0 \Gamma_s \bigl( \nabla
\sigma^j b\bigr) (X_s) \,dB^j_s
\]
in $L^p$ for any $p\ge1$.
So it follows from the H\"older inequality that
\[
\biggl\vert\mE\biggl[ n \nabla f\bigl( Z^{\theta}_t
\bigr) I_{2,j} (t) - \nabla f( X_t ) \frac{T}{2}
\Lambda_t \intl^t_0 \Gamma_s
\bigl( \nabla\sigma^j b\bigr) (X_s) \,dB^j_s
\biggr] \biggr\vert\rightarrow0
\]
as $n \rightarrow\infty$. The other two terms can be studied in
similar way.
This completes the proof of the theorem.
\end{pf}

%
\begin{remark}
Theorem~\ref{thmm5.1} may be used to construct a Richard extrapolation
scheme with error bound $o(n^{-1})$.
\end{remark}

\section{Rate of convergence for the Euler scheme}\label{sec10}
In this section, we apply our approach based on Malliavin calculus
developed in Section~\ref{section strong conv} to
study the rate of convergence of the naive Euler scheme
defined in (\ref{e.1.2}).
Our first result is the rate of the strong convergence of the naive
Euler scheme.
As we will see, the weak rate of convergence and the rate of strong
convergence are the same for the naive Euler scheme.
We still use $X^n$ to represent the naive Euler scheme (\ref{e.1.2}).
This will not cause confusion since we will only deal with this scheme
in this section.
%

\begin{theorem}\label{t.6.1}
Let $X$ and $X^n$ be the processes defined in (\ref{e.1.1}) and (\ref
{e.1.2}), respectively. Suppose that $b \in C_b^1(\mathbb{R}^d;
\mathbb
{R}^{d })$ and $\sigma\in C_b^2(\mathbb{R}^d; \mathbb{R}^{d\times
m})$. Then for each $p \geq1$, we have
\[
n^{2H-1} \sup_{t \in[0, T]} \mE\bigl(\bigl|X_t -
X^n_t\bigr|^p\bigr)^{{1}/p} \leq C.
\]
If we assume $b \in C_b^3(\mathbb{R}^d; \mathbb{R}^{d })$ and $\sigma
\in C_b^4(\mathbb{R}^d; \mathbb{R}^{d\times m})$, then as $n $ tends to
infinity,
\[
n^{2H-1}\bigl(X_t - X^n_t\bigr)
\rightarrow\frac{T^{2H-1}}{2} \sum_{j=1}^m
\int^t_0 \Lambda_t
\Gamma_s \bigl(\nabla\sigma^j\sigma^j\bigr)
(X_s) \,ds,
\]
where $\Lambda$ is the solution to linear equation (\ref{eqn 46}) and $
\Gamma_t = \Lambda^{-1}_t $, and the convergence holds in $L^p$ for all
$p\ge1$.
\end{theorem}

\begin{pf}We let $Y_t = X_t - X^n_t$, $t \in[0, T]$. Then
as in the proof of Theorem~\ref{thmm4.1}, we can derive the
decomposition of $Y_t$
\begin{eqnarray*}
Y_t &=& \Lambda^n_t \intl^t_0
\Gamma^n_s b_3 (s) \Biggl[ b
\bigl(X^n_{\eta(s)}\bigr) \bigl(s- \eta(s)\bigr) + \sum
_{l=1} ^m\sigma^l
\bigl(X^n_{\eta(s) } \bigr) \bigl( B^l_{s}
- B^l_{\eta(s )
} \bigr) \Biggr] \,ds
\\
&&{}+ \sum_{j=1}^m \intl^t_0
\Lambda_t^n \Gamma_s^n
b^j_2 (s) \bigl(s- \eta(s)\bigr) \,dB^j_s
\\
&&{}+ \sum_{i, j =1}^m \intl^t_0
\Lambda_t^n \Gamma_s^n \sigma
_2^{j,i}(s) \bigl( B^i_{s} -
B^i_{\eta(s ) } \bigr) \,dB^{j}_s
\\
&=:& I_1 (t) +I_2(t) +I_3(t)
+I_4(t),
\end{eqnarray*}
where $ \Lambda^n$, $ \Gamma^n$, $b_2^j(s)$, $\sigma_2^{j, i}(s)$ and
$b_3(s)$ are the same terms as those defined in the proof of
Theorem~\ref{thmm4.1} with the scheme $X^n$ replaced by the classical Euler
scheme~(\ref{e.1.2}).

It is clear that $\|I_1 (t) \|_p \leq Cn^{-1}$. On the other hand,
estimates (\ref{11.a}) and (\ref{11.b}) of Lemma~\ref{lem11.3} imply that
$\|I_2 (t) \|_p \leq Cn^{-1}$ and $\| I_3 (t) \|_p \leq Cn^{-1}$.
Finally, as in the proof of (\ref{e 11.10}) in Lemma~\ref{lem11.5} we
obtain $\| I_4 (t) \|_p \leq C n^{1-2H}$.
This completes the proof of the first part of the theorem.

Applying the integration by parts to $I_4(t)$ yields
\begin{eqnarray*}
&& \int^t_0 \Lambda^n_t
\Gamma^n_s \sigma_2^{j, i} (s)
\bigl(B^i_s - B^i_{\eta(s)}\bigr)
\,dB^j_s
\\
&&\qquad= \int_0^t \Lambda^n_t
\Gamma^n_s \sigma_2^{j, i} (s)
\bigl(B^i_s - B^i_{\eta(s)}\bigr)
\delta B^j_s
\\
&&\qquad\quad{} + \alpha_H \int^t_0 \int
^t_0 D_r^j \bigl[
\Lambda^n_t \Gamma^n_s
\sigma_2^{j, i} (s) \bigr] \bigl(B^i_s
- B^i_{\eta(s)}\bigr) \mu(ds \,dr )
\\
&&\qquad\quad{} +\delta_{ij} \alpha_H \int^t_0
\int^t_0 \Lambda^n_t
\Gamma^n_s \sigma_2^{j, i} (s)
\mathbf{1}_{[\eta(s), s]}(r) \mu(ds \,dr )
\\
&&\qquad=: A^1_n (t) + A^2_n (t) +
A^3_n (t).
\end{eqnarray*}
From (\ref{11.c}) we have $\|A^1_n (t) \|_p \leq C\gamma_n^{-1}$.
Applying (\ref{11.b}) with $F_u $ replaced by
\[
\int_0^t D_r^j
\bigl[ \Lambda^n_t \Gamma^n_s
\sigma_2^{j, i} (s) \bigr] |r-u|^{2H-2} \,dr
\]
we obtain $\|A^2_n (t) \|_p \leq Cn^{-1}$. So it suffices to identify
the limit of $n^{2H-1}A^3_n(t)$ in $L^p$.
It follows from Lemma~\ref{lem3.5} and Remark~\ref{remark4.3} that
\[
\bigl\llVert\Lambda^n_t \Gamma^n_s
\sigma_2^{j,j} (s) - \Lambda_t
\Gamma_s \bigl(\nabla\sigma^j\sigma^j\bigr)
(X_s) \bigr\rrVert_p \leq C n^{1-2\beta}.
\]
Therefore, $n^{2H-1}A^3_n(t)$, and the quantity
\[
n^{2H-1} \int^t_0 \int
^t_0 \Lambda_t
\Gamma_s \bigl(\nabla\sigma^j\sigma^j\bigr)
(X_s) \mathbf{1}_{[\eta(s), s]}(r) |r-s|^{2H-2} \,ds \,dr
\]
converges to the same value in $L^p$.
The theorem now follows by noticing that
\begin{eqnarray*}
&&n^{2H-1} \int^t_0 \int
^t_0 \Lambda_t
\Gamma_s \bigl(\nabla\sigma^j\sigma^j\bigr)
(X_s) \mathbf{1}_{[\eta(s), s]}(r) |r-s|^{2H-2} \,ds \,dr
\\
&&\qquad= n^{2H-1} \int^t_0
\Lambda_t \Gamma_s \bigl( \nabla\sigma^j
\sigma^j\bigr) (X_s) \frac{(s - \eta(s ) )^{2H-1}}{ 2H-1 } \,ds
\\
&&\qquad\rightarrow\frac{T^{2H-1}}{2 \alpha_H} \int^t_0
\Lambda_t \Gamma_s \bigl(\nabla\sigma^j
\sigma^j\bigr) (X_s) \,ds,
\end{eqnarray*}
in $L^p$ for all $p\ge1$.
\end{pf}

As a consequence of the above theorem, we can deduce the following result.
%

\begin{cor}
Let $X$ and $X^n$ be the processes defined in (\ref{e.1.1}) and (\ref
{e.1.2}), respectively. Suppose that $b \in C_b^3 (\mathbb{R}^d;
\mathbb
{R}^{d })$, $\sigma\in C_b^4(\mathbb{R}^d; \mathbb{R}^{d\times m})$
and $f \in C_b^2(\mathbb{R}^d) $. Let $\Lambda$ be defined in (\ref{eqn
46}). Then we have the following $L^p$-convergence as $n
\rightarrow\infty$ for all $p \geq1$:
\[
{ n^{2H-1} } \bigl[ f\bigl(X_t^n \bigr) -
f(X_t ) \bigr] \rightarrow\frac{T^{2H-1}}{2} \sum
_{j=1}^m \int^t_0
\nabla f(X_t) \Lambda_t \Gamma_s \bigl(
\nabla\sigma^j\sigma^j\bigr) (X_s) \,ds.
\]
\end{cor}

\begin{pf}
We can write
\[
{ n^{2H-1} } \bigl[ f\bigl(X^n_t
\bigr) - f(X_t) \bigr] = { n^{2H-1} } \biggl( \int
^1_0 \nabla f\bigl( Z^{\theta}_t
\bigr) \,d \theta\biggr) \bigl(X_t^n-X_t
\bigr), %
\]
where we denote $Z^{\theta}_t = \theta X_t + (1 - \theta)X^n_t$, $ t
\in[0, T] $.
Then the result follows from Theorem~\ref{t.6.1}, the convergence of
$X_t^n$ to $X_t$ and the assumption on $f$.
\end{pf}

The above corollary implies the following weak approximation result:
\begin{eqnarray*}
&& \lim_{n\rightarrow\infty} n^{2H-1} \bigl\{ \mathbb{E}
\bigl[f(X_t)\bigr] - \mathbb{E}\bigl[f\bigl(X^n_t
\bigr) \bigr] \bigr\}
\\
&&\qquad= \frac{ T^{2H-1} }{ 2 } \sum_{j=1}^m
\int^t_0 \mathbb{E}\bigl[ \nabla
f(X_t) \Lambda_t \Gamma_s \bigl(\nabla
\sigma^j\sigma^j\bigr) (X_s) \bigr] \,ds.
\end{eqnarray*}

\begin{appendix}\label{app}
\section*{Appendix}
\subsection{Estimates of a Young integral}\label{section estimate
stochastic integral}
In this section, we give an estimate on the pathwise integral using
fractional calculus.
%

\begin{lemma}\label{lem7.1}
Let $z=\{z_t, t \in[0, T]\}$ be a H\"older continuous function with
index $\beta\in(0, 1)$. Suppose that $f \dvtx\mathbb{R} ^{l+m}
\rightarrow\mathbb{R}$ is continuously differentiable. We denote by
$\nabla_xf$ the $l$-dimensional vector with coordinates $ \frac
{\partial f} {\partial x_i}$, $i=1,\ldots, l$, and by $\nabla_yf$ the
$m$-dimensional vector with coordinates $ \frac{\partial f} {\partial
x_{l+i}}$, $i=1,\ldots, m$.
Consider processes $x=\{x_t, t\in[0,T]\}$ and $y=\{y_t, t\in[0,T]\}$
with dimensions $l$ and $m$, respectively, such that $ \| x\|_{0,T,
\beta'}$ and $ \| y\|_{0,T,\beta',n} $ are finite for each $n\ge1$,
where $\beta' \in(0, 1)$ is such that $\beta'+\beta>1$. Then we have
the following estimates:
\begin{longlist}[(ii)]
\item[(i)] for any $s,t \in[0, T]$ such that $s\le t$ and $s = \eta
(s)$, we have
\begin{eqnarray*}
\label{equation estimate frac int in lemma} 
\biggl\vert\int^t_s
f(x_r, y_{\eta(r)})\,dz_r \biggr\vert& \leq&
K_1 \sup_{r \in[s, t ]} \bigl|f( x_{r },
y_{\eta(r )} )\bigr| \|z \| _{\beta} (t-s)^{\beta}
\\
& &{}+ K _2 \sup_{r_1, r_2 \in[s, t ]} \bigl|\nabla_x f(
x_{r_1 }, y_{\eta
(r_2) } )\bigr| \| x\|_{s, t, \beta'} \|z
\|_{\beta}(t-s)^{ \beta+\beta'
}
\\
& &{}+ K_3 \sup_{r_1, r_2 \in[s, t]} \bigl|\nabla_y
f(x_{r_1 }, y_{r_2} )\bigr| \| y\|_{s, t, \beta', n} \|z
\|_{\beta} (t-s)^{ \beta+\beta' }, %
\end{eqnarray*}
where the $K_i$, $i=1,2,3$, are constants depending on $\beta$ and
$\beta'$;

\item[(ii)]if the function $f$ only depends on the first $l$
variables, then the above estimate holds for all $0\le s\le t\le T$.
\end{longlist}
\end{lemma}

\begin{pf}
Take $\alpha$ such that $ \beta' >
\alpha>1-\beta$.
Let $s, t \in[0, T]$ be such that $s = \eta(s)$ and $s\le t$.
Applying the fractional integration by parts formula in
Proposition~\ref{prop2.1}, we obtain
%
%
\begin{equation}
\label{equation integral by part fractional} %
\biggl\vert\int^t_s
f(x_r, y_{\eta(r)})\,dz_r \biggr\vert\leq\int
^t_s \bigl\vert D^{\alpha}_{s+}
f(x_r, y_{\eta(r)}) \bigr\vert\bigl\vert
D^{1-\alpha}_{t-} (z _r - z _t) \bigr
\vert\,dr. %
\end{equation}
By the definition of fractional differentiation in (\ref{e.2.2}) and
taking into account that $\alpha+\beta-1>0$, we can show that
%
%
\begin{equation}
\label{equation estimate fractional derivative of fBM} %
\bigl| D^{1-\alpha}_{t-} (z
_r - z _t) \bigr| \leq K_ 0 \|z
\|_{\beta} (t -r)^{\alpha+\beta-1},\qquad s \leq r \leq t, %
\end{equation}
where $ K_ 0= \frac{\beta}{( \beta+\alpha-1)\Gamma(\alpha) } $.
On the other hand, using (\ref{e.2.1}) we obtain
%
%
\begin{eqnarray}
\label{equation estimate fractional derivative of
sigma X^n} 
&& \bigl| D^{\alpha}_{s+}f(x_r,
y_{\eta(r)}) \bigr|\nonumber
\\
&&\qquad \leq\frac{1}{ \Gamma(1-\alpha)} \biggl[ \frac{| f(x_r, y_{\eta
(r)}) |}{ ( r -s )^{ \alpha}} + \alpha\int
^r_s \frac{| f(x_r, y_{\eta(r)}) - f(x_u, y_{\eta(u)})
|}{( r-u)^{ \alpha+1 }} \,du \biggr]
\nonumber
\\
&&\qquad \leq\frac{1}{ \Gamma(1-\alpha)}
\nonumber
\\[-8pt]
\\[-8pt]
\nonumber
&&\qquad\quad{}\times \biggl[ \sup_{r \in[s, t ]} \bigl| f(
x_r, y_{\eta(r)} )\bigr| ( r - s )^{ - \alpha}
\\
&&\hspace*{15pt}\qquad\quad{} + \alpha\sup_{r_1, r_2 \in[s, t ]} \bigl|\nabla_x f(
x_{r_1}, y_{\eta(r_2)} )\bigr| \| x\|_{s, t, \beta'} \int
^r_s ( r - u)^{\beta'- \alpha-1 } \,du
\nonumber
\\
&&\hspace*{15pt}\qquad\quad{} + \alpha\sup_{r_1, r_2 \in[s, t ]}\bigl |\nabla_y f(
x_{r_1}, y_{r_2} )\bigr| \| y\|_{s, t, \beta',n } \int
^r_s \frac{| \eta(r)-\eta
(u)|^{\beta'}}{( r - u)^{ \alpha+1 }} \,du \biggr].
\nonumber
\end{eqnarray}
Inequalities (\ref{equation integral by part fractional}), (\ref
{equation estimate fractional derivative of sigma X^n}) and (\ref
{equation estimate fractional derivative of fBM}) together imply
\begin{eqnarray*}
&& \biggl\vert\int^t_s f(x_r,
y_{\eta(r)})\,dz_r \biggr\vert
\\
&&\qquad \leq\frac{1}{ \Gamma(1-\alpha)} \int^t_s \biggl[
\sup_{r \in[s, t
]} \bigl| f( x_r, y_{\eta(r)} )\bigr| ( r - s
)^{ - \alpha}
\\
&&\hspace*{61pt}\qquad\quad{} + \alpha\sup_{r_1, r_2 \in[s, t ]} \bigl|\nabla_x f(
x_{r_1 }, y_{\eta
(r_2 )} )\bigr| \| x\|_{s, t, \beta'} \\
&&\hspace*{61pt}\qquad\quad{}\times\int
^r_s ( r - u)^{\beta'- \alpha
-1 } \,du
\\
&&\hspace*{61pt}\qquad\quad{} + \alpha\sup_{r_1, r_2 \in[s, t ]} \bigl|\nabla_y f(
x_{r_1}, y_{r_2} )\bigr| \| y\|_{s, t, \beta',n }\\
&&\hspace*{134pt}\qquad{}\times \int
^r_s \frac{| \eta(r)-\eta(u)|^{\beta'}}{(
r-u)^{ \alpha+1 }} \,du \biggr]
\\
&&\qquad\hspace*{57pt}\quad{}\times K_0 \|z \|_{\beta} (t -r)^{\alpha+\beta-1} \,dr
\\
&&\qquad \leq K_1 \sup_{r \in[s, t ]} \bigl| f( x_{r},
y_{\eta(r)} )\bigr| \|z \|_{\beta} (t-s)^{\beta}
\\
&&\qquad\quad{} + K _2 \sup_{r_1, r_2 \in[s, t ]} \bigl|\nabla_x f(
x_{r_1}, y_{\eta
(r_2)} )\bigr| \| x\|_{s, t, \beta' } \|z
\|_{\beta}(t-s)^{ \beta+\beta'
}
\\
&&\qquad\quad{} + K_3 \sup_{r_1, r_2 \in[s, t ]} \bigl|\nabla_y f(
x_{r_1}, y_{r_2 } )\bigr| \| y\|_{s, t, \beta', n} \|z
\|_{\beta} (t-s)^{ \beta+\beta' },
\end{eqnarray*}
where $ {K_1= K_0 \frac{\Gamma(\alpha+ \beta)} {\Gamma
(\beta+1)}}$, $ {K_2 =K_0 \frac{\alpha\Gamma(\alpha+\beta)
\Gamma(\beta'-\alpha+1)}{\Gamma(1-\alpha)\Gamma( \beta+ \beta
'+1)(\beta'-\alpha)}}$,
$ {K_3= K_0K_4 \frac{\alpha}{\Gamma(1-\alpha)}}
$ and $K_4 $ is the constant in Lemma~\ref{Lemma estimate discrete
integral}. This completes the proof. \hfill$\Box$

%
\begin{lemma}\label{Lemma estimate discrete integral}
Let $\beta$, $\beta'$ and $\alpha$ be such that $\beta'>\alpha
>1-\beta
$. Then for any $s, t \in[0, T]$ such that $s<t$, $ s = \eta(s)$,
there exists a constant $K_4$ depending on $\alpha$, $\beta$ and $T$,
such that
\[
\int^t_s (t - r)^{ \alpha+\beta-1 }
\int^r_s \frac{| \eta(r)-\eta
(u)|^{\beta' }}{( r-u)^{ \alpha+1 }} \,du \,dr \leq
K_4 (t-s)^{ \beta+
\beta'}. %
\]
\end{lemma}

\begin{pf} Without loss of generality, we let $T=1$. Note
that when $\eta(s)=s < t \leq\eta(s)+ \frac{1}n$, the double integral
equals zero. In the following we will assume $ t > \eta(s) + \frac{1}n$.

We first write
\begin{eqnarray*}
&&\int^t_s (t - r)^{ \alpha+\beta-1 } \int
^r_s \frac{| \eta(r ) -
\eta
(u)|^{\beta'}}{ ( r-u)^{ \alpha+1} }\,du\,dr
\\
&&\qquad= \int^t_{ \eta(s) +{1}/{n} } (t - r)^{ \alpha+\beta-1 } \int
^{\eta(r)}_{\eta(s)} \frac{| \eta(r ) - \eta(u)|^{\beta'}}{ ( r-u)^{
\alpha+1} }\,du\,dr
\\
&&\qquad= \int^t_{ \eta(s) +{1}/{n} } (t - r)^{ \alpha+\beta-1 } \biggl(
\int^{\eta(r)}_{\eta(r) -{1}/n} + \int^{\eta(r) -{1}/n}_{\eta
(s)}
\biggr) \frac{| \eta(r ) - \eta(u)|^{\beta'}}{ ( r-u)^{ \alpha+1}
}\,du\,dr
\\
&&\qquad:= J_1 +J_2.
\end{eqnarray*}

On one hand, notice that in the term $J_{ 2}$ we always have $r - u
>\frac{1}n $, and thus $ \eta(r ) - \eta(u) \leq r-u +\frac{1}n
\leq
2(r-u)$. Therefore,
\begin{eqnarray*}
J_{ 2} &\leq&\int^t_{ \eta(s) +{1}/{n} } (t -
r)^{ \alpha+\beta-1 } \int^{\eta(r) - {1}/n}_{\eta(s)}
\frac{2^{\beta'}(r-u)^{\beta'}}{ (
r-u)^{ \alpha+1} }\,du \,dr
\\
& \leq& K(t-s)^{ \beta+ \beta'}.
\end{eqnarray*}

On the other hand,
\begin{eqnarray*}
J_1& =& \int^t_{ \eta(s) +{1}/{n} } (t -
r)^{ \alpha+\beta-1 } \int^{\eta(r)}_{\eta(r)- {1}/n}
\frac{| \eta(r ) - \eta(u)|^{\beta
'}}{ (
r-u)^{ \alpha+1} }\,du \,dr
\\
&\leq& K n^{-\beta'} (t - s)^{ \alpha+\beta-1 } \int^t_{ \eta(s)
+{1}/{n} }
\biggl[\frac{1}{(r - \eta(r))^{\alpha}} - \frac{1}{(r -
\eta
(r) + {1}/{n})^{\alpha}}\biggr] \,dr
\\
&\leq& K n^{-\beta'} (t - s)^{ \alpha+\beta-1 } \int^t_{ \eta(s)
+{1}/{n} }
\frac{1}{(r - \eta(r))^{\alpha}} \,dr
\\
&\leq& K n^{-\beta'} (t - s)^{ \alpha+\beta-1 } \frac{(\eta(t)
+
{1}/{n})- ( \eta(s) + {1}/{n} ) }{1/n}
n^{\alpha-1}
\\
&\leq& K(t-s)^{ \beta+ \beta' }.
\end{eqnarray*}
The lemma is now proved.
\end{pf}

\subsection{Estimates for some special Young and Skorohod integrals}
In this section we derive estimates for some specific Young and
Skorohod integrals.
We fix $n \in\mathbb{N}$ and consider the uniform partition on $[0, T]$.

%
\begin{lemma}\label{lem11.3}
Let $B=\{B_t, t\in[0,T]\}$ be a one-dimensional fBm with Hurst
parameter $H>\frac{1}2$. Fix $\nu\ge0$ and $p \ge\frac{1}H$.
Let $F=\{F_t, t \in[0, T]\}$ be a stochastic process whose
trajectories are H\"older continuous of order $\gamma>1-H$ and such
that $F_t \in\mathbb{D}^{1,q}$, $t\in[0,T]$, for some $q>p$.
For any $\rho>1$ we set
\[
F_{1,\rho} = \sup_{s,t \in[0,T]} \bigl( \|F_t
\|_\rho\vee\|D_sF_t \| _\rho\bigr).
\]
Then there exists a constant $C$ (independent of $F$) such that the
following inequalities hold for all $0\le s<t \le T$:
%
%
\begin{eqnarray}
\biggl\llVert\int_s^t F_u
\bigl( u - {\eta(u) } \bigr)^{\nu} \,dB_u \biggr\rrVert
_{p} &\leq& C n^{-\nu} (t-s)^H
F_{1,p}, \label{11.a}
\\
\biggl\llVert\int_s^t F_u
(B_u - B_{\eta(u) } ) \,du \biggr\rrVert_{p} & \leq&
C n^{-1} (t-s)^H F_{1,q}. \label{11.b}
\end{eqnarray}
\end{lemma}

\begin{pf*}{Proof of (\ref{11.a})}
Applying (\ref{e.2.3}) we can decompose the Young integral as the sum of
a Skorohod integral plus a complementary term,
%
%
\begin{eqnarray}
\label{11.10e} && \int_s^t F_u
\bigl( u - {\eta(u) } \bigr)^{\nu} \,dB_u \nonumber\\
&&\qquad= \int
_s^t F_u \bigl( u - {\eta(u) }
\bigr)^{\nu} \delta B_u 
\\
&&\qquad\quad{} + \al_H \int_s^t \int
_0^T \bigl(u- \eta(u)\bigr)^{\nu}
D_r F_u |r-u|^{2H-2} \,dr \,du.\nonumber
\end{eqnarray}
It follows from (\ref{eqn 2.10}) that the $L^p$-norm of the first
integral of the right-hand side of (\ref{11.10e}) is bounded by $ C
n^{-\nu} (t-s)^H F_{1,p}$. On the other hand, from Minkowski's
inequality it follows that the $L^p$-norm of the second integral is
less than or equal to
$ C n^{-\nu} (t-s) F_{1,p}$.
These estimates imply (\ref{11.a}) because $(t-s) \le(t-s)^H T^{1-H}$.
\end{pf*}

\begin{pf*}{Proof of (\ref{11.b})}
If $t-s \leq\frac{1}n$, we can write
\begin{eqnarray*}
\biggl\llVert\int_s^t F_u(B_u-B_{\eta(u ) }
) \,du \biggr\rrVert_p & \leq& \int_s^t
\bigl\llVert F_u(B_u-B_{\eta(u ) } ) \bigr\rrVert
_p \,du
\\
& \leq& C\sup_{t \in[0, T]} \|F_t\|_{q}
n^{-H} (t-s)
\\
&\leq& C\sup_{t \in[0, T]}\|F_t\|_{q}
n^{-1} (t-s)^H,
\end{eqnarray*}
where the first inequality follows from Minkowski's inequality and the
second one from H\"older's inequality.
Suppose that $t-s\geq\frac{1}n$. Applying Fubini's theorem for the Young
integral, we obtain
\[
\int_s^t F_u(B_u-B_{\eta(u ) }
) \,du = \int_{\eta(s ) }^t \biggl( \int
_{v}^{\varepsilon(v)} \mathbf{1}_{[s,t]}(u)F_u
\,du \biggr) \,dB_v.
\]
Applying (\ref{11.a}) with $\nu=0$ we obtain
%
%
\begin{eqnarray}
\biggl\llVert\int_{\eta(s ) }^t \biggl( \int
_{v}^{\varepsilon(v)} \mathbf{1}_{[s,t]}(u)F_u
\,du \biggr) \,dB_v \biggr\rrVert_p &\le& C \bigl(t -
\eta(s) \bigr)^H n^{-1} F_{1,p}
\nonumber
\\[-8pt]
\\[-8pt]
\nonumber
& \le& C (t - s)^H n^{-1} F_{1,p}.
\end{eqnarray}
This completes the proof of (\ref{11.b}).
\end{pf*}

%
\begin{lemma} \label{lem11.4}
Let $B=\{ B_t, t \in[0,T] \}$ be an $m$-dimensional fBm with Hurst
parameter $H>\frac{1}2$. Fix $p\ge\frac{1}H$.
Let $F=\{F_t, t \in[0, T]\}$ be a stochastic process such that $F_t
\in\mathbb{D}^{2,q}$, $t\in[0,T]$, for some $q>p$. For any $\rho>1$
we set
\[
F_{2,\rho} = \sup_{r,s,t \in[0,T]} \bigl( \|F_t
\|_\rho\vee\| D_sF_t \| _\rho\vee
\llVert{D}_{r} D_{s}F_{t} \rrVert
_\rho\bigr).
\]
Set also
\[
F_{* } = \sup_{r,s,t \in[0, T] } \bigl( \vert
F_{t}\vert\vee\vert D_{s}F_{t}\vert
\vee\vert{D}_{r} D_{s}F_{t} \vert\bigr).
\]
Then there exists a constant $C$ (independent of $F$) such that the
following holds for all $0\le s<t \le T$, $i,j=1, \ldots, m$:
%
%
\begin{eqnarray}
\Biggl\llVert\sum_{k= \lfloor{ns}/ T \rfloor}^{ \lfloor{nt}/
T \rfloor}
F_{t_k} \int_{t_k \vee s}^{t_{k+1}\wedge t} \int
_{t_k}^{u} \delta B^i_v
\delta{B}^j_u \Biggr\rrVert_{p} & \leq& C
\gamma_n^{-1} (t-s)^{{1}/2 }\|F_{ *}
\|_{q}, \label{11.c}
\\
\Biggl\llVert\sum_{k= \lfloor{ns}/ T \rfloor}^{ \lfloor{nt}/
T \rfloor}
F_{t_k} \int_{t_k \vee s}^{t_{k+1}\wedge t} \int
_{t_k}^{u} \delta B^i_v
\delta{B}^j_u \Biggr\rrVert_{p} & \leq& C
n^{-H} (t-s)^{ H } F_{2,q }. \label{11.d}
\end{eqnarray}
\end{lemma}

\begin{pf}
Using $(\ref{e1})$, we can write
%
%
\begin{eqnarray}
\label{11.20e} 
&&\sum_{k = \lfloor{ns} /T \rfloor}^{ \lfloor{nt}/ T
\rfloor}
F_{t_k} \int_{t_k \vee s}^{t_{k+1}\wedge t} \int
_{t_k}^{u} \delta B^i_v
\delta{B}^j_u \nonumber\\
&&\qquad= \int_{ s}^{ t}
F_{\eta(u) } \bigl( B^i_u - B^i_{ \eta(u) }
\bigr) \delta{B}^j_u
\\
&&\qquad\quad{} + \alpha_H \int_s^t \int
_0^T { D}^j_r
F_{\eta(u)} \bigl( B^i_u - B^i_{ \eta(u) }
\bigr) \mu(drdu).
\nonumber
\end{eqnarray}
Applying (\ref{11.b}) to the second integral of the right-hand side of
(\ref{11.20e}) with $F_u $ replaced by $\int_0^T {D}^j_r F_{\eta(u)}
|r-u|^{2H-2} \,dr$ (notice that here we do not need the H\"older
continuity of the integrand for the Young integral to be well defined) yields
%
%
\begin{eqnarray}
\label{11.18e} && \biggl\llVert\int_s^t \int
_0^T {D}_r^j
F_{\eta(u)} \bigl( B^i_u - B^i_{ \eta(u) }
\bigr) \mu(dr\,du) \biggr\rrVert_p
\nonumber\\
&&\qquad\leq C n^{-1} (t-s)^H F_{2,q} \sup
_{u\in[0, T] } \int_0^T
|r-u|^{2H-2} \,dr
\\
&&\qquad\leq C n^{-1} (t-s)^H F_{2,q}.
\nonumber
\end{eqnarray}
This implies both estimates (\ref{11.c}) and (\ref{11.d}).

Applying (\ref{e1}) to the first summand on the right-hand side of
(\ref{11.20e}) yields
%
%
\begin{eqnarray}
\label{11.18 e} 
&& \int_s^t
F_u \bigl(B^i_u - B^i_{\eta(u) }
\bigr) \delta{B}^j_u
\nonumber
\\[-8pt]
\\[-8pt]
\nonumber
&&\qquad = \int_s^t \int_{\eta(u)}^u
F_u \delta B^i_v \delta{B}^j_u
+ \al_H \int_s^t \biggl\{ \int
_0^T \int_{\eta(u) }^u
D^i_v F_u \mu(drdv) \biggr\}
\delta{B}^j_u.
\end{eqnarray}
Now we apply (\ref{eqn 2.10}) to the second term of the right-hand side
of (\ref{11.18 e}), and we obtain
%
%
\begin{eqnarray}
\label{11.19e} && \biggl\llVert\int_s^t
\biggl\{ \int_0^T \int_{\eta(u) }^u
D^i_v F_u \mu(dr\,dv) \biggr\}
\delta{B}^j_u \biggr\rrVert_p\nonumber
\\
& &\qquad \leq C F_{2, p} \biggl\llVert\mathbf{1}_{[s, t ] } (u) \int
_0^T \int_{\eta(u) }^u
\mu(dr\,dv) \biggr\rrVert_{L^{{1}/H} ( [0, T] ) }
\\
&&\qquad\leq C F_{2, p} n^{-1}(t-s)^H.
\nonumber
\end{eqnarray}
Again, this inequality implies both estimates (\ref{11.c}) and (\ref{11.d}).

It remains to estimate the term $
I_{s,t}:= \int_s^t \int_{\eta(u)}^u F_u \delta B^i_v \delta{B}^j_u$.
It follows from (\ref{eqn 2.10}) that
\begin{eqnarray*}
\| I_{s,t} \|_p &\le& C F_{2, p} \bigl\|
\mathbf{1}_{[s,t]}(u) \mathbf{1}_{[\eta(u), u]} (v) \bigr\|_{L^{{1}/H}
([0,T]^2 )}
\\
&\le& C F_{2, p} n^{-H} (t-s)^{H},
\end{eqnarray*}
which completes the proof of (\ref{11.d}).

To derive (\ref{11.c}) we need a more accurate estimate.

Meyer's inequality implies that
\begin{eqnarray*}
\| I_{s,t}\| _{p } &\leq& C \bigl[ \bigl\llVert\bigl
\llVert\mathbf{1}_{[s, t]}(u) \mathbf{1}_{[\eta(u), u]}(v)
F_u \bigr\rrVert_{\cH
^{\otimes2}} \bigr\rrVert_p
\\
&&{}+ \bigl\llVert\bigl\llVert\mathbf{1}_{[s, t]}(u)
\mathbf{1}_{[\eta(u), u]}(v) D_{r}F_u \bigr\rrVert
_{\cH^{\otimes3}} \bigr\rrVert_p
\\
&&{}+ \bigl\llVert\bigl\llVert\mathbf{1}_{[s, t]}(u)
\mathbf{1}_{[\eta(u), u]}(v) D_{r'}D_{r}F_u
\bigr\rrVert_{\cH^{\otimes4}}\bigr\rrVert_p \bigr]
\\
&\leq& C \|F_{* } \|_{p} \bigl\llVert
\mathbf{1}_{[s, t]}(u) \mathbf{1}_{[\eta(u),
u]}(v) \bigr\rrVert
_{\cH^{\otimes2}}.
\end{eqnarray*}
Therefore, to complete the proof, it suffices to show that
%
%
\begin{eqnarray}
\label{e3.1} 
&& \bigl\| \mathbf{1}_{[s, t]}(u)
\mathbf{1}_{[\eta(u), u]}(v) \bigr\|_{\cH
^{\otimes2}}^2\nonumber
\\
& &\qquad= \alpha_H^2 \int_s^t
\int_s^t \int_{\eta(u')}^{ u'}
\int_{\eta(u)}^{ u} \mu\bigl(dv\,dv'
\bigr) \mu\bigl(du\,du'\bigr)
\\
& &\qquad \leq(t-s)\gamma_n^{-2}.
\nonumber
\end{eqnarray}
In the case $t-s \geq\frac{1}n$,
\begin{eqnarray*}
&& \int_s^t \int_s^t
\int_{\eta(u')}^{ u'} \int_{\eta(u)}^{ u}
\mu\bigl(dv\,dv'\, du\,du'\bigr)
\\
&&\qquad \leq\sum_{k,k' = \lfloor{ns}/ T \rfloor}^{ \lfloor{nt}/ T
\rfloor} \int
_{t_{k' } }^{t_{{k' } +1} } \int_{t_k}^{t_{k+1} }
\int_{t_{k' } }^{ u'} \int_{t_k}^{ u}
\mu\bigl(dv\,dv'\, du\,du'\bigr)
\\
& &\qquad\leq\sum_{k = \lfloor{ns}/ T \rfloor}^{ \lfloor{nt}/ T
\rfloor} \sum
_{p= 1-n}^{n-1} \int_{t_{k+p } }^{t_{{k+p } +1} }
\int_{t_k}^{t_{k+1} } \int_{t_{k+p }
}^{ u'}
\int_{t_k}^{ u} \mu\bigl(dv\,dv'\,
du\,du'\bigr)
\\
& &\qquad= n^{-4H} \sum_{k = \lfloor{ns}/ T \rfloor}^{ \lfloor{nt}/ T
\rfloor}
\sum_{p= 1-n}^{n-1} Q(p )
\\
&&\qquad \leq C(t-s)\gamma_n^{-2},
\end{eqnarray*}
where we recall that $Q(p)$ is defined in Section~\ref{sec2.3}, and
inequality (\ref{e3.1}) follows.

In the case $t-s \leq\frac{1}{n}$, we have the raw estimate
\[
\int_s^t \int_s^t
\int_{\eta(u')}^{ u'}\int_{\eta(u)}^{ u}
\mu\bigl(dr\,dr'\, du\,du'\bigr) \leq\frac{1}{n^{2H}}
\int_s^t \int_s^t
\mu\bigl( du\,du'\bigr) = \frac
{1}{n^{2H}}(t-s)^{2H}
\]
and
\[
{n^{-2H}}(t-s)^{2H} \leq(t-s) \gamma_n^{-2}.
\]
So (\ref{e3.1}) is also true for this case. The proof of the lemma is
now complete.
\end{pf}
%

\begin{lemma} \label{lem11.5}
Let $B=\{B_t, t
\in[0,T]\}$ be a one-dimensional fBm with Hurst parameter $H>\frac
{1}2$. Suppose that $F=\{F_t, t\in[0,T\} $, $G=\{G_t, t\in[0,T]\} $ are
processes that are H\"older continuous of order $\beta\in(\frac
{1}2,H) $. Then there exists a constant $C$ (not depending on $F$ or $G$)
such that for all $0\le s<t \le T$, $\nu\ge0$,
%
%
\begin{eqnarray}\label{e 11.10}
&& \biggl\vert\int_s^t F_u
(G_u- G_{\eta(u)}) (B_u - B_{\eta(u)} )
\,dB_u \biggr\vert
\nonumber
\\[-8pt]
\\[-8pt]
\nonumber
&& \qquad\leq C\bigl( \|F \|_{\infty} + \|F\|_{ \beta} \bigr) \|G\|_{\beta}
\|B\|_{\beta}^2 n^{1-3\beta} (t-s)^{\beta}
\end{eqnarray}
and
%
%
\begin{eqnarray}\label{e 11.11}
&& \biggl\vert\int_s^t F_u
(G_u- G_{\eta(u)}) \bigl(u-{\eta(u)} \bigr)^{\nu}
\,dB_u \biggr\vert
\nonumber
\\[-8pt]
\\[-8pt]
\nonumber
&&\qquad \leq C\bigl( \|F \|_{\infty} + \|F\|_{ \beta} \bigr) \|G\|_{\beta}
\|B\|_{\beta} n^{1-2\beta-\nu} (t-s)^{\beta}.
\nonumber
\end{eqnarray}
\end{lemma}
\begin{pf*}{Proof of (\ref{e 11.10})}
We assume first that $s, t \in[t_k, t_{k+1} ]$ for some $k = 0, 1,
\ldots, n-1$. By Lemma~\ref{lem7.1}(ii),
%
%
\begin{eqnarray}
\label{e11.7} && \biggl\vert\int_s^t
F_u (G_u- G_{\eta(u)}) (B_u-B_{\eta(u)})
\,dB_u \biggr\vert
\nonumber
\\
&&\qquad \leq K_1 \sup_{u \in[s, t] } \bigl| F_u
(G_u- G_{t_k}) (B_u-B_{ t_k}) \bigr| \|B
\| _{\beta} (t-s)^{\beta}
\nonumber
\\
&&\qquad\quad{} + K_2 \sup_{u \in[s, t] } \bigl[ \bigl\vert
F_u (G_u- G_{ t_k }) \bigr\vert\|B
\|_{\beta}^2 (t-s)^{ 2 \beta}
\nonumber
\\[-8pt]
\\[-8pt]
\nonumber
&&\hspace*{57pt}\qquad\quad{} + \bigl\vert F_u (B_u-B_{t_k }) \bigr
\vert\|G\|_{\beta} \|B\|_{\beta} (t-s)^{2\beta}
\\
&&\hspace*{57pt}\qquad\quad{} + \bigl\vert(G_u- G_{t_k}) (B_u-B_{ t_k})
\bigr\vert\|F\|_{\beta} \|B\| _{\beta} (t-s)^{2\beta}
\bigr]
\nonumber
\\
&&\qquad\leq C \kappa_\beta(F,G) n^{-2\beta} (t-s)^{\beta},\nonumber
\end{eqnarray}
where $ \kappa_\beta(F,G)=( \|F \|_{\infty} + \| F\|_{ \beta} ) \|
G\|
_{\beta} \|B\|_{\beta}^2$.
In the general case, we can write
\begin{eqnarray*}
&& \biggl\vert\int_s^t F_u
(G_u- G_{\eta(u)}) (B_u-B_{\eta(u)})
\,dB_u \biggr\vert
\nonumber
\\
&&\qquad= \Biggl\vert\Biggl(\int_s^{\ep(s ) } + \sum
_{k = \lfloor{ns}/{T}
\rfloor+1}^{ \lfloor{nt}/{T} \rfloor} \int_{t_k}^{t_{k+1 } }+
\int_{\eta(t)}^{t} \Biggr) F_u
(G_u- G_{\eta
(u)}) (B_u-B_{\eta(u)})
\,dB_u \Biggr\vert
\nonumber
\\
&&\qquad\leq C \kappa_\beta(F,G) n^{-2 \beta} \Biggl[ \bigl(\ep(s) - s
\bigr)^{\beta
}+\bigl(t-\eta(t) \bigr)^{\beta} + \sum
_{k=\lfloor{ns}/{T} \rfloor+1}^{ \lfloor{nt}/{T}
\rfloor
} (T/n)^{\beta} \Biggr]
\nonumber
\\
&&\qquad\leq C \kappa_\beta(F,G) n^{-2\beta} \bigl[ \bigl(\ep(s) - s
\bigr)^{\beta} +\bigl(t-\eta(t) \bigr)^{\beta} + \bigl(\eta(t) -
\ep(s)\bigr) n^{1-\beta} \bigr]
\nonumber
\\
&&\qquad\leq C \kappa_\beta(F,G)n^{1-3 \beta} (t-s)^{\beta},
\nonumber
\end{eqnarray*}
where the first inequality follows from (\ref{e11.7}).
\end{pf*}

\begin{pf*}{Proof of (\ref{e 11.11})} This estimate can be
proved by following the lines of the proof of (\ref{e 11.10}) and
noticing the fact that $(u-\eta(u ) )^{\nu}$ has finite $\nu$-H\"older
seminorm on $(t_k, t_{k+1})$ for each $k=1, \ldots, n-1$.
\end{pf*}
\noqed\end{pf}
\end{appendix}


\section*{Acknowledgment}
We wish to thank the referee for many useful comments and suggestions.
%
%





\printaddresses

\begin{thebibliography}{34}

\bibitem{Aldous}
%
\begin{barticle}[mr]
\bauthor{\bsnm{Aldous},~\bfnm{D.~J.}\binits{D.~J.}} \AND
\bauthor{\bsnm{Eagleson},~\bfnm{G.~K.}\binits{G.~K.}}
(\byear{1978}).
\btitle{On mixing and stability of limit theorems}.
\bjournal{Ann. Probab.}
\bvolume{6}
\bpages{325--331}.
\bid{mr={0517416}}
\end{barticle}
%

\bptok{imsref}%
\endbibitem

\bibitem{CH}
%
\begin{barticle}[mr]
\bauthor{\bsnm{Cambanis},~\bfnm{Stamatis}\binits{S.}} \AND
\bauthor{\bsnm{Hu},~\bfnm{Yaozhong}\binits{Y.}}
(\byear{1996}).
\btitle{Exact convergence rate of the {E}uler--{M}aruyama scheme, with
application to sampling design}.
\bjournal{Stochastics Stochastics Rep.}
\bvolume{59}
\bpages{211--240}.
\bid{issn={1045-1129}, mr={1427739}}
\end{barticle}
%

\bptok{imsref}%
\endbibitem

\bibitem{CNP}
%
\begin{barticle}[mr]
\bauthor{\bsnm{Corcuera},~\bfnm{Jos{\'e}~Manuel}\binits{J.~M.}},
\bauthor{\bsnm{Nualart},~\bfnm{David}\binits{D.}} \AND
\bauthor{\bsnm{Podolskij},~\bfnm{Mark}\binits{M.}}
(\byear{2014}).
\btitle{Asymptotics of weighted random sums}.
\bjournal{Commun. Appl. Ind. Math.}
\bvolume{6}
\bpages{e--486, 11}.
\bid{doi={10.1685/journal.caim.486}, issn={2038-0909}, mr={3277312}}
\bptnote{check volume, check pages, check year}%
\end{barticle}
%

\bptok{imsref}%
\endbibitem

\bibitem{DNT}
%
\begin{barticle}[mr]
\bauthor{\bsnm{Deya},~\bfnm{A.}\binits{A.}},
\bauthor{\bsnm{Neuenkirch},~\bfnm{A.}\binits{A.}} \AND
\bauthor{\bsnm{Tindel},~\bfnm{S.}\binits{S.}}
(\byear{2012}).
\btitle{A {M}ilstein-type scheme without L\'evy area terms for {SDE}s
driven by fractional {B}rownian motion}.
\bjournal{Ann. Inst. Henri Poincar\'e Probab. Stat.}
\bvolume{48}
\bpages{518--550}.
\bid{doi={10.1214/10-AIHP392}, issn={0246-0203}, mr={2954265}}
\end{barticle}
%

\bptok{imsref}%
\endbibitem

\bibitem{FR}
%
\begin{barticle}[mr]
\bauthor{\bsnm{Friz},~\bfnm{Peter}\binits{P.}} \AND
\bauthor{\bsnm{Riedel},~\bfnm{Sebastian}\binits{S.}}
(\byear{2014}).
\btitle{Convergence rates for the full {G}aussian rough paths}.
\bjournal{Ann. Inst. Henri Poincar\'e Probab. Stat.}
\bvolume{50}
\bpages{154--194}.
\bid{doi={10.1214/12-AIHP507}, issn={0246-0203}, mr={3161527}}
\end{barticle}
%

\bptok{imsref}%
\endbibitem

\bibitem{FV}
%
\begin{bbook}[mr]
\bauthor{\bsnm{Friz},~\bfnm{Peter~K.}\binits{P.~K.}} \AND
\bauthor{\bsnm{Victoir},~\bfnm{Nicolas~B.}\binits{N.~B.}}
(\byear{2010}).
\btitle{Multidimensional Stochastic Processes as Rough Paths: Theory
and Applications}.
\bseries{Cambridge Studies in Advanced Mathematics}
\bvolume{120}.
\bpublisher{Cambridge Univ. Press},
\blocation{Cambridge}.
\bid{doi={10.1017/CBO9780511845079}, mr={2604669}}
\end{bbook}
%

\bptok{imsref}%
\endbibitem

\bibitem{Hu}
%
\begin{bincollection}[mr]
\bauthor{\bsnm{Hu},~\bfnm{Yaozhong}\binits{Y.}}
(\byear{1996}).
\btitle{Strong and weak order of time discretization schemes of
stochastic differential equations}.
In \bbooktitle{S\'eminaire de {P}robabilit\'es, {XXX}}.
\bseries{Lecture Notes in Math.}
\bvolume{1626}
\bpages{218--227}.
\bpublisher{Springer},
\blocation{Berlin}.
\bid{doi={10.1007/BFb0094650}, mr={1459485}}
\end{bincollection}
%

\bptok{imsref}%
\endbibitem

\bibitem{HuNu}
%
\begin{bincollection}[mr]
\bauthor{\bsnm{Hu},~\bfnm{Yaozhong}\binits{Y.}} \AND
\bauthor{\bsnm{Nualart},~\bfnm{David}\binits{D.}}
(\byear{2007}).
\btitle{Differential equations driven by H\"older continuous functions
of order greater than $1/2$}.
In \bbooktitle{Stochastic Analysis and Applications}.
\bseries{Abel Symp.}
\bvolume{2}
\bpages{399--413}.
\bpublisher{Springer},
\blocation{Berlin}.
\bid{doi={10.1007/978-3-540-70847-6_17}, mr={2397797}}
\bptnote{check pages}%
\end{bincollection}
%

\bptok{imsref}%
\endbibitem

\bibitem{hunu09}
%
\begin{barticle}[mr]
\bauthor{\bsnm{Hu},~\bfnm{Yaozhong}\binits{Y.}} \AND
\bauthor{\bsnm{Nualart},~\bfnm{David}\binits{D.}}
(\byear{2009}).
\btitle{Stochastic heat equation driven by fractional noise and local time}.
\bjournal{Probab. Theory Related Fields}
\bvolume{143}
\bpages{285--328}.
\bid{doi={10.1007/s00440-007-0127-5}, issn={0178-8051}, mr={2449130}}
\end{barticle}
%

\bptok{imsref}%
\endbibitem

\bibitem{Jacod}
%
\begin{barticle}[mr]
\bauthor{\bsnm{Jacod},~\bfnm{Jean}\binits{J.}} \AND
\bauthor{\bsnm{Protter},~\bfnm{Philip}\binits{P.}}
(\byear{1998}).
\btitle{Asymptotic error distributions for the {E}uler method for
stochastic differential equations}.
\bjournal{Ann. Probab.}
\bvolume{26}
\bpages{267--307}.
\bid{doi={10.1214/aop/1022855419}, issn={0091-1798}, mr={1617049}}
\end{barticle}
%

\bptok{imsref}%
\endbibitem

\bibitem{JacodShiryaev}
%
\begin{bbook}[mr]
\bauthor{\bsnm{Jacod},~\bfnm{Jean}\binits{J.}} \AND
\bauthor{\bsnm{Shiryaev},~\bfnm{Albert~N.}\binits{A.~N.}}
(\byear{1987}).
\btitle{Limit Theorems for Stochastic Processes}.
\bseries{Grundlehren der Mathematischen Wissenschaften}
\bvolume{288}.
\bpublisher{Springer},
\blocation{Berlin}.
\bid{doi={10.1007/978-3-662-02514-7}, mr={0959133}}
\end{bbook}
%

\bptok{imsref}%
\endbibitem

\bibitem{Kloeden}
%
\begin{bbook}[mr]
\bauthor{\bsnm{Kloeden},~\bfnm{Peter~E.}\binits{P.~E.}} \AND
\bauthor{\bsnm{Platen},~\bfnm{Eckhard}\binits{E.}}
(\byear{1992}).
\btitle{Numerical Solution of Stochastic Differential Equations}.
\bseries{Applications of Mathematics (New York)}
\bvolume{23}.
\bpublisher{Springer},
\blocation{Berlin}.
\bid{doi={10.1007/978-3-662-12616-5}, mr={1214374}}
\end{bbook}
%

\bptok{imsref}%
\endbibitem

\bibitem{Kurtz}
%
\begin{bincollection}[mr]
\bauthor{\bsnm{Kurtz},~\bfnm{Thomas~G.}\binits{T.~G.}} \AND
\bauthor{\bsnm{Protter},~\bfnm{Philip}\binits{P.}}
(\byear{1991}).
\btitle{Wong--{Z}akai corrections, random evolutions, and simulation
schemes for {SDE}s}.
In \bbooktitle{Stochastic Analysis}
\bpages{331--346}.
\bpublisher{Academic Press},
\blocation{Boston, MA}.
\bid{mr={1119837}}
\end{bincollection}
%

\bptok{imsref}%
\endbibitem

\bibitem{Ly}
%
\begin{barticle}[mr]
\bauthor{\bsnm{Lyons},~\bfnm{Terry}\binits{T.}}
(\byear{1994}).
\btitle{Differential equations driven by rough signals. I. {A}n
extension of an inequality of L. C. {Y}oung}.
\bjournal{Math. Res. Lett.}
\bvolume{1}
\bpages{451--464}.
\bid{doi={10.4310/MRL.1994.v1.n4.a5}, issn={1073-2780}, mr={1302388}}
\end{barticle}
%

\bptok{imsref}%
\endbibitem

\bibitem{LQ}
%
\begin{bbook}[mr]
\bauthor{\bsnm{Lyons},~\bfnm{Terry}\binits{T.}} \AND
\bauthor{\bsnm{Qian},~\bfnm{Zhongmin}\binits{Z.}}
(\byear{2002}).
\btitle{System Control and Rough Paths}.
\bpublisher{Oxford Univ. Press},
\blocation{Oxford}.
\bid{doi={10.1093/acprof:oso/9780198506485.001.0001}, mr={2036784}}
\end{bbook}
%

\bptok{imsref}%
\endbibitem

\bibitem{MMV}
%
\begin{barticle}[mr]
\bauthor{\bsnm{M{\'e}min},~\bfnm{Jean}\binits{J.}},
\bauthor{\bsnm{Mishura},~\bfnm{Yulia}\binits{Y.}} \AND
\bauthor{\bsnm{Valkeila},~\bfnm{Esko}\binits{E.}}
(\byear{2001}).
\btitle{Inequalities for the moments of {W}iener integrals with
respect to a fractional {B}rownian motion}.
\bjournal{Statist. Probab. Lett.}
\bvolume{51}
\bpages{197--206}.
\bid{doi={10.1016/S0167-7152(00)00157-7}, issn={0167-7152}, mr={1822771}}
\end{barticle}
%

\bptok{imsref}%
\endbibitem

\bibitem{Mishura}
%
\begin{bbook}[mr]
\bauthor{\bsnm{Mishura},~\bfnm{Yuliya~S.}\binits{Y.~S.}}
(\byear{2008}).
\btitle{Stochastic Calculus for Fractional {B}rownian Motion and
Related Processes}.
\bseries{Lecture Notes in Math.}
\bvolume{1929}.
\bpublisher{Springer},
\blocation{Berlin}.
\bid{doi={10.1007/978-3-540-75873-0}, mr={2378138}}
\end{bbook}
%

\bptok{imsref}%
\endbibitem

\bibitem{Neuenkirch}
%
\begin{barticle}[mr]
\bauthor{\bsnm{Neuenkirch},~\bfnm{Andreas}\binits{A.}} \AND
\bauthor{\bsnm{Nourdin},~\bfnm{Ivan}\binits{I.}}
(\byear{2007}).
\btitle{Exact rate of convergence of some approximation schemes
associated to {SDE}s driven by a fractional {B}rownian motion}.
\bjournal{J. Theoret. Probab.}
\bvolume{20}
\bpages{871--899}.
\bid{doi={10.1007/s10959-007-0083-0}, issn={0894-9840}, mr={2359060}}
\end{barticle}
%

\bptok{imsref}%
\endbibitem

\bibitem{NTU}
%
\begin{barticle}[mr]
\bauthor{\bsnm{Neuenkirch},~\bfnm{A.}\binits{A.}},
\bauthor{\bsnm{Tindel},~\bfnm{S.}\binits{S.}} \AND
\bauthor{\bsnm{Unterberger},~\bfnm{J.}\binits{J.}}
(\byear{2010}).
\btitle{Discretizing the fractional L\'evy area}.
\bjournal{Stochastic Process. Appl.}
\bvolume{120}
\bpages{223--254}.
\bid{doi={10.1016/j.spa.2009.10.007}, issn={0304-4149}, mr={2576888}}
\end{barticle}
%

\bptok{imsref}%
\endbibitem

\bibitem{NNT}
%
\begin{barticle}[mr]
\bauthor{\bsnm{Nourdin},~\bfnm{Ivan}\binits{I.}},
\bauthor{\bsnm{Nualart},~\bfnm{David}\binits{D.}} \AND
\bauthor{\bsnm{Tudor},~\bfnm{Ciprian~A.}\binits{C.~A.}}
(\byear{2010}).
\btitle{Central and non-central limit theorems for weighted power
variations of fractional {B}rownian motion}.
\bjournal{Ann. Inst. Henri Poincar\'e Probab. Stat.}
\bvolume{46}
\bpages{1055--1079}.
\bid{doi={10.1214/09-AIHP342}, issn={0246-0203}, mr={2744886}}
\end{barticle}
%

\bptok{imsref}%
\endbibitem

\bibitem{Peccati2}
%
\begin{bbook}[mr]
\bauthor{\bsnm{Nourdin},~\bfnm{Ivan}\binits{I.}} \AND
\bauthor{\bsnm{Peccati},~\bfnm{Giovanni}\binits{G.}}
(\byear{2012}).
\btitle{Normal Approximations with {M}alliavin Calculus: From Stein's
Method to Universality}.
\bseries{Cambridge Tracts in Mathematics}
\bvolume{192}.
\bpublisher{Cambridge Univ. Press},
\blocation{Cambridge}.
\bid{doi={10.1017/CBO9781139084659}, mr={2962301}}
\end{bbook}
%

\bptok{imsref}%
\endbibitem

\bibitem{Nualart2}
%
\begin{bincollection}[mr]
\bauthor{\bsnm{Nualart},~\bfnm{David}\binits{D.}}
(\byear{2003}).
\btitle{Stochastic integration with respect to fractional {B}rownian
motion and applications}.
In \bbooktitle{Stochastic Models ({M}exico {C}ity, 2002)}.
\bseries{Contemp. Math.}
\bvolume{336}
\bpages{3--39}.
\bpublisher{Amer. Math. Soc.},
\blocation{Providence, RI}.
\bid{doi={10.1090/conm/336/06025}, mr={2037156}}
\end{bincollection}
%

\bptok{imsref}%
\endbibitem

\bibitem{Nualart}
%
\begin{bbook}[mr]
\bauthor{\bsnm{Nualart},~\bfnm{David}\binits{D.}}
(\byear{2006}).
\btitle{The {M}alliavin Calculus and Related Topics},
\bedition{2nd} ed.
\bpublisher{Springer},
\blocation{Berlin}.
\bid{mr={2200233}}
\end{bbook}
%

\bptok{imsref}%
\endbibitem

\bibitem{Nualart3}
%
\begin{barticle}[mr]
\bauthor{\bsnm{Nualart},~\bfnm{David}\binits{D.}} \AND
\bauthor{\bsnm{Peccati},~\bfnm{Giovanni}\binits{G.}}
(\byear{2005}).
\btitle{Central limit theorems for sequences of multiple stochastic integrals}.
\bjournal{Ann. Probab.}
\bvolume{33}
\bpages{177--193}.
\bid{doi={10.1214/009117904000000621}, issn={0091-1798}, mr={2118863}}
\end{barticle}
%

\bptok{imsref}%
\endbibitem

\bibitem{NuRa}
%
\begin{barticle}[mr]
\bauthor{\bsnm{Nualart},~\bfnm{David}\binits{D.}} \AND
\bauthor{\bsnm{R{\u{a}}{\c{s}}canu},~\bfnm{Aurel}\binits{A.}}
(\byear{2002}).
\btitle{Differential equations driven by fractional {B}rownian motion}.
\bjournal{Collect. Math.}
\bvolume{53}
\bpages{55--81}.
\bid{issn={0010-0757}, mr={1893308}}
\end{barticle}
%

\bptok{imsref}%
\endbibitem

\bibitem{NuSau}
%
\begin{barticle}[mr]
\bauthor{\bsnm{Nualart},~\bfnm{David}\binits{D.}} \AND
\bauthor{\bsnm{Saussereau},~\bfnm{Bruno}\binits{B.}}
(\byear{2009}).
\btitle{Malliavin calculus for stochastic differential equations
driven by a fractional {B}rownian motion}.
\bjournal{Stochastic Process. Appl.}
\bvolume{119}
\bpages{391--409}.
\bid{doi={10.1016/j.spa.2008.02.016}, issn={0304-4149}, mr={2493996}}
\end{barticle}
%

\bptok{imsref}%
\endbibitem

\bibitem{Peccati}
%
\begin{bincollection}[mr]
\bauthor{\bsnm{Peccati},~\bfnm{Giovanni}\binits{G.}} \AND
\bauthor{\bsnm{Tudor},~\bfnm{Ciprian~A.}\binits{C.~A.}}
(\byear{2005}).
\btitle{Gaussian limits for vector-valued multiple stochastic integrals}.
In \bbooktitle{S\'eminaire de {P}robabilit\'es {XXXVIII}}.
\bseries{Lecture Notes in Math.}
\bvolume{1857}
\bpages{247--262}.
\bpublisher{Springer},
\blocation{Berlin}.
\bid{mr={2126978}}
\end{bincollection}
%

\bptok{imsref}%
\endbibitem

\bibitem{PT1}
%
\begin{barticle}[mr]
\bauthor{\bsnm{Pipiras},~\bfnm{Vladas}\binits{V.}} \AND
\bauthor{\bsnm{Taqqu},~\bfnm{Murad~S.}\binits{M.~S.}}
(\byear{2000}).
\btitle{Integration questions related to fractional {B}rownian motion}.
\bjournal{Probab. Theory Related Fields}
\bvolume{118}
\bpages{251--291}.
\bid{doi={10.1007/s440-000-8016-7}, issn={0178-8051}, mr={1790083}}
\end{barticle}
%

\bptok{imsref}%
\endbibitem

\bibitem{PT2}
%
\begin{barticle}[mr]
\bauthor{\bsnm{Pipiras},~\bfnm{Vladas}\binits{V.}} \AND
\bauthor{\bsnm{Taqqu},~\bfnm{Murad~S.}\binits{M.~S.}}
(\byear{2001}).
\btitle{Are classes of deterministic integrands for fractional
{B}rownian motion on an interval complete?}
\bjournal{Bernoulli}
\bvolume{7}
\bpages{873--897}.
\bid{doi={10.2307/3318624}, issn={1350-7265}, mr={1873833}}
\end{barticle}
%

\bptok{imsref}%
\endbibitem

\bibitem{Rosenblatt}
%
\begin{bbook}[mr]
\bauthor{\bsnm{Rosenblatt},~\bfnm{Murray}\binits{M.}}
(\byear{2011}).
\btitle{Selected Works of {M}urray {R}osenblatt}.
\bpublisher{Springer},
\blocation{New York}.
\bid{mr={2742596}}
\end{bbook}
%

\bptok{imsref}%
\endbibitem

\bibitem{TT}
%
\begin{barticle}[mr]
\bauthor{\bsnm{Talay},~\bfnm{Denis}\binits{D.}} \AND
\bauthor{\bsnm{Tubaro},~\bfnm{Luciano}\binits{L.}}
(\byear{1990}).
\btitle{Expansion of the global error for numerical schemes solving
stochastic differential equations}.
\bjournal{Stoch. Anal. Appl.}
\bvolume{8}
\bpages{483--509}.
\bid{doi={10.1080/07362999008809220}, issn={0736-2994}, mr={1091544}}
\bptnote{check pages}%
\end{barticle}
%

\bptok{imsref}%
\endbibitem

\bibitem{Tudor}
%
\begin{barticle}[mr]
\bauthor{\bsnm{Tudor},~\bfnm{Ciprian~A.}\binits{C.~A.}}
(\byear{2008}).
\btitle{Analysis of the {R}osenblatt process}.
\bjournal{ESAIM Probab. Stat.}
\bvolume{12}
\bpages{230--257}.
\bid{doi={10.1051/ps:2007037}, issn={1292-8100}, mr={2374640}}
\end{barticle}
%

\bptok{imsref}%
\endbibitem

\bibitem{Yo}
%
\begin{barticle}[mr]
\bauthor{\bsnm{Young},~\bfnm{L.~C.}\binits{L.~C.}}
(\byear{1936}).
\btitle{An inequality of the H\"older type, connected with {S}tieltjes
integration}.
\bjournal{Acta Math.}
\bvolume{67}
\bpages{251--282}.
\bid{doi={10.1007/BF02401743}, issn={0001-5962}, mr={1555421}}
\end{barticle}
%

\bptok{imsref}%
\endbibitem

\bibitem{Za}
%
\begin{barticle}[mr]
\bauthor{\bsnm{Z{\"a}hle},~\bfnm{M.}\binits{M.}}
(\byear{1998}).
\btitle{Integration with respect to fractal functions and stochastic
calculus. {I}}.
\bjournal{Probab. Theory Related Fields}
\bvolume{111}
\bpages{333--374}.
\bid{doi={10.1007/s004400050171}, issn={0178-8051}, mr={1640795}}
\end{barticle}
%

\bptok{imsref}%
\endbibitem
\end{thebibliography}
\end{document}